\documentclass[10pt]{amsart}
\usepackage[margin=2.5cm]{geometry}
\usepackage{tikz-cd}
\usepackage{graphicx}					
\usepackage{amssymb}
\usepackage{amsmath}
\usepackage{subfig, graphicx}
\usepackage{mathrsfs}
\usepackage{amssymb, amsthm, indentfirst}
\usepackage{relsize}
\usepackage{hyperref}
\usepackage{stmaryrd}
\usepackage{lineno}

\numberwithin{equation}{section}
\usepackage{color}
\theoremstyle{plain}
\newtheorem{thm}{Theorem}[section]

\newtheorem{lemma}[thm]{Lemma}
\newtheorem{prop}[thm]{Proposition}
\newtheorem{corollary}[thm]{Corollary}
\theoremstyle{definition}

\newtheorem{rmk}[thm]{Remark}

\newtheorem*{hypothesis*}{Hypothesis}

\author{SHIH-YU CHEN}
\address{Institute of Mathematics~\\Academia Sinica~\\ 6F, Astronomy-Mathematics Building, No.\,1, Sec.\,4, Roosevelt Road, Taipei 10617, Taiwan, ROC}
\email{sychen0626@gate.sinica.edu.tw}	
	
\def\SL{{\rm{SL}}}
\def\GL{{\rm{GL}}}
\def\GSp{{\rm GSp}}
\def\PGSp{{\rm PGSp}}
\def\Sp{{\rm Sp}}

\def\o{\frak{o}}
\def\c{\frak{c}}
\def\A{{\mathbb A}}
\def\C{{\mathbb C}}
\def\F{{\mathbb F}}
\def\R{{\mathbb R}}
\def\Q{{\mathbb Q}}
\def\Z{{\mathbb Z}}

\def\<{\langle}
\def\>{\rangle}

\def\x{\times}

\def\bp{\begin{pmatrix}}
\def\ep{\end{pmatrix}}

\def\<{\langle}
\def\>{\rangle}

\def\Ad{\operatorname{Ad}}
\def\diag{\operatorname{diag}}
\def\GL{\operatorname{GL}}
\def\GSO{\operatorname{GSO}}
\def\GSp{\operatorname{GSp}}

\def\O{\operatorname{O}}

\def\SL{\operatorname{SL}}
\def\SO{\operatorname{SO}}
\def\Sp{\operatorname{Sp}}

\def\Sym{\operatorname{Sym}}
\def\tr{\operatorname{tr}}

\def\AA{\mathbb{A}}

\def\GG{\mathbb{G}}

\def\calF{\mathcal{F}}

\def\calS{\mathcal{S}}
\def\calZ{\mathcal{Z}}

\def\1{\mathbf{1}}

\def\itPi{\mathit{\Pi}}
\def\itPsi{\mathit{\Psi}}

\title{Algebraicity of critical values of adjoint $L$-functions for $\GSp_4$}

\begin{document}
\begin{abstract}
We prove an algebraicity result for certain critical value of adjoint $L$-functions for $\GSp_4$ over a totally real number field in terms of the Petersson norm of normalized generic cuspidal newforms on $\GSp_4$.
This is a generalization of our previous result \cite{CI2019}.
\end{abstract}

\maketitle

\section{Introduction}

\subsection{Main result}

Let $f$ be a normalized elliptic newform of weight $\kappa \geq 2$ and level $\Gamma_1(N)$. Denote by $L(s,f,{\rm Ad})$ the completed adjoint $L$-function of $f$.
By the result of Sturm \cite{Sturm1989}, the algebraicity of 
$L(1,f,{\rm Ad})$ is expressed in terms of the Petersson norm
\[
 \Vert f \Vert = \int_{\Gamma_0(N) \backslash \mathfrak{H}}
 |f(\tau)|^2 \operatorname{Im}(\tau)^{\kappa-2} \, d \tau.
\]
More precisely, we have
\[
\sigma\left( \frac{L(1,f,{\rm Ad})}{\Vert f \Vert}\right) = \frac{L(1,{}^\sigma\!f,{\rm Ad})}{\Vert {}^\sigma\!f \Vert}
\]
for all $\sigma \in {\rm Aut}(\C)$, as predicted by Deligne's conjecture \cite{Deligne1979}.
The purpose of this paper is to prove an analogue of it for $\GSp_4$.

We give a description of our main result. Let $\itPi=\bigotimes_v \itPi_v$ be an irreducible globally generic cuspidal automorphic representation of $\GSp_4(\A_\F)$ with central character $\omega_\itPi$ over a totally real number field $\F$. 
Denote by 
\[
L(s,\itPi,{\rm Ad}) = \prod_{v}L(s,\itPi_v,{\rm Ad})
\]
the adjoint $L$-function of $\itPi$, where $v$ runs through the places of $\F$.
Assume $\itPi_\infty = \bigotimes_{v \mid \infty}\itPi_v$ is a discrete series representation of $\GSp_4(\F_\infty)$. Then
\[
\itPi_{\infty}\vert_{\Sp_4(\F_\infty)} = \bigoplus_{v \mid \infty}D_{(\lambda_{1,v},\lambda_{2,v})}\oplus D_{(-\lambda_{2,v},-\lambda_{1,v})},
\]
where $D_{(\lambda_{1,v},\lambda_{2,v})}$ is the discrete series representation of $\Sp_4(\F_v) \simeq \Sp_4(\R)$ with Blattner parameter $(\lambda_{1,v},\lambda_{2,v})$ such that $2-\lambda_{1,v} \leq \lambda_{2,v} \leq -1$ for each real place $v$. Here we follow \cite{Moriyama2004} for the choice of the Cartan subalgebra in $\frak{sp}_4(\R)$ and the positive systems.
Let $f = \bigotimes_v f_{ v} \in \itPi$ be a non-zero cusp form satisfying the following conditions:
\begin{itemize}
\item $f_{v}\mbox{ is a paramodular newform of $\itPi_v$ for all finite places }v$;
\item $f_{v}\mbox{ is a lowest weight vector of the minimal }{\rm U}(2)\mbox{-type of }D_{(-\lambda_{2,v},-\lambda_{1,v})}$ for all real places $v$.
\end{itemize}
These conditions characterize $f \in \itPi$ up to scalars. Let $W_{f}$ be the Whittaker function of $f$ with respect to $\psi_U$ defined by
\[
W_{f}(g)=\int_{U(\F)\backslash U(\A_\F)}f(ug)\overline{\psi_U(u)}\,du^{\rm Tam}.
\]
Here $U$ is the standard maximal unipotent subgroup of $\GSp_4$, $\psi_U$ is the standard non-degenerate character of $U(\F)\backslash U(\A_\F)$ (see \S\,\ref{S:notation} for precise definition), and $du^{\rm Tam}$ is the Tamagawa measure on $U(\A_\F)$.
We may decompose $W_f = \prod_v W_v$ as a product of local Whittaker functions of $\itPi_v$ with respect to $\psi_{U,v}$. We normalize $f$ as follows:
\begin{itemize}
\item for finite place $v$, we have
\[
W_v(\diag( \varpi_v^{-3\frak{c}_v}, \varpi_v^{-2\frak{c}_v}, 1, \varpi_v^{-\frak{c}_v} ))=1,
\]
where $\varpi_v$ is a uniformizer of $\F_v$ and $\frak{c}_v$ is the valuation of the different ideal of $\F_v$;
\item for real place $v$, we have
\begin{align}\label{E:normalization}
\begin{split}
W_v(1)
&=e^{-2\pi}\int_{c_1-\sqrt{-1}\infty}^{c_1+\sqrt{-1}\infty}\frac{ds_1}{2\pi \sqrt{-1}}\,\int_{c_2-\sqrt{-1}\infty}^{c_2+\sqrt{-1}\infty}\frac{ds_2}{2\pi \sqrt{-1}}\,2^{-s_1-s_2}\Gamma_\R(s_1+\lambda_{1,v}+1)\Gamma_\R(-s_2-\lambda_{2,v})\\
&\quad\quad\quad\quad\quad\quad\quad\quad\quad\quad\quad\quad\times\Gamma_\R(s_1+s_2+\lambda_{1,v}-\lambda_{2,v}+2)\Gamma_\R(s_1+s_2+\lambda_{1,v}+\lambda_{2,v}+2),
\end{split}
\end{align}
where $\Gamma_\R(s)=\pi^{-s/2}\Gamma(\tfrac{s}{2})$ and $c_1,c_2 \in \R$ satisfy $c_1+c_2+\lambda_{1,v}+\lambda_{2,v}+2>0$ and $c_1+\lambda_{1,v}+1>0>c_2+\lambda_{2,v}$. 
\end{itemize}
We write $f=f_\itPi$ and call it the normalized newform of $\itPi$.
We defined the normalized newform $f_{\itPi^\vee}$ of $\itPi^\vee$ in a similar way.
Let $\Vert f_\itPi \Vert$ be the Petersson norm of $f_\itPi$ defined by 
\[
\Vert f_\itPi \Vert=\int_{\A_\F^\times\GSp_4(\F)\backslash \GSp_4(\A_\F)}f_\itPi(g)f_{\itPi^\vee}(g \cdot{\rm diag}(-1,-1,1,1)_\infty)\,dg^{\rm Tam}.
\]
Here $dg^{\rm Tam}$ is the Tamagawa measure on $\A_\F^\times\backslash\GSp_4(\A_\F)$.

For $\sigma \in {\rm Aut}(\C)$, let ${}^\sigma\!\itPi$ be the irreducible admissible representation of $\GSp_4(\A_\F)$ defined by
\[
{}^\sigma\!\itPi = {}^\sigma\!\itPi_{\infty} \otimes {}^\sigma\!\itPi_f,
\]
where ${}^\sigma\!\itPi_f$ is the $\sigma$-conjugate of $\itPi_f = \bigotimes_{v\nmid \infty} \itPi_v$ and ${}^\sigma\!\itPi_{\infty}$ is the representation of $\GSp_4(\F_\infty)$ so that its $v$-component is equal to $\itPi_{\sigma^{-1} \circ v}$.
Assume further that $\itPi$ is motivic, that is, there exists ${\sf w}\in\Z$ such that $|\omega_\itPi| = |\mbox{ }|_{\A_\F}^{\sf w}$ and 
\[
\lambda_{1,v} - \lambda_{2,v} \equiv {\sf w}\,({\rm mod}\,2)
\]
for all real places $v$.
In Lemma \ref{L:globally generic} below, we show that ${}^\sigma\!\itPi$ is cuspidal automorphic and globally generic.
The rationality field $\Q(\mathit{\Pi})$ of $\mathit{\Pi}$ is the fixed field of
$\left\{\sigma \in {\rm Aut}(\C) \, \vert \, {}^\sigma\!\itPi = \itPi \right\}$ and is a number field.

The following theorem is our main result on the algebraicity of the critical adjoint $L$-value $L(1,\itPi,{\rm Ad})$ in terms of the Petersson norm of the normalized newform of $\itPi$.

\begin{thm}\label{T:main} 
Let $\itPi$ be an irreducible motivic globally generic cuspidal automorphic representation of $\GSp_4(\A_\F)$.
For $\sigma \in {\rm Aut}(\C)$, we have
\[
\sigma\left(\frac{L(1,\itPi,{\rm Ad})}{\zeta_\F(2)\zeta_\F(4)\cdot\Vert f_\itPi \Vert}\right) = \frac{L(1,{}^\sigma\!\itPi,{\rm Ad})}{\zeta_\F(2)\zeta_\F(4)\cdot\Vert f_{{}^\sigma\!\itPi} \Vert}.
\]
Here $\zeta_\F(s)$ is the completed Dedekind zeta function of $\F$. In particular, we have
\[
\frac{L(1,\itPi,{\rm Ad})}{\zeta_\F(2)\zeta_\F(4)\cdot\Vert f_\itPi \Vert} \in \Q(\itPi).
\]
\end{thm}

\begin{rmk}
In \cite{CI2019}, we compute the ratio explicitly when $\omega_\itPi$ is trivial and the paramodular conductor of $\itPi$ is square-free. The theorem can be regarded as a generalization of \cite{CI2019}.
\end{rmk}

\begin{rmk}
The Petersson norm $\Vert f_\itPi \Vert $ can be factorized into product of periods which are obtained by comparing the rational structures via the Whittaker model and via the coherent cohomology (cf.\,\cite{HK1992}).
We expect these periods to capture the transcendental part of the critical values of certain automorphic $L$-functions. Moreover, the expected period relation implies that Theorem \ref{T:main} is compatible with Deligne's conjecture for $L(1,\itPi,{\rm Ad})$.
This is an ongoing project considered by the author.
\end{rmk}

\subsection{An outline of the proof}

The first step is to show that for any non-trivial additive character $\psi$ of $\F\backslash\A_\F$, we have
\begin{align}\label{E:outline 1}
\Vert f_\itPi \Vert = C\cdot \frac{L(1,\itPi,{\rm Ad})}{\zeta_\F(2)\zeta_\F(4)}\cdot\prod_v C_{\psi_v}(\itPi_v)
\end{align}
for some non-zero constant $C_{\psi_v}(\itPi_v)$ depending only on $\itPi_v$ and $\psi_v$, and some constant $C \in \Q^\times$ depending only on $\F$ and the type of $\itPi$ (stable or endoscopic). 
This equality can be proved by proceeding exactly as in the proof of \cite[Proposition 5.4]{CI2019} (see also \cite[\S\,1.2]{CI2019} for brief outline), subject to the well-definedness of $C_{\psi_v}(\itPi_v)$. 
We have
\[
C_{\psi_v}(\itPi_v) = \frac{1}{L(1,\itPi_v,{\rm Ad})}\cdot\frac{\mathcal{Z}_v(1,\itPi_v,\mathcal{F}_{\psi_v})}{Z_v(\tfrac{1}{2},\itPi_v,F_{\psi_v})},
\]
where $\mathcal{F}_{\psi_v}$ and $F_{\psi_v}$ are sections in degenerate principal series representations of $\GSp_8(\F_v)$, and $\mathcal{Z}_v(s,\itPi_v,\mathcal{F}_{\psi_v})$ and $Z_v(s,\itPi_v,F_{\psi_v})$ are local zeta integral for $\itPi_v \times \itPi_v^\vee$ and doubling local zeta integral for $\itPi_v$ defined and studied by Jiang \cite{Jiang1996} and Piatetski-Shapiro and Rallis \cite{LNM1254}, respectively. 
We prove that the constant is well-defined by showing that $\mathcal{Z}_v$ and $Z_v$ converge absolutely for ${\rm Re}(s) \geq 1$ and ${\rm Re}(s)\geq \tfrac{1}{2}$, respectively, and there exists a local Siegel--Weil section $F_{\psi_v}$ such that $Z_v(\tfrac{1}{2},\itPi_v,F_{\psi_v})$ is non-zero.
The next step is to prove
\begin{align}\label{E:outline 2}
\sigma \left( \prod_{v \nmid \infty} C_{\psi_v}(\itPi_v)\right) = \prod_{v \nmid \infty}C_{\psi_v}({}^\sigma\!\itPi_v)
\end{align}
for all $\sigma \in {\rm Aut}(\C)$.
First we show that for all $v \nmid \infty$ and $\sigma \in {\rm Aut}(\C)$, we have
\begin{align}\label{E:outline 3}
\sigma Z_v(\tfrac{1}{2},\itPi_v,F_{\psi_v}) = Z_v(\tfrac{1}{2},{}^\sigma\!\itPi_v,F_{{}^\sigma\!\psi_v}),\quad \sigma \mathcal{Z}_v(1,\itPi_v,\mathcal{F}_{\psi_v}) = \mathcal{Z}_v(1,{}^\sigma\!\itPi_v,\mathcal{F}_{{}^\sigma\!\psi_v}),
\end{align}
which imply that 
\begin{align}\label{E:outline 4}
\sigma C_{\psi_v}(\itPi_v) = C_{{}^\sigma\!\psi_v}({}^\sigma\!\itPi_v).
\end{align}
For $\sigma \in {\rm Aut}(\C)$, in general $\bigotimes_{v \nmid \infty}{}^\sigma\!\psi_v$ is not the finite part of a non-trivial additive character of $\F \backslash \A_\F$.
Nonetheless, since we have freedom to vary $\psi$ in (\ref{E:outline 1}), we show that (\ref{E:outline 2}) holds by a global argument together with (\ref{E:outline 4}).
By (\ref{E:outline 1}) and (\ref{E:outline 2}), we then have
\begin{align}\label{E:outline 5}
\sigma\left(\frac{L(1,\itPi,{\rm Ad})}{\zeta_\F(2)\zeta_\F(4)\cdot\Vert f_\itPi \Vert}\cdot\prod_{v \mid \infty} C_{\psi_v}(\itPi_v)\right) = \frac{L(1,{}^\sigma\!\itPi,{\rm Ad})}{\zeta_\F(2)\zeta_\F(4)\cdot\Vert f_{{}^\sigma\!\itPi} \Vert} \cdot\prod_{v \mid \infty} C_{\psi_v}(\itPi_v).
\end{align}
for all $\sigma \in {\rm Aut}(\C)$.
Finally, we show that 
\begin{align}\label{E:outline 6}
\prod_{v \mid \infty} C_{\psi_v}(\itPi_v) \in \Q^\times.
\end{align}
This is a local problem, but we address it by a global argument.
Based on the Rallis inner product formula \cite{GQT2014} and archimedean computations, we prove that Theorem \ref{T:main} holds when $\itPi$ is endoscopic. Choose an endoscopic irreducible globally generic cuspidal automorphic representation $\itPi'$ of $\GSp_4(\A_\F)$ such that $\itPi'_\infty = \itPi_\infty$. 
We thus obtain (\ref{E:outline 6}) by comparing Theorem \ref{T:main} with (\ref{E:outline 5}) for $\itPi'$.
(Strictly speaking, we need only to consider endoscopic lifts for $\F=\Q$.)

This paper is organized as follows. In \S\,\ref{S:Petersson norm}, we recall the definition of the local zeta integrals $\mathcal{Z}_v$ and $Z_v$, and state the precise form of (\ref{E:outline 1}) in Proposition \ref{P:Petersson norm}. The proposition holds subject to Lemmas \ref{L:ab doubling}-(1) and \ref{L:ab2}-(1) on the convergence of the local zeta integrals. In \S\,\ref{S:proof of thm}, we prove (\ref{E:outline 2}) in Proposition \ref{T:3.6} subject to Lemmas \ref{L:ab doubling}-(2) and \ref{L:ab2}-(2) on the Galois equivariant property (\ref{E:outline 3}) of the local zeta integrals. In \S\,\ref{S:endoscopic}, we prove in Theorem \ref{T:endoscopic case} that Theorem \ref{T:main} holds when $\itPi$ is endoscopic and $\F=\Q$. The proposition is proved based on the Rallis inner product formula and the arithmeticity of global theta lifting in Proposition \ref{P:arithmeticity}.
As we sketched above, Theorem \ref{T:main} follows from Propositions \ref{P:Petersson norm}, \ref{T:3.6}, and Theorem \ref{T:endoscopic case}.
The context of \S\,\ref{S:local zeta integrals} is purely local and we prove Lemmas \ref{L:ab doubling} and \ref{L:ab2} in this section.

\subsection{Notation}\label{S:notation}

Fix a totally real number field $\F$. Let $\o$ 
and $\frak{D}$ be the ring of integers
and the absolute discriminant of $\F$, respectively.
Let $\A=\A_\F$ be the ring of adeles of $\F$ and $\A_f$ be its finite part. We denote by $\hat{\o}$ the closure of $\o$ in $\A_f$. For a finite dimensional vector space $V$ over $\F$, let $\mathcal{S}(V(\A))$ be the space of Schwartz functions on $V(\A)$. We will write $v$ for places of $\F$ and $\infty$ the archimedean place of $\Q$.

Let $v$ be a place of $\F$. If $v$ is a finite place, let $\frak{o}_v$, $\varpi_v$, and $q_v$ be the maximal compact subring of $\F_v$, a generator of the maximal ideal of $\frak{o}_v$, and the cardinality of $\frak{o}_v / \varpi_v\frak{o}_v$, respectively. Let $|\mbox{ }|_v$ be the absolute value on $\F_v$ normalized so that $|\varpi_v|_v = q_v^{-1}$. 
If $v$ is a real place, let $|\mbox{ }|_v=|\mbox{ }|$ be the ordinary absolute value on $\F_v \simeq \R$.
For a character $\chi$ of $\F_v^\times$, let $e(\chi) \in \R$ be the exponent of $\chi$ defined so that $|\chi| = |\mbox{ }|_v^{e(\chi)}$. When $v$ is finite, let $c(\chi) \in \Z_{\geq 0}$ be the smallest integer so that $\chi$ is trivial on $1+\varpi_v^{c(\chi)}\frak{o}_v$.

Let $\zeta(s) = \zeta_\F(s) = \prod_{v}\zeta_v(s)$ be the completed Dedekind zeta function of $\F$, where $v$ ranges over the places of $\F$ and 
\[
\zeta_v(s) = \begin{cases}
(1-q_v^{-s})^{-1} & \mbox{ if $v$ is finite},\\
\pi^{-s/2}\Gamma(s/2) & \mbox{ if $v$ is real}.
\end{cases}
\]
Here $\Gamma(s)$ is the gamma function.

Let $\psi_0=\bigotimes_v\psi_{v,0}$ be the standard additive character of $\Q\backslash \A_\Q$ defined so that
\begin{align*}
\psi_{p,0}(x) & = e^{-2\pi \sqrt{-1}\,x} \mbox{ for }x \in \Z[p^{-1}],\\
\psi_{\infty,0}(x) & = e^{2\pi \sqrt{-1}\,x} \mbox{ for }x \in \R.
\end{align*}
Let $\psi = \bigotimes_v \psi_v$ be a non-trivial additive character of $\F \backslash \A$. We say $\psi$ is standard if $\psi = \psi_0\circ{\rm tr}_{\F/\Q}$.
In this case, $\psi_v$ is called the standard additive character of $\F_v$. For $a \in \F^{\times}$ (resp.\,$a \in \F_v^{\times}$), let $\psi^a$ (resp.\,$\psi_v^a$) be the additive character of $\A$ (resp.\,$\F_v$) defined by $\psi^a(x) = \psi(ax)$ (resp.\,$\psi_v^a(x) = \psi_v(ax)$). For each finite place $v$, let $\varpi_v^{c_v}\o_v$ be the largest fractional ideal of $\F_v$ on which $\psi_v$ is trivial. The absolute conductor ${\rm cond}(\psi)$ of $\psi$ is the $\hat{\o}$-submodule of $\A_f$ defined by
\[
{\rm cond}(\psi) = \prod_{v \nmid \infty} \varpi_v^{|c_v|}\o_v.
\]
We write $v \mid {\rm cond}(\psi)$ if $c_v \neq 0$.

If $S$ is a set, then we let $\mathbb{I}_S$ be the characteristic function of $S$. Let ${\rm M}_{n , m}$ be the matrix algebra of $n$ by $m$ matrices.
Let $\GSp_{2n}$ and $\Sp_{2n}$ be the symplectic similitude group and symplectic group, respectively, defined by
\[
 \GSp_{2n} =
 \left\{ g \in \GL_{2n} \, \left| \, g
 \begin{pmatrix}
  0      & \1_n \\
  - \1_n & 0
 \end{pmatrix}
 {}^t \! g = 
  \nu(g) \begin{pmatrix}
  0      & \1_n \\
  - \1_n & 0
 \end{pmatrix}, \,\nu(g) \in \GL_1
 \right. \right\},\quad \Sp_{2n} = {\rm ker}(\nu).
\]
Let
\[
{\bf B}=\left.\left\{\bp t_1 & * & *&* \\0&t_2&*&*\\0&0&\nu t_1^{-1}&0\\0&0& *&\nu t_2^{-1}  \ep\in \GSp_4  \mbox{ }\right\vert\mbox{ } t_1,t_2,\nu \in \GL_1\right\}
\]
be the standard Borel subgroup of $\GSp_4$ and $U$ be its unipotent radical. Let $\bf T \subset {\bf B}$ be the standard maximal tours of $\GSp_4$. For a non-trivial additive character $\psi$ of $\F \backslash \A$, let $\psi_U$ be the associated additive character of $U(\F) \backslash U(\A)$ defined by
\[
\psi_U\left(\bp 1 & x & *&* \\0&1&*&y\\0&0&1&0\\0&0&-x&1 \ep\right)= \psi(-x-y).
\]
We call $\psi_U$ standard if $\psi$ is standard. 
Similar notation apply to additive character $\psi_v$ of $\Q_v$.
In $\GL_2$, let $B$ be the Borel subgroup consisting of upper triangular matrices, and put
\[
{\bf a}(\nu) = \bp \nu & 0 \\ 0 & 1 \ep,\quad {\bf d}(\nu) =   \bp 1 & 0 \\ 0 & \nu \ep,\quad {\bf m}(t)= \bp t & 0 \\ 0 & t^{-1}\ep,\quad {\bf n}(x) = \bp 1 & x \\ 0 & 1\ep,\quad {\bf w} = \bp 0 & -1 \\ 1 & 0\ep
\]
for $\nu,t \in \GL_1$ and $x \in \mathbb{G}_a$. 
Let $v$ be a finite place of $\F$ and $c \in \Z_{\geq 0}$. For $n \geq 2c$, the quasi-paramodular group ${\rm K}_1(\varpi_v^n;c)$ of level $\varpi_v^n\o_v$ is the open compact subgroup of $\GSp_4(\F_v)$ consisting of $g \in \GSp_4(\F_v)$ such that $\nu(g) \in \o_v^\times$ and 
\[g \in \bp \o_v & \o_v & \varpi_v^{-n+c}\o_v & \o_v\\
\varpi_v^{n-c}\o_v & \o_v & \o_v & \o_v\\
\varpi_v^{n}\o_v & \varpi_v^{n-c}\o_v & 1+\varpi_v^c\frak{o}_v & \varpi_v^{n-c}\o_v\\
\varpi_v^{n-c}\o_v & \o_v & \o_v & \o_v
 \ep.\]

Let $\sigma \in {\rm Aut}(\C)$. Define the $
\sigma$-linear action on $\C(X)$, which is the field of formal Laurent series in variable $X$ over $\C$, as follows:
\[
{}^\sigma\!P(X) = \sum_{n \gg -\infty}^\infty\sigma(a_n) X^n
\]
for $P(X) = \sum_{n \gg -\infty}^\infty a_n X^n \in \C(X)$. For a complex representation $\itPi$ of a group $G$ on the space $\mathcal{V}_\itPi$ of $\itPi$, let ${}^\sigma\!\itPi$ of $\itPi$ be the representation of $G$ defined
\begin{align}\label{E:sigma action}
{}^\sigma\!\itPi(g) = t \circ \itPi(g) \circ t^{-1},
\end{align}
where $t:\mathcal{V}_\itPi \rightarrow \mathcal{V}_\itPi$ is a $\sigma$-linear isomorphism. Note that the isomorphism class of ${}^\sigma\!\itPi$ is independent of the choice of $t$. We call ${}^\sigma\!\itPi$ the $\sigma$-conjugate of $\itPi$. When $v$ is a finite place and $\varphi$ is a complex-valued function on $\F_v^n$ or $(\F_v^\times)^n$, we define ${}^\sigma\!\varphi(x) = \sigma(\varphi(x))$ for $x \in \F_v^n$ or $x \in (\F_v^\times)^n$.

\subsection{Measures}\label{SS:measures}
Let $v$ be a place of $\F$. 
If $v$ is finite, we normalize the Haar measures on $\F_v$ and $\F_v^\times$ so that ${\rm vol}(\o_v)=1$ and ${\rm vol}(\o_v^\times)=1$, respectively. If $v$ is real, we normalize the Haar measures on $\F_v \simeq \R$ and $\F_v^\times \simeq \R^\times$ so that ${\rm vol}([1,2])=1$ and ${\rm vol}([1,2])=\log 2$, respectively.
Let $m$ be a positive integer. Let $dg_v$ be the Haar measure on $\GL_m(\F_v)$ defined as follows: For $\phi \in L^1(\GL_m(\F_v))$, we have 
\begin{align}\label{E:standard measure on GL}
\int_{\GL_m(\F_v)}\phi(g_v)\,dg_v = \prod_{1 \leq i<j \leq m}\int_{\F_v}du_{ij}\prod_{1\leq i \leq m}\int_{k_v^{\times}}d^{\times}t_i\int_{K_v}dk\,\phi\left(\bp t_1 & u_{12} & \cdots &u_{1m} \\ 0 & t_2 & \cdots & u_{2m} \\ \vdots & \vdots & \ddots& \vdots \\ 0 & 0 & \cdots & t_m\ep k\right)\prod_{1\leq i \leq m}|t_i|_v^{-m+i},
\end{align}
where 
\begin{align*}
K_v = \begin{cases}    
\GL_m(\frak{o}_v)& \mbox{ if $v$ is finite},\\
{\rm O}(m) & \mbox{ if $v$ is real},
\end{cases}
\end{align*}
and ${\rm vol}(K_v)=1$. 
Let $H$ be a connected reductive linear algebraic group defined and split over $\F_v$. Fix a Chevalley basis of ${\rm Lie}(H)$. The basis determines a top differential form on $H$ over $\Z$ which is unique up to $\pm 1$. The top differential form together with the Haar measure on $\F_v$ determines a Haar measure on $H(\F_v)$ as explained in \cite[\S\,6]{Vos1996}. We call it the local Tamagawa measure on $H(\F_v)$.
For any compact group $K$, we take the Haar measure on $K$ such that ${\rm vol}(K)=1$.


\subsection{Weil representation}

Let $(V,(\mbox{ },\mbox{ }))$ be a non-degenerate quadratic space of even dimension $m$ over $\F$. Define the orthogonal similitude group ${\rm GO}(V)$ by
\[
{\rm GO}(V) = \{ h \in \GL(V)\mbox{ }\vert \mbox{ }(hx,hy)=\nu(h)(x,y) \mbox{ for }x,y \in V \},
\]
here $\nu : {\rm GO}(V) \rightarrow \GL_1$ is the scale map. Let
\[
{\rm GSO}(V) = \{h \in {\rm GO}(V) \mbox{ }\vert \mbox{ }\det(h)=\nu(h)^{m/2} \}.
\]
Let ${\rm O}(V)$ and ${\rm SO}(V)$ be the orthogonal group and special orthogonal group defined by
\begin{align*}
{\rm O}(V) &= \{ h \in {\rm GO}(V) \mbox{ }\vert \mbox{ } \nu(h)=1 \},\\
{\rm SO}(V) &= \{ h \in {\rm GO(V)} \mbox{ }\vert \mbox{ } \det(h)=\nu(h)=1 \}.
\end{align*}
Let $\psi = \bigotimes_v \psi_v$ be a non-trivial additive character of $\F \backslash \A$. We denote by $\omega_{\psi,V,n} = \bigotimes_v \omega_{\psi_v,V,n}$ the Weil representation of 
$\Sp_{2n}(\A) \times {\rm O}(V)(\A)$ on $\mathcal{S}(V^n(\A))$ with respect to $\psi$ (cf.\,\cite[\S\,5]{Kudla1994} and \cite[\S\,4.2]{Ichino2005} for explicit formulas). Let $S(V^n(\A))$ be the subspace of $\mathcal{S}(V^n(\A))$
consisting of functions which correspond to polynomials in the Fock model
at the archimedean places. Let 
\[{\rm G}(\Sp_{2n} \times {\rm O}(V)) = \left\{ (g,h) \in \GSp_{2n}\times {\rm O}(V) \mbox{ }\vert\mbox{ }\nu(g)=\nu(h) \right\}.\]
We extend $\omega_{\psi,V,n}$ to a representation of ${\rm G}(\Sp_{2n} \times {\rm O}(V))(\A)$ as follows:
\[
\omega_{\psi,V,n}(g,h)\varphi = \omega_{\psi,V,n}\left (g \bp {\bf 1}_n & 0 \\ 0 & \nu(g)^{-1}{\bf 1}_n\ep,1 \right )L(h)\varphi
\]
for $(g,h) \in {\rm G}(\Sp_{2n} \times {\rm O}(V))(\A)$ and $\varphi \in \mathcal{S}(V^n(\A))$. Here
\[
L(h)\varphi(x)=\vert \nu(h)  \vert_\A^{-nm/4}\varphi(h^{-1}x).
\]

\subsection{Adjoint $L$-functions}

Let $\itPi=\bigotimes_v \itPi_v$ be an irreducible globally generic cuspidal automorphic representation of $\GSp_4(\A)$. 
By \cite{AH2006} and \cite[Theorem 12.1]{GT2011}, $\itPi$ has a strong functorial lift $\itPsi$ to $\GL_4(\A)$.
By \cite[Theorem 12.1]{GT2011}, 
either $\itPsi$ is cuspidal or $\itPsi = \tau_1 \boxplus \tau_2$ for some irreducible cuspidal automorphic representations $\tau_1$ and $\tau_2$
of $\GL_2(\AA)$ with equal central character such that $\tau_1 \neq \tau_2$.
We say that $\itPi$ is stable (resp.~endoscopic)
if $\itPsi$ is cuspidal (resp.~non-cuspidal).

Recall that the dual group of $\GSp_4$ is $\GSp_4(\C)$.
Let ${\rm Ad}$ denote the adjoint representation of $\GSp_4(\C)$ on $\frak{pgsp}_4(\C)$, and $\rm std$ the composition of the projection $\GSp_4(\C) \rightarrow \PGSp_4(\C)$ with the standard
representation of $\PGSp_4(\C) \simeq {\rm SO}_5(\C)$ on $\C^5$.
Let $S$ be a finite set of places of $\F$
including the archimedean places such that,
for $v \notin S$, $\itPi_v$ is unramified.
Then the partial adjoint and standard $L$-functions of $\itPi$ are defined as the Euler products
\[
L^S(s,\itPi,{\rm Ad}) = \prod_{v \notin S} L(s,\itPi_v,{\rm Ad}),\quad L^S(s,\itPi,{\rm std}) = \prod_{v \notin S} L(s,\itPi_v,{\rm std})
\]

for $s \in \C$, which are absolutely convergent for ${\rm Re}(s)$ sufficiently large.
Also, we have 
\begin{align*}
 L^S(s, \itPsi, \Sym^2 \otimes \,\omega_\itPi^{-1}) & = L^S(s, \itPi, \Ad), \\
 L^S(s, \itPsi, \wedge^2\otimes \omega_\itPi^{-1}) & = \zeta^S(s) L^S(s, \itPi, {\rm std}).
\end{align*}
In particular, $L^S(s,\itPi,{\rm Ad})$ and $L^S(s,\itPi,{\rm std})$ admit meromorphic continuations to $\C$.
(In a more general context, the meromorphic continuation of $L^S(s,\itPi,{\rm std})$ was established by Piatetski-Shapiro and Rallis \cite{LNM1254} much earlier.)
By \cite[Theorem 12.1]{GT2011}, 
$L^S(s, \itPsi, \wedge^2 \otimes \omega_\itPi^{-1})$ has a simple (resp.~double) pole at $s=1$
if $\itPi$ is stable (resp.~endoscopic).
Hence $L^S(s, \itPi, {\rm std})$ is holomorphic and non-zero
(resp.~has a simple pole) at $s=1$ if $\itPi$ is stable (resp.~endoscopic).
In particular, $L^S(s, \itPi, \Ad)$ is holomorphic and non-zero at $s=1$.

For any place $v$ of $\F$, we denote by $\phi_{\itPi_v} : L_{\F_v} \rightarrow \GSp_4(\C)$ the local $L$-parameter attached to $\itPi_v$ by the local Langlands correspondence established by Gan and Takeda \cite{GT2011} if $v$ is finite and by Langlands \cite{Langlands1989} if $v$ is real. 
Here $L_{\F_v}$ is the Weil--Deligne group of $\F_v$ if $v$ is finite but the Weil group of $\F_v$ if $v$ is real.
Since $\itPsi_v$ is essentially unitary and generic (and hence ``almost tempered''), the adjoint $L$-factor 
\[
 L(s,\itPi_v,{\rm Ad}) = L(s,{\rm Ad}\circ \phi_{\itPi_v})
\]
defined as in \cite[\S\,3]{Tate1979} is holomorphic at $s=1$.
In fact, the same holds for any irreducible admissible generic  representation of $\GSp_4(\F_v)$ (see \cite[Conjecture 2.6]{GP1992}, \cite{AS2008}, \cite{GT2011}, \cite[Proposition B.1]{GI2016}).
Hence the completed adjoint $L$-function $L(s, \itPi, \Ad)$ is holomorphic and non-zero at $s=1$.

\section{Formula for Petersson norms}\label{S:Petersson norm}

Let $H = \GSp_8$ and
\[
 {\bf G} = \{ (g_1, g_2) \in \GSp_4 \times \GSp_4 \, | \, \nu(g_1) = \nu(g_2) \}.
\]
Denote by $Z_{H}$ the center of $H$.
We identify $\bf G$ with its image under the embedding
\begin{align}\label{E:doubling embedd.}
 {\bf G}  \longrightarrow H, \quad
 \left(
 \begin{pmatrix}
  a_1 & b_1 \\
  c_1 & d_1
 \end{pmatrix},
 \begin{pmatrix}
  a_2 & b_2 \\
  c_2 & d_2
 \end{pmatrix}
 \right)  \longmapsto
 \begin{pmatrix}
  a_1 & 0    & b_1 & 0 \\
  0   & a_2  & 0   & -b_2 \\
  c_1 & 0    & d_1 & 0 \\
  0   & -c_2 & 0   & d_2 
 \end{pmatrix}.
\end{align}
Let $V_{3,3} = \F^{6}$ be the space of column vectors
equipped with a non-degenerate symmetric bilinear form $(\ ,\ )$ given by
\[
 (x, y) = {}^t \! x 
 \begin{pmatrix}
  0    & \1_3 \\
  \1_3 & 0
 \end{pmatrix}
 y
\]
for $x,y \in V_{3,3}$. With respect to the standard basis of $V_{3,3}$, we identify ${\rm GO}(V_{3,3})$ with the split orthogonal similitude group
\[
 {\rm GO}_{3,3} = \left\{ h \in \GL_{6} \, \left| \, {}^t \! h
 \begin{pmatrix}
  0    & \1_3 \\
  \1_3 & 0
 \end{pmatrix}
 h =\nu(h) 
 \begin{pmatrix}
  0    & \1_3 \\
  \1_3 & 0
 \end{pmatrix},\,\nu(h) \in \GL_1
 \right. \right\}.
\]
For a non-trivial additive character $\psi_v$ of $\F_v$, we write $\omega_{\psi_v}=\omega_{\psi_v,V_{3,3},4}$ for the Weil representation of $\Sp_8(\F_v) \times {\rm O}_{3,3}(\F_v)$ on $\mathcal{S}(V_{3,3}^4(\F_v))$ with respect to $\psi_v$.

Let $\itPi = \bigotimes_v \itPi_v$ be an irreducible globally generic cuspidal automorphic representation of $\GSp_4(\A)$ with central character $\omega_\itPi$. 
We assume 
\[
\itPi_{v}\vert_{\Sp_4(\F_v)} = D_{(\lambda_{1,v},\lambda_{2,v})}\oplus D_{(-\lambda_{2,v},-\lambda_{1,v})}
\]
for each real place $v$,
where $D_{(\lambda_{1,v},\lambda_{2,v})}$ is the discrete series representation of $\Sp_4(\F_v)\simeq \Sp_4(\R)$ with Blattner parameter $(\lambda_{1,v},\lambda_{2,v}) \in \Z^2$ such that $2-\lambda_{1,v} \leq \lambda_{2,v} \leq -1$. For each finite place $v$, by the newform theory of Robert--Schmidt \cite[Theorem 7.5.4]{RS2007} and Okazaki \cite[Main Theorem]{Okazaki2019}, there exists a smallest non-negative integer $n_v \geq 2c(\omega_{\itPi_v})$ such that $\itPi_v^{{\rm K}_1(\varpi_v^{n_v};\,c(\omega_{\itPi_v}))} \neq 0$. In this case, we have ${\rm dim}\,\itPi_v^{{\rm K}_1(\varpi_v^{n_v};\,c(\omega_{\itPi_v}))} =1$. A paramodular newform of $\itPi_v$ is a non-zero vector in this one-dimensional space. The paramodular conductor ${\rm cond}(\itPi)$ of $\itPi$ is the $\hat{\o}$-submodule of $\A_f$ defined by
\[
{\rm cond}(\itPi) = \prod_{v \nmid \infty}\varpi_v^{n_v}\o_v.
\]
We write $v \mid {\rm cond}(\itPi)$ if $n_v >0$. For a place $v$ and a non-trivial additive character $\psi_v$ of $\F_v$, we denote by $\mathcal{W}(\itPi_v,\psi_{U,v})$ the space of Whittaker functions of $\itPi_v$ with respect to $\psi_{U,v}$.

\subsection{Doubling local zeta integrals}\label{SS:local zeta 1}

Let $P$ be the standard Siegel parabolic subgroup of $H$ defined by 
\[
P = \left\{ \left.
 \begin{pmatrix}
  a & * \\
  0 & \nu {}^t \! a^{-1}
 \end{pmatrix}
 \in H \, \right| \, a \in \GL_4, \, \nu \in \GL_1 \right\}.
\]
Let $v$ be a place of $\F$. Denote by $I_v(s)$ the degenerate principal series representation
${\rm Ind}_{P(\F_v)}^{H(\F_v)}(\delta_P^{s/5})$ of $H(\F_v)$. Here $\delta_P$ is the modulus character of $P(\F_v)$ given by
\[
\delta_P \left(\begin{pmatrix}
  a & * \\
  0 & \nu {}^t \! a^{-1}
 \end{pmatrix} \right) = |\det(a)|_v^5|\nu|_v^{-10}.
\]
Let $\phi_v$ be a matrix coefficient of $\itPi_v$ and $F_v \in I_v(s)$ be a holomorphic section. We define the local zeta integral 
\begin{align}\label{E:local zeta1}
Z_v(s,\phi_v, F_v) = \int_{\Sp_4(\F_v)} F_v(\delta (g_v,1), s) \phi_v(g_v) \, dg_v,
\end{align}
where
\[
 \delta =
 \begin{pmatrix}
  0 & 0 & -\frac{1}{2} \1_2 & \frac{1}{2} \1_2 \\
  \frac{1}{2} \1_2 & \frac{1}{2} \1_2 & 0 & 0 \\
  \1_2 & - \1_2 & 0 & 0 \\
  0 & 0 & \1_2 & \1_2
 \end{pmatrix}.
\]
Here $dg_v$ is a Haar measure on $\Sp_4(\F_v)$ normalized so that ${\rm vol}(\Sp_4(\o_v),dg_v)=1$ if $v$ is finite and is the local Tamagawa measure if $v$ is real.
Note that
\begin{align}\label{E:diagonal}
\delta (g,g) \delta^{-1}\in P(\F_v)
\end{align}
for all $g \in \GSp_4(\F_v)$ and all places $v$.

Let $v$ be a finite place. Let $\phi_v$ be a matrix coefficient of $\itPi_v$ and $F_v \in I(s)$. For $\sigma \in {\rm Aut}(\C)$, we define the matrix coefficient ${}^\sigma\!\phi$ of ${}^\sigma\!\itPi_v$ (cf.\,Lemma \ref{L:sigma matrix coeff.}) and ${}^\sigma\!F$ by
\[
{}^\sigma\!\phi_v(g) = \sigma(\phi_v(g)),\quad {}^\sigma\!F_v(h,s) = \sigma(F_v(h,s))
\]
for $g \in \GSp_4(\F_v)$ and $h \in H(\F_v)$. Note that ${}^\sigma\! F_v\vert_{s=\tfrac{n}{2}} \in I_v(\tfrac{n}{2})$ for all odd integers $n$.

\begin{lemma}\label{L:ab doubling}

Let $\phi_v$ be a matrix coefficient of $\itPi_v$ and $F_v \in I_v(s)$ be a holomorphic section.

(1) The integral $Z_v(s, \phi_v, F_v)$ is absolutely convergent for ${\rm Re}(s) \geq \tfrac{1}{2}$. 

(2) Assume $v$ is finite. For $\sigma \in {\rm Aut}(\C)$, we have
\[
\sigma Z_v(\tfrac{1}{2},\phi_v,F_v) = Z_v(\tfrac{1}{2},{}^\sigma\!\phi_v,{}^\sigma\!F_v).
\]

\end{lemma}

\begin{proof}
The assertions will be proved in Proposition \ref{P:ab doubling} below.
\end{proof}

We recall the local Siegel-Weil sections. Let $\psi_v$ be a non-trivial additive character of $\F_v$. Define a $H(\F_v)$-intertwining map
\begin{align}\label{E:SW section}
 \mathcal{S}(V_{3,3}^4(\F_v))  \longrightarrow I_v (\tfrac{1}{2}), \quad
 \varphi  \longmapsto F_{\psi_v}(\varphi)
\end{align}
by 
\[
 F_{\psi_v}({\varphi}) (g, \tfrac{1}{2}) = \omega_{\psi_v}(g, h) \varphi(0),
\]
where $\nu(g)=\nu(h)$. We extend $F_{\psi_v}({\varphi})$ to a holomorphic section $F_{\psi_v}({\varphi})$ of $I_v(s)$
such that its restriction to ${\bf K}_v$ is independent of $s$, where ${\bf K}_v$ is the maximal compact subgroup of $H(\F_v)$ defined by
\[
 {\bf K}_v =
 \begin{cases}
  H(\o_v)         & \text{if $v$ is finite,} \\
  H(\R) \cap \O(2n) & \text{if $v$ is real.}
 \end{cases}
\]

\begin{lemma}\label{lem:nonvanish}
Let $\psi_v$ be a non-trivial additive character of $\F_v$. There exists $\varphi \in S(V_{3,3}^4(\F_v))$ such that 
\[
 Z_{v} (\tfrac{1}{2}, \phi_{v}, F_{\psi_v}(\varphi)) \ne 0.
\]
\end{lemma}

\begin{proof}
The assertion was proved for real $v$ in \cite[Lemma 5.3]{CI2019}. We assume $v$ is finite. Let $\itPi_0$ be an irreducible component of $\itPi_v\vert_{\Sp_4(\F_v)}$.
Fix a bilinear equivariant pairing $\langle \ , \ \rangle$ on $\itPi_0 \times \itPi_0^\vee$.
For $f_1\in \itPi_0$ and $f_2 \in \itPi_0^\vee$,
define a matrix coefficient $\phi_{f_1 \otimes f_2}$ of $\itPi_0$ by
\[
 \phi_{f_1 \otimes f_2}(g) = \langle \itPi_0 (g) f_1, f_2 \rangle.
\]
Then the local zeta integral $Z_v(\tfrac{1}{2}, \phi_{f_1 \otimes f_2}, F)$
is absolutely convergent by Lemma \ref{L:ab doubling}-(1), and defines
an $\Sp_4(\F_v) \times \Sp_4(\F_v)$-intertwining map
\[
\ell: I_v(\tfrac{1}{2})
  \longrightarrow \itPi_0^{\vee} \otimes \itPi_0, \quad
 F_v  \longmapsto \left[ f_1 \otimes f_2 \mapsto
 Z_v (\tfrac{1}{2}, \phi_{f_1 \otimes f_2}, F_v) \right].
\]

Let $V_0$ be the quaternionic quadratic space of dimension $4$ over $\F_v$. Let $V_1$ (resp.~$V_2$) be the split quadratic space (resp.~quaternionic quadratic space) of dimension $6$ over $\F_v$. For $i=1,2$, let $R(V_i)$ be the image of the $\Sp_8(\F_v)$-intertwining map
\[
 S(V_i^4(\F_v))  \longrightarrow I_v (\tfrac{1}{2}), \quad \varphi  \longmapsto F_{\psi_v}^{(i)}({\varphi}),
\]
where $F_{\psi_v}^{(i)}({\varphi})(g, \tfrac{1}{2}) = \omega_{\psi_v, V_i, 4}(g, 1) \varphi(0)$. By \cite{KR1992}, 
\begin{align}\label{E:KR relation}
\begin{split}
I_v(\tfrac{1}{2})= R(V_1)+R(V_2),\quad 
R(V_0) \simeq R(V_2) / R(V_1) \cap R(V_2).
\end{split}
\end{align}
By 
\cite[Proposition 7.2.1]{KR1994}, the intertwining map $\ell$ is non-zero. If $\ell \vert_{R(V_1)}$ is zero, then it follows from (\ref{E:KR relation}) that $\ell \vert_{R(V_2)}$ must be non-zero and defines a non-zero element in 
\[
{\rm Hom}_{\Sp_4(\F_v)\times \Sp_4(\F_v)}(R(V_0),\itPi_0^\vee \otimes \itPi_0).
\]
Therefore the local theta lift of $\itPi_0$ to ${\rm O}(V_0)(\F_v)$ is non-zero by \cite[Proposition 3.1]{HKS1996}. This contradicts the genericity of $\itPi_0$ (cf.\,\cite[Corollary 4.2-(i)]{GT2011b}). Hence $\ell \vert_{R(V_1)}$ is non-zero. This completes the proof.

\end{proof}

\subsection{Local zeta integrals for $\GSp_4 \times \GSp_4$}\label{SS:local zeta 2}

Let $\mathcal{P}$ be the parabolic subgroup of $H$ defined by
\[
\mathcal{P} = \left\{ \left.
 \begin{pmatrix}
  a & *  & * & * \\
  0 & a' & * & b' \\
  0 & 0 & \nu {}^t \! a^{-1} & 0 \\
  0 & c' & * & d'
 \end{pmatrix}
 \in H \, \right| \,
 \begin{array}{l}
  a \in \GL_3 \\
  \begin{pmatrix}
   a' & b' \\
   c' & d'
  \end{pmatrix}
  \in \GL_2 \\
  \nu = a' d' - b' c' \in \GL_1
 \end{array}
 \right\}. 
\]
Let $v$ be a place of $\F$ and $\psi_v$ be a non-trivial additive character of $\F_v$. Denote by $\mathcal{I}_v(s)$ the degenerate principal series representation
${\rm Ind}_{\mathcal{P}(\F_v)}^{H(\F_v)}(\delta_{\mathcal P}^{s/6})$ of $H(\F_v)$. Here $\delta_{\mathcal{P}}$ is the modulus character of $\mathcal{P}(\F_v)$ given by
\[
\delta_{\mathcal P} \left(\begin{pmatrix}
  a & *  & * & * \\
  0 & a' & * & b' \\
  0 & 0 & \nu {}^t \! a^{-1} & 0 \\
  0 & c' & * & d'
 \end{pmatrix}
 \right) = |\det(a)|_v^6|\nu|_v^{-9}.
\]
Let $W_{1,v} \in \mathcal{W}(\itPi_v,\psi_{U,v})$ and $W_{2,v} \in \mathcal{W}(\itPi_v^\vee,\psi_{U,v}^{-1})$ be Whittaker functions of $\itPi_v$ and $\itPi_v^\vee$ with respect to $\psi_{U,v}$ and ${\psi}_{U,v}^{-1}$, respectively, and $\mathcal{F}_v \in \mathcal{I}_v(s)$ be a holomorphic section. We define the local zeta integral
\begin{align}\label{E:local zeta2}
 \calZ_v(s, W_{1,v}, W_{2,v}, \calF_v) = 
 \int_{Z_H(\F_v) \tilde{U}(\F_v) \backslash {\bf G}(\F_v)}
 \calF_v(\eta g_v, s) (W_{1,v}\otimes W_{2,v})(g_v) \, d\bar{g}_v,
\end{align}
where
\[
 \tilde{U} = \left\{ \left. \left( u
  \begin{pmatrix}
   1 & 0 & x & 0 \\
   0 & 1 & 0 & 0 \\
   0 & 0 & 1 & 0 \\
   0 & 0 & 0 & 1
 \end{pmatrix},
 u \right) \, \right| \, u \in U, \, x \in \GG_a \right\}
\]
and 
\[
 \eta =
 \begin{pmatrix}
   1 & 0 & 0  & 0  & 0 & 0 & 0 & 0 \\
   0 & 1 & 0  & 0  & 0 & 0 & 0 & 0 \\
   0 & 0 & 0  & 0  & 0 & 0 & 0 & 1 \\
   0 & 0 & 0  & 0  & 0 & 0 & 1 & 0 \\
   0 & 0 & 0  & 0  & 1 & 0 & 1 & 0 \\
   0 & 0 & 0  & 0  & 0 & 1 & 0 & 1 \\
   0 & 1 & 0  & -1 & 0 & 0 & 0 & 0 \\
   1 & 0 & -1 & 0  & 0 & 0 & 0 & 0
 \end{pmatrix}.
\]
Here $d\bar{g}_v$ is the quotient measure normalized so that ${\rm vol}(Z_H(\o_v) \tilde{U}(\o_v) \backslash {\bf G}(\o_v),d\bar{g}_v)=1$ if $v$ is finie and is the quotient measure defined by the local Tamagawa measures on $Z_H(\F_v)\backslash {\bf G}(\F_v)$ and $\tilde{U}(\F_v)$ if $v$ is real.

Let $v$ be a finite place. Let $W_v \in \mathcal{W}(\itPi_v,\psi_{U,v})$ be a Whittaker function and $\mathcal{F}_v \in \mathcal{I}_v(s)$. 
For $\sigma \in {\rm Aut}(\C)$, we define ${}^\sigma W_v \in \mathcal{W}({}^\sigma\!\itPi_v,{}^\sigma\!\psi_{U,v})$ and ${}^\sigma\!\mathcal{F}_v \in \mathcal{I}_v(s)$ by
\[
{}^\sigma W_v(g) = \sigma(W_v(g)),\quad {}^\sigma\!\mathcal{F}_v(h,s) = \sigma(\mathcal{F}_v(h,s))
\]
for $g \in \GSp_4(\F_v)$ and $h \in H(\F_v)$. Note that ${}^\sigma\!\mathcal{F}_v\vert_{s=m} \in \mathcal{I}_v(m)$ for all odd integers $m$.

\begin{lemma}\label{L:ab2}

Let $W_{1,v} \in \mathcal{W}(\itPi_v,\psi_{U,v})$ and $W_{2,v} \in \mathcal{W}(\itPi_v^\vee,{\psi}_{U,v}^{-1})$ be Whittaker functions and $\mathcal{F}_v \in \mathcal{I}_v(s)$ be a holomorphic section.

(1) The integral $\mathcal{Z}_v(s, W_{1,v}, W_{2,v}, \mathcal{F}_v)$ is absolutely convergent for ${\rm Re}(s) \geq 1$. 

(2) Assume $v$ is finite. For $\sigma \in {\rm Aut}(\C)$, we have
\[
\sigma \mathcal{Z}_v(1,W_{1,v},W_{2,v},\mathcal{F}_v) = \mathcal{Z}_v(1, {}^\sigma W_{1,v}, {}^\sigma W_{2,v},{}^\sigma\!\mathcal{F}_v).
\]

\end{lemma}

\begin{proof}
The assertions will be proved in Proposition \ref{P:ab2} below.
\end{proof}

Let $v$ be a place of $\F$ and $\psi_v$ a non-trivial additive character of $\F_v$. When $v$ is finite, let $\varpi_v^{d_v}\o_v$ be the largest fractional ideal of $\F_v$ on which $\psi_v$ is trivial. Let ${\bf K}_v'$ and $K_v'$ be the maximal compact subgroups of ${\rm GO}_{3,3}(\F_v)$ and ${\rm O}_{3,3}(\F_v)$, respectively, defined by
\begin{align*}
{\bf K}_v' &= \begin{cases}
  \bp {\bf 1}_3 & 0 \\ 0 & \varpi_v^{-d_v} {\bf 1}_3\ep {\rm GO}_{3,3}(\o_v)\bp {\bf 1}_3 & 0 \\ 0 & \varpi_v^{d_v} {\bf 1}_3\ep & \text{if $v$ is finite,} \\
  {\rm GO}_{3,3}(\R) \cap \O(6) & \text{if $v$ is real}.
  \end{cases}\\
K_v' &= {\bf K}_v' \cap {\rm O}_{3,3}(\F_v).   
\end{align*}
Note that ${\bf K}_v'$ and $K_v'$ depend on the additive character $\psi_v$ when $v$ is finite. Let 
\[
f_{\psi_v}^\circ \in {\rm Ind}_{P'(\F_v)}^{{\rm GO}_{3,3}(\F_v)}(\delta_{P'}^{s/2})
\]
be the ${\bf K}_v'$-invariant section such that $f_{\psi_v}^\circ(1,s)=1$, where $P'$ is the standard Siegel parabolic subgroup of ${\rm GO}_{3,3}$.
Define the partial Fourier transform
\[
 \calS(V_{3,3}^4(\F_v))
  \longrightarrow \calS({\rm M}_{3 , 8}(\F_v)), \quad
 \varphi \longmapsto \hat{\varphi},
\]
where
\begin{align}\label{Partial Fourier}
 \hat{\varphi}(u, v) = \int_{{\rm M}_{3 , 4}(\F_v)} \varphi
 \begin{pmatrix}
  x \\
  u
 \end{pmatrix}
 \psi_v( \tr(v {}^t \! x)) \, dx
\end{align}
for $u, v \in {\rm M}_{3 , 4}(\F_v)$ and the Haar measure on $M_{3,4}(\F_v) \simeq \F_v^{12}$ is the product measure of the Haar measure on $\F_v$ defined in \S\,\ref{SS:measures}.
Let $\hat{\omega}_{\psi_v}$ be the representation of ${\rm G}(\Sp_8 \times {\rm O}_{3,3})(\F_v)$
on $\calS({\rm M}_{3 , 8}(\F_v))$ defined by 
\[
\hat{\omega}_{\psi_v}(g,h)\hat{\varphi} = (\omega_{\psi_v}(g,h)\varphi)^{\mathlarger{\hat{}}}.
\]
Define a $H(\F_v)$-intertwining map 
\begin{align}\label{E:theta integral}
\mathcal{S}(V_{3,3}^4(\F_v))  \longrightarrow \mathcal{I}_v (s), \quad
 \varphi  \longmapsto \mathcal{F}_{\psi_v}(\varphi)
\end{align}
by
\[
 \calF_{\psi_v}(\varphi)(g, s) = \int_{\GL_3(\F_v)} \int_{K_v'} \hat{\omega}_{\psi_v}(g, k_vh)
 \hat{\varphi} (0_{3 \times 4}, {}^t \! a_v, 0_{3 \times 1}) f_{\psi_v}^\circ(k_vh,s)
 |\det (a_v)|_{v}^{s+3} \,dk_v \, da_v,
\]
where $\nu(h)=\nu(g)$. The Haar measure $da_v$ on $\GL_3(\F_v)$ is normalized as in (\ref{E:standard measure on GL}). Note that the integral is absolutely convergent for ${\rm Re}(s)>-1$ and admits meromorphic continuation to $s \in \C$ (cf.\,\cite{GJ1972}). Therefore, $\mathcal{F}_{\psi_v}(\varphi)$ defines a meromorphic section of $\mathcal{I}_v(s)$ which is holomorphic for ${\rm Re}(s)>-1$.

\subsection{Petersson norms and adjoint $L$-values}\label{SS:2.3}

Let $f_\itPi = \bigotimes_v f_v \in \itPi$ and $f_{\itPi^\vee} = \bigotimes_v f_v^\vee \in \itPi^\vee$ be the normalized newforms of $\itPi$ and $\itPi^\vee$. Choose local bilinear equivariant pairings $\langle \ , \ \rangle_v$ on $\itPi_v \times \itPi_v^\vee$ such that
\[
\Vert f_\itPi \Vert = \prod_{v \mid \infty} \<f_v,\itPi_v^\vee({\rm diag}(-1,-1,1,1))f_v^\vee\>_v\cdot \prod_{v\nmid \infty} \langle f_v, f_v^\vee \rangle_v,
\]
and define a matrix coefficient $\phi_v$ of $\itPi_v$ by
\begin{align*}
 \phi_v(g) = \begin{cases}
\displaystyle{ \frac{\langle \itPi_v (g) f_v, f_v^\vee \rangle_v}{\langle f_v, f_v^\vee \rangle_v}} & \mbox{ if $v$ is finite},\\
\displaystyle{\frac{\langle \itPi_v (g) f_v, \itPi_v^\vee({\rm diag}(-1,-1,1,1))f_v^\vee \rangle_v}{\langle f_v, \itPi_v^\vee({\rm diag}(-1,-1,1,1))f_v^\vee \rangle_v}} & \mbox{ if $v$ is real}.
 \end{cases}
\end{align*}
Let $v$ be a place of $\F$ and $\psi_v$ a non-trivial additive character of $\F_v$. Let $W_{\psi_v} \in \mathcal{W}(\itPi_v,\psi_{U,v})$ be a non-zero Whittaker function satisfying the following condition:
\begin{itemize}
\item $W_{\psi_v}$ is a paramodular newform if $v$ is finite;
\item $W_{\psi_v}$ is a lowest weight vector of the minimal ${\rm U}(2)$-type of $D_{(-\lambda_{2,v},-\lambda_{1,v})}$ if $v$ is real.
\end{itemize}
The condition characterize $W_{\psi_v}$ up to scalars, we normalize it so that
\begin{itemize}
\item $W_{\psi_v}(\diag( a_v^3,a_v^2,1,a_v ))=1$ if $v$ is finite;
\item $W_{\psi_v}((\diag( a_v^3,a_v^2,1,a_v ))$ is normalized as in (\ref{E:normalization}) if $v$ is real.
\end{itemize}
Here $a_v \in \F_v^\times$ is chosen so that
\begin{itemize}
\item $a_v\o_v$ is the largest fractional ideal of $\F_v$ on which $\psi_v$ is trivial if $v$ is finite;
\item $\psi_v^{a_v}$ is the standard additive character of $\F_v$ if $v$ is real.
\end{itemize}
Note that $W_{\psi_v}(\diag( a_v^3,a_v^2,1,a_v )) \neq 0$ for finite place $v$ by the results of Robert--Schmidt \cite{RS2007} and Okazaki \cite{Okazaki2019}, hence the normalization is valid.
We call $W_{\psi_v}$ the normalized Whittaker newform of $\itPi_v$ with respect to $\psi_{U,v}$.
By definition, for $a \in \F_v^\times$, we have
\[
W_{\psi_v^a} (g) = W_{\psi_v}({\rm diag}(a^3,a^2,1,a)g)
\]
for $g \in \GSp_4(\F_v)$.
We define the normalized Whittaker newform $W_{\psi_v}^\vee$ of $\itPi_v^\vee$ with respect to $\psi_v$ in a similar way. 

\begin{prop}\label{P:Petersson norm}
Let $\psi$ be a non-trivial additive character of $\F\backslash \A$. We have 
\[
\Vert f_\itPi \Vert = 2^c \cdot \frak{D}^{-13}\cdot \frac{L(1,\itPi,{\rm Ad})}{\zeta_\F(2)\zeta_\F(4)}\cdot\prod_v C_{\psi_v}(\itPi_v).
\]
Here 
\begin{align*}
c & =\begin{cases}
      1 & \mbox{ if $\itPi$ is stable},\\
      2 & \mbox{ if $\itPi$ is endoscopic},
     \end{cases}
\end{align*}
and $C_{\psi_v}(\itPi_v)$ is a non-zero constant depending only on $\itPi_v$ and $\psi_v$ given by
\begin{align}\label{E:local constants1}
\begin{split}
C_{\psi_v}(\itPi_v) & =  \zeta_v(1)^{-1}\zeta_v(3)^{-1}\zeta_v(4)L(1,\itPi_v,{\rm Ad})^{-1}\\
&\times\begin{cases}
\zeta_v(2)^{-1}\zeta_v(4)^{-1}\cdot\displaystyle{\frac{\calZ_v(1, W_{\psi_v}, W_{\psi_v^{-1}}^\vee, \calF_{\psi_v}({\varphi_v}))}{Z_v \left(\frac{1}{2}, \phi_{v}, F_{\psi_v}(\varphi_{v})\right)}} & \mbox{ if $v$ is finite},\\
\displaystyle{\frac{\calZ_v(1, W_{\psi_v}, \itPi_v^\vee({\rm diag}(-1,-1,1,1))W_{\psi_v^{-1}}^\vee, \calF_{\psi_v}({\varphi_v}))}{Z_v \left(\frac{1}{2}, \phi_{v}, F_{\psi_v}(\varphi_{v})\right)} }& \mbox{ if $v$ is real},
\end{cases}
\end{split}
\end{align}
where $\varphi_v \in S(V_{3,3}^4(\F_v))$ is any Schwartz function such that $Z_v\left(\frac{1}{2},\phi_v,F_{\psi_v}(\varphi_v)\right) \neq 0$. Moreover, if $v$ is a finite place such that $\itPi_v$ and $\psi_v$ are unramified, then we have $C_{\psi_v}(\itPi_v)=1$.
\end{prop}

\begin{proof}
The assertion can be proved by proceeding exactly as in the proof of \cite[Proposition 5.4]{CI2019} with Lemmas 4.2, 4.4, and 5.3 in \cite{CI2019} replaced by Lemmas \ref{L:ab doubling}-(1), \ref{L:ab2}-(1), and \ref{lem:nonvanish}, respectively. Note that the constant $C_{\psi_v}(\itPi_v)$ is well-defined, since $L(s,\itPi_v,{\rm Ad})$ is holomorphic and non-zero at $s=1$ by the genericity of $\itPi_v$. Also note that the factors $\frak{D}^{-13}$ and $\zeta_v(2)^{-1}\zeta_v(4)^{-1}$ for finite places $v$ comparing with \cite[Proposition 5.4]{CI2019} are due to different normalization of Haar measures used to define the local zeta integrals $Z_v$ and $\mathcal{Z}_v$ for finite places $v$, and the partial Fourier transform $\hat{\varphi}_v$.
\end{proof}

\section{Proof of Main Theorem}\label{S:proof of thm}

We keep the notation of \S\,\ref{S:Petersson norm}.
Assume further that $\itPi$ is motivic, that is, there exists ${\sf w}\in\Z$ such that $|\omega_\itPi| = |\mbox{ }|_{\A}^{\sf w}$ and 
\[
\lambda_{1,v} - \lambda_{2,v} \equiv {\sf w}\,({\rm mod}\,2)
\]
for all real places $v$.
In other words, we have $\omega_{\itPi_v} = {\rm sgn}^{\sf w}|\mbox{ }|_v^{{\sf w}}$ for all real places $v$.
In the following lemma, we show that being globally generic is an arithmetic property of an irreducible cuspidal automorphic representation of $\GSp_4(\A)$ of motivic discrete type at the archimedean places.

\begin{lemma}\label{L:globally generic}
Assume $\itPi$ is motivic.
For $\sigma \in {\rm Aut}(\C)$, the representation ${}^\sigma\!\itPi$ is an irreducible motivic globally generic cuspidal automorphic representation of $\GSp_4(\A)$.
\end{lemma}

\begin{proof}

Fix $\sigma \in {\rm Aut}(\C)$. We denote by $\itPsi = \bigotimes_v{\itPsi}_v$ the strong functorial lift of $\itPi$ to $\GL_4(\A)$. By \cite[Theorem 12.1]{GT2011}, the automorphic representation
$\itPsi \boxtimes \omega_\itPi$
is equal to the global theta lift of $\itPi$ from $\GSp_4(\A)$ to $\GSO_{3,3}(\A)$, where we identify $\GSO_{3,3}$ with
\[
(\GL_4 \times \GL_1)/ \{(a{\bf 1}_4,a^{-2})\,\vert\,a\in \GL_1\}.
\]
For a real place $v$, by the assumption on $\itPi_v$, we have 
\[
\itPsi_v \simeq {\rm Ind}_{P_{2,2}(\R)}^{\GL_4(\R)}(D(\lambda_{1,v}+ \lambda_{2,v})\boxtimes D(\lambda_{1,v}- \lambda_{2,v})) \otimes |\mbox{ }|_v^{{\sf w}/2},
\]
where $P_{2,2}$ is the standard maximal parabolic subgroup of $\GL_4$ of type $(2,2)$ and $D(\kappa)$ denotes the discrete series representation of $\GL_2(\R)$ with minimal weight $\kappa \geq 2$ and central character ${\rm sgn}^\kappa$. In particular, we see that $\itPsi$ is regular algebraic in the sense of Clozel \cite[\S\,1.2.3 and Definition 3.12]{Clozel1990}. Therefore, ${}^\sigma\!\itPsi = {}^\sigma\!\itPsi_\infty \otimes {}^\sigma\!\itPsi_f$ is automorphic by \cite[Th\'eor\`eme 3.13]{Clozel1990}, where ${}^\sigma\!\itPsi_f$ is the $\sigma$-conjugate of $\itPsi_f = \bigotimes_{v \nmid \infty}\itPsi_v$ and ${}^\sigma\!\itPsi_{\infty}$ is the representation of $\GL_4(\F_\infty)$ so that its $v$-component is equal to $\itPsi_{\sigma^{-1} \circ v}$. We claim that there exists an irreducible globally generic cuspidal automorphic representation $\itPi^\sharp = \bigotimes_v \itPi_v^\sharp$ of $\GSp_4(\A)$ with central character ${}^\sigma\!\omega_\itPi = \omega_{\itPi_\infty}\cdot{}^\sigma\!\omega_{\itPi_f}$ such that ${}^\sigma\!\itPsi$ is the strong functorial lift of $\itPi^\sharp$. Assume the claim holds. For any finite place $v$, ${}^\sigma\!\itPsi_v \boxtimes {}^\sigma\!\omega_{\itPi_v}$ is the local theta lift of both ${}^\sigma\!\itPi_v$ and $\itPi_v^\sharp$ from $\GSp_4(\F_v)$ to $\GSO_{3,3}(\F_v)$. Indeed, since $\itPsi_v \boxtimes \omega_{\itPi_v}$ is the local theta lift of $\itPi_v$ from $\GSp_4(\F_v)$ to $\GSO_{3,3}(\F_v)$, it follows from \cite[Proposition 1.9]{Roberts2001} and \cite[Proposition 5.7]{Morimoto2012} that ${}^\sigma\!\itPsi_v \boxtimes {}^\sigma\!\omega_{\itPi_v}$ 
is the local theta lift of ${}^\sigma\!\itPi_v$ from $\GSp_4(\F_v)$ to $\GSO_{3,3}(\F_v)$. Also note that ${}^\sigma\!\itPsi_\infty \boxtimes \omega_{\itPi_\infty}$ is the local theta lift of both ${}^\sigma\!\itPi_\infty$ and $\itPi_\infty^\sharp$ from $\GSp_4(\F_\infty)$ to $\GSO_{3,3}(\F_\infty)$ (cf.\,\cite{Paul2005}). By the Howe duality principle, we have
\[
{}^\sigma\!\itPi_v \simeq \itPi_v^\sharp,\quad {}^\sigma\!\itPi_\infty \simeq \itPi_\infty^\sharp
\]
for all finite places $v$. 
We conclude that $\itPi^\sigma \simeq \itPi^\sharp$ is an irreducible globally generic cuspidal automorphic representation of $\GSp_4(\A)$.


It remains to verify the claim. Firstly we consider the case when $\itPi$ is stable. In this case, ${}^\sigma\!\itPsi$ is cuspidal by \cite[Theorem 3.13]{Clozel1990} and $L(s,\itPsi,\wedge^2\otimes \omega_\itPi^{-1})$ has a pole at $s=1$. By \cite[Theorem 12.1]{GT2011}, the claim holds if and only if $L(s,{}^\sigma\!\itPsi,\wedge^2\otimes {}^\sigma\!\omega_\itPi^{-1})$ has a pole at $s=1$. The last assertion was proved by Gan in \cite[Theorem 3.6.2]{GR2014} based on the result on functorial lifts from ${\rm GSpin}_4(\A)$ to $\GL_4(\A)$. Now we assume $\itPi$ is endoscopic. In this case, we write $\itPsi = \tau_1 \boxplus \tau_2$ for some non-isomorphic irreducible cuspidal automorphic representations $\tau_1$ and $\tau_2$ of $\GL_2(\A)$ with equal central character $\omega_\itPi$ such that
\[
\tau_{1,v} = D(\lambda_{1,v}+\lambda_{2,v})\otimes|\mbox{ }|_v^{{\sf w}/2},\quad \tau_{2,v} = D(\lambda_{1,v}-\lambda_{2,v})\otimes|\mbox{ }|_v^{{\sf w}/2}
\]
for all real places $v$. In particular, we see that $\tau_1$ and $\tau_2$ are regular algebraic.
Therefore, ${}^\sigma\!\itPsi = {}^\sigma\!\tau_1 \boxplus {}^\sigma\!\tau_2$ is an isobaric automorphic representation of $\GL_4(\A)$. We then take $\itPi^\sharp$ be the global theta lift of ${}^\sigma\!\tau_1 \boxtimes {}^\sigma\!\tau_2^\vee$ from ${\rm GSO}_{2,2}(\A)$ to $\GSp_4(\A)$, where we identify ${\rm GSO}_{2,2}$ with 
\[
(\GL_2 \times \GL_2)/\{(a{\bf 1}_2,a{\bf 1}_2)\,\vert\,a\in\GL_1\}
\]
as in \S\,\ref{SS:global theta lifts}.
Then $\itPi^\sharp$ is globally generic cuspidal and ${}^\sigma\!\itPsi$ is the strong functorial lift of $\itPi^\sharp$ by \cite[Theorem 12.1]{GT2011}.
This completes the proof.
\end{proof}

\begin{lemma}\label{L:rationality field}
The rationality field $\Q(\itPi)$ of $\itPi$ is equal to the fixed field of $\left\{\sigma \in {\rm Aut}(\C) \, \vert \, {}^\sigma\!\itPi_f = \itPi_f \right\}$ and is a number field.
\end{lemma}

\begin{proof}
Let $\itPsi = \bigotimes_v{\itPsi}_v$ be the strong functorial lift of $\itPi$ to $\GL_4(\A)$. The rationality field $\Q(\itPsi)$ of $\itPsi$ is the fixed field of $\left\{\sigma \in {\rm Aut}(\C) \, \vert \, {}^\sigma\!\itPsi = \itPsi \right\}$.
By the result of Clozel \cite[Th\'eor\`eme 3.13]{Clozel1990}, $\Q(\itPsi)$ is a number field. It follows from the strong multiplicity one theorem for isobaric automorphic representations \cite[Theorem 4.4]{JacquetShalika1981a} that $\Q(\itPsi)$ is equal to the fixed field of $\left\{\sigma \in {\rm Aut}(\C) \, \vert \, {}^\sigma\!\itPsi_f = \itPsi_f \right\}$.

Let $\sigma \in {\rm Aut}(\C)$. For each place $v$, we have explained in the proof of Lemma \ref{L:globally generic} that ${}^\sigma\!{\itPsi}_v$ is the functorial lift of ${}^\sigma\!\itPi_v$ to $\GL_4(\F_v)$ via the local theta correspondence. It then follows from the Howe duality principle that ${}^\sigma\!{\itPsi}_v = {\itPsi}_v$ if and only if ${}^\sigma\!\itPi_v = \itPi_v$. Therefore, we conclude that $\Q(\itPi) = \Q(\itPsi)$ and is equal to the fixed field of $\left\{\sigma \in {\rm Aut}(\C) \, \vert \, {}^\sigma\!\itPi_f = \itPi_f \right\}$.
This completes the proof.
\end{proof}

In the following lemma, we prove the Galois equivariant property of the local adjoint $L$-functions. The argument is standard and we give a proof for convenience of the reader (cf.\,\cite[Proposition 3.17]{Raghuram2009} and \cite[Proposition 5.4]{Morimoto2012}).

\begin{lemma}\label{L:3.1}
Let $v$ be a finite place. For $\sigma \in {\rm Aut}(\C)$, we have
\[
\sigma L(1,\itPi_v,{\rm Ad}) = L(1,{}^\sigma\!\itPi_v,{\rm Ad}).
\]
\end{lemma}

\begin{proof}

For $n \geq 1$, let ${\rm Irr}(\GL_n)$ be the set of isomorphism classes of irreducible admissible representations of $\GL_n(\F_v)$ and $\Phi(\GL_n)$ the set of equivalence classes of admissible $n$-dimensional representation of the Weil--Deligne group $L_{\F_v}$ of $\F_v$. Let
\[
{\rm Irr}(\GL_n) \longrightarrow \Phi(\GL_n) , \quad \itPsi \longmapsto \phi_\itPsi
\]
be the local Langlands correspondence established in \cite{HT2001} and \cite{Henniart1999}. Let $\sigma \in {\rm Aut}(\C)$. By \cite[Lemme 4.6]{Clozel1990} and \cite[Propri\'et\'e 3, \S\,7]{Henniart2001}, we have 
\begin{align*}
{}^\sigma\! L(s+\tfrac{n-1}{2}, \phi_\itPsi) & = L(s+\tfrac{n-1}{2}, \phi_{{}^\sigma\!\itPsi}),\\
{}^\sigma\!\phi_{\itPsi} &= \phi_{{}^\sigma\!\itPsi} \otimes \chi_\sigma^{n-1}.
\end{align*}
Here we regard the $L$-functions as rational functions in $q_v^{-s}$ and the $\sigma$-linear action is defined as in \S\,\ref{S:notation}, and 
$\chi_\sigma$ is the quadratic character of $L_{\F_v}$ associated to the quadratic character $\sigma(|\mbox{ }|_v^{1/2})\cdot|\mbox{ }|_v^{-1/2}$ of $\F_v^\times$.

Let $\itPsi_v \in {\rm Irr}(\GL_4)$ be the local functorial lift of $\itPi_v$ to $\GL_4(\F_v)$. We see from the proof of Lemma \ref{L:globally generic} that ${}^\sigma\!\itPsi_v$ is the local functorial lift of ${}^\sigma\!\itPi_v$ to $\GL_4(\F_v)$. Let ${\rm Sym}^2 : \GL_4(\C) \rightarrow \GL_{10}(\C)$ be the symmetric square representation. It is easy to verify that
\begin{align*}
{}^\sigma\!({\rm Sym}^2 \circ \phi) & = {\rm Sym}^2 \circ {}^\sigma\!\phi,\\
{\rm Sym}^2 \circ (\phi \otimes \chi) &= ({\rm Sym}^2\circ \phi) \otimes \chi^2 
\end{align*}
for any $\phi \in \Phi(\GL_4)$ and character $\chi$ of $L_{\F_v}$. In particular, we have
\[
{}^\sigma\!({\rm Sym}^2 \circ \phi_{\itPsi_v}) = {\rm Sym}^2\circ \phi_{{}^\sigma\!\itPsi_v}.
\]
Therefore, we deduce that
\begin{align*}
{}^\sigma\! L(s,\itPi_v,{\rm Ad}) & = {}^\sigma\!L(s,{\rm Sym}^2\circ \phi_{\itPsi_v})\\
& = {}^\sigma\! L(s+\tfrac{1}{2}, ({\rm Sym}^2\circ \phi_{\itPsi_v}) \otimes |\mbox{ }|_v^{-1/2})\\
& = L(s+\tfrac{1}{2}, {}^\sigma\!({\rm Sym}^2\circ \phi_{\itPsi_v}) \otimes \sigma(|\mbox{ }|_v^{-1/2})\chi_\sigma)\\
& = L(s+\tfrac{1}{2}, ({\rm Sym}^2\circ \phi_{{}^\sigma\!\itPsi_v})\otimes |\mbox{ }|_v^{-1/2})\\
& = L(s, {\rm Sym}^2\circ \phi_{{}^\sigma\!\itPsi_v})\\
& = L(s,{}^\sigma\!\itPi_v,{\rm Ad}).
\end{align*}
We obtain the assertion by evaluating at $s=1$. This completes the proof.
\end{proof}

\begin{lemma}\label{L:3.2}
Let $\psi_v$ be a non-trivial additive character of $\F_v$ and $\varphi \in \mathcal{S}(V_{3,3}^4(\F_v))$.

(1) Let $a \in \F_v^\times$ so that $a \in \frak{o}_v^\times$ if $v$ is finite. We have
\[
F_{\psi_v^{a^2}}(\varphi) = |a|_v^{12}\cdot F_{\psi_v}(\varphi'), \quad \mathcal{F}_{\psi_v^{a^2}}(\varphi) = |a|_v^{-3s-9}\cdot \mathcal{F}_{\psi_v}(\varphi'),
\]
where
\[
\varphi ' = \omega_{\psi_v}\left( \bp a^{-1} {\bf 1}_4 & 0 \\ 0 & a {\bf 1}_4\ep,1\right)\varphi.
\]

(2) Assume $v$ is finite. For $\sigma \in {\rm Aut}(\C)$, we have
\[
{}^\sigma\!F_{\psi_v}(\varphi)(g,\tfrac{1}{2}) = F_{{}^\sigma\!\psi_v}({}^\sigma\!\varphi)(g,\tfrac{1}{2}),\quad {}^\sigma\!\mathcal{F}_{\psi_v}(\varphi)(g,1) = \mathcal{F}_{{}^\sigma\!\psi_v}({}^\sigma\!\varphi)(g,1).
\]
\end{lemma}

\begin{proof}
Note that we have
\begin{align}\label{E:Weil rep1}
\omega_{\psi_v^{a^2}}= \omega_{\psi_v}\left( \bp a{\bf 1}_4 & 0 \\ 0 & a^{-1}{\bf 1}_4\ep,1\right)\circ \omega_{\psi_v}\circ \omega_{\psi_v}\left( \bp a^{-1}{\bf 1}_4 & 0 \\ 0 & a{\bf 1}_4\ep ,1\right).
\end{align}
Therefore, one can easily verify that
\begin{align}\label{E:Weil rep2}
\hat{\omega}_{\psi_v^{a^2}}(g,h)\hat{\varphi}(u,v) = \hat{\omega}_{\psi_v}(g,h)\hat{\varphi}'(au,av)
\end{align}
for $(g,h) \in {\rm G}(\Sp_8 \times {\rm O}_{3,3})(\F_v)$ and $u,v \in {\rm M}_{3,4}(\F_v)$. The first assertion for $F$ follows immediately from (\ref{E:Weil rep1}).
By definition, we have $f_{\psi_v}^\circ = f_{\psi_v^{a^2}}^\circ$ if $v$ is real. The assumption $a \in \o_v^\times$ if $v$ is finite implies that $f_{\psi_v}^\circ = f_{\psi_v^{a^2}}^\circ$.
The first assertion for $\mathcal{F}$ thus follows from (\ref{E:Weil rep2}).

Assume $v$ is finite. Let $m$ be a positive integer and $\phi \in \mathcal{S}({\rm M}_{m,m}(\F_v))$. The integral
\[
Z(s,\phi) = \int_{\GL_m(\F_v)} \phi(a) |\det(a)|_v^{s}\,da
\]
is absolutely convergent for ${\rm Re}(s) > m-1$. Moreover, we can deduce from the proof of \cite[Proposition 1.1]{GJ1972} that $Z(s,\phi)$ defines a rational function in $q_v^{-s}$ and satisfies the Galois equivariant property
\begin{align}\label{E:GJ integral}
\sigma Z(s,\phi) = Z(s,{}^\sigma\!\phi).
\end{align}
By the explicit formula for the Weil representation, we have
\begin{align}\label{E:Weil rep3}
{}^\sigma\!(\omega_{\psi_v}(g,h)\varphi) = \omega_{{}^\sigma\!\psi_v}(g,h){}^\sigma\!\varphi,\quad {}^\sigma\!(\hat{\omega}_{\psi_v}(g,h)\hat{\varphi}) = \hat{\omega}_{{}^\sigma\!\psi_v}(g,h){}^\sigma\!\hat{\varphi}
\end{align}
for $(g,h) \in {\rm G}(\Sp_{8} \times {\rm O}(V_{3,3}))(\F_v)$.
Also note that ${}^\sigma\!f_{\psi_v}^\circ(g,1) = f_{{}^\sigma\!\psi_v}^\circ(g,1)$ by definition. The second assertion then follows immediately from (\ref{E:GJ integral}) and (\ref{E:Weil rep3}). This completes the proof.

\end{proof}

\begin{lemma}\label{L:3.3}
Let $\psi_v$ be a non-trivial additive character of $\F_v$ and $a \in \F_v^\times$.
Assume either $v$ is finite and $\itPi_v$ is unramified, or $v$ is real and $a>0$.
Then we have
\begin{align*}
C_{\psi_v^a}(\itPi_v) = |a|_v^{-8}\cdot C_{\psi_v}(\itPi_v). 
\end{align*}

\end{lemma}

\begin{proof}
First assume $v$ is finite and $\itPi_v$ is unramified. Let $F_v^\circ \in I_v(s)$ and $\mathcal{F}_v^\circ \in \mathcal{I}_v(s)$ be the $H(\o_v)$-invariant sections such that $F_v^\circ(1,s)=\mathcal{F}_v^\circ(1,s)=1$. 
Let $\varpi_v^{d_v}\o_v$ be the largest fractional ideal of $\F_v$ on which $\psi_v$ is trivial. For $a \in \F_v^\times$, define $\varphi=\varphi_a \in S(V_{3,3}^4(\F_v))$ by 
\[
\varphi \bp x \\ y \ep = \mathbb{I}_{{\rm M}_{3,4}(a^{-1}\varpi_v^{d_v}\o_v)}(x) \mathbb{I}_{{\rm M}_{3,4}(\o_v)}(y).
\]
Then the partial Fourier transform (\ref{Partial Fourier}) of $\varphi$ with respect to $\psi_v^a \circ {\rm tr}$ is equal to 
\[
\hat{\varphi} = |a|_v^{-12}q_v^{-12d_v}\cdot\mathbb{I}_{{\rm M}_{3,8}(\o_v)}.
\]
Note that
\[
\omega_{\psi_v^a}(k,k')\varphi = \varphi 
\]
for $(k,k') \in (H(\o_v) \times {\bf K}_v')\cap {\rm G}(\Sp_{8}\times {\rm O}_{3,3})(\F_v)$, where
\begin{align*}
{\bf K}_v' = \bp {\bf 1}_3 & 0 \\ 0 & a \varpi_v^{-d_v}{\bf 1}_3\ep {\rm GO}_{3,3}(\o_v)\bp {\bf 1}_3 & 0 \\ 0 & a^{-1}\varpi_v^{d_v} {\bf 1}_3\ep.
\end{align*} 
Therefore, we have
\[
F_{\psi_v^a}(\varphi) = F_{\psi_v^a}(\varphi)(1,\tfrac{1}{2})\cdot F_v^\circ,\quad \mathcal{F}_{\psi_v^a}(\varphi) = \mathcal{F}_{\psi_v^a}(\varphi)(1,s)\cdot \mathcal{F}_v^\circ.
\]
Note that $F_{\psi_v^a}(\varphi)(1,\tfrac{1}{2})=\varphi(0)=1$ and
\begin{align*}
\mathcal{F}_{\psi_v^a}(\varphi)(1,s) &= |a|_v^{-12}q_v^{-12d_v} \int_{\GL_3(\F_v)} \mathbb{I}_{{\rm M}_{3,3}(\o_v)}(t) |\det(t)|_v^{s+3}\,dt\\
&=|a|_v^{-12}q_v^{-12d_v}\zeta_v(s+1)\zeta_v(s+2)\zeta_v(s+3).
\end{align*}
On the other hand, a change of variables $g \mapsto ({\rm diag}(a^3,a^2,1,a),{\rm diag}(a^3,a^2,1,a))^{-1}g$ shows that
\begin{align*}
\mathcal{Z}_v(s,W_{\psi_v^a},W_{\psi_v^{-a}}^\vee,\mathcal{F}) & = |a|_v^{-3s/2+11/2} \mathcal{Z}_v(s,W_{\psi_v},W_{\psi_v^{-1}}^\vee,\mathcal{F}).
\end{align*}
for any holomorphic section $\mathcal{F}$ of $\mathcal{I}_v(s)$.
We conclude that
\begin{align*}
Z_v(s,\phi_v,F_{\psi_v^a}(\varphi)) &= Z_v (s,\phi_v,F_v^\circ),\\
\mathcal{Z}_v(s,W_{\psi_v^a},W_{\psi_v^{-a}}^\vee,\mathcal{F}_{\psi_v^a}(\varphi)) & =  |a|_v^{-3s/2-13/2}q_v^{-12d_v}\zeta_v(s+1)\zeta_v(s+2)\zeta_v(s+3)\mathcal{Z}_v(s,W_{\psi_v},W_{\psi_v^{-1}}^\vee,\mathcal{F}_v^\circ).
\end{align*}
The assertion then follows immediately. 

Now we assume $v$ is real and $a>0$. 
Put $\widetilde{W}_{\xi_v} = \itPi_v^\vee({\rm diag}(-1,-1,1,1))W_{\xi_v}^\vee$ for any additive character $\xi_v$ of $\F_v$.
Let $\varphi \in S(V_{3,3}^4(\F_v))$ such that $Z_v \left(\frac{1}{2}, \phi_v, F_{\psi_v^a}(\varphi)\right) \neq 0$. By Lemma \ref{L:3.2}-(1), we have
\[
F_{\psi_v^{a}}(\varphi) = |a|_v^{6}\cdot F_{\psi_v}(\varphi'),\quad \mathcal{F}_{\psi_v^{a}}(\varphi) = |a|_v^{-3s/2-9/2}\cdot \mathcal{F}_{\psi_v}(\varphi'),
\]
where 
\[
\varphi ' = \omega_{\psi_v}\left( \bp \sqrt{a}^{-1}\, {\bf 1}_4 & 0 \\ 0 & \sqrt{a}\, {\bf 1}_4\ep,1\right)\varphi.
\]
Similarly, a change of variables shows that
\begin{align*}
\mathcal{Z}_v(s,W_{\psi_v^a},\widetilde{W}_{\psi_v^{-a}},\mathcal{F}) & = |a|_v^{-3s/2+11/2} \mathcal{Z}_v(s,W_{\psi_v},\widetilde{W}_{\psi_v^{-1}},\mathcal{F})
\end{align*}
for any holomorphic section $\mathcal{F}$ of $\mathcal{I}_v(s)$. We conclude that
\begin{align*}
Z_v(s,\phi_v,F_{\psi_v^a}(\varphi)) &= |a|_v^6Z_v (s,\phi_v,F_{\psi_v}(\varphi')),\\
\mathcal{Z}_v(s,W_{\psi_v^a},\widetilde{W}_{\psi_v^{-a}},\mathcal{F}_{\psi_v^a}(\varphi)) & =  |a|_v^{-3s+1}\mathcal{Z}_v(s,W_{\psi_v},\widetilde{W}_{\psi_v^{-1}},\mathcal{F}_{\psi_v}(\varphi')).
\end{align*}
This completes the proof.

\end{proof}

\begin{lemma}\label{L:3.4}
Let $\psi$ be a non-trivial additive character of $\F\backslash \A$. 
Let $S$ be any finite set of finite places of $\F$ containing the divisors of ${\rm cond}(\itPi)$.
For $a \in  \bigcap_{v \in S}\o_v^\times \cap \F_{>0}^\times$, we have
\[
\prod_{v \in S}C_{\psi_v^a}(\itPi_v) = \prod_{v \in S}C_{\psi_v}(\itPi_v).
\]
Here $\F_{>0}^\times$ denotes the set of totally positive elements in $\F$.
\end{lemma}

\begin{proof}
By Proposition \ref{P:Petersson norm}, we have
\[
\prod_v C_{\psi_v^a}(\itPi_v) = \prod_v C_{\psi_v}(\itPi_v).
\]
On the other hand, since $a$ is totally positive, by Lemma \ref{L:3.3} we have
\[
\prod_{v \notin S}C_{\psi_v^a}(\itPi_v) =  \prod_{v \notin S} |a|_v^{-8}\cdot C_{\psi_v}(\itPi_v).
\]
Since $a \in \bigcap_{v\in S}\o_v^\times \cap \F^\times$, we have $\prod_{v \notin S}|a|_v = |a|_\A=1$. This completes the proof.

\end{proof}

In the following theorem we prove the Galois equivariance property of the product of local constant $C_{\psi_v}(\itPi_v)$ over finite places.

\begin{prop}\label{T:3.6}
Let $\psi$ be a non-trivial additive character of $\F\backslash \A$. For $\sigma \in {\rm Aut}(\C)$, we have
\[
\sigma\left(\prod_{v \nmid \infty}  C_{\psi_v}(\itPi_v)\right) = \prod_{v \nmid \infty} C_{\psi_v}({}^\sigma\!\itPi_v).
\]
\end{prop}

\begin{proof}
Let $v$ be a place of $\F$ and $\xi_v$ a non-trivial additive character of $\F_v$. For $\varphi_v \in S(V_{3,3}^4(\F_v))$ such that $Z_v(\tfrac{1}{2},\phi_v,F_{\xi_v}(\varphi_v)) \neq 0$, define 
\begin{align*}
C_{\xi_v}(\itPi_v,\varphi_v)  &=  \zeta_v(1)^{-1}\zeta_v(3)^{-1}\zeta_v(4)L(1,\itPi_v,{\rm Ad})^{-1}\\
&\times\begin{cases}
\zeta_v(2)^{-1}\zeta_v(4)^{-1}\cdot\displaystyle{\frac{\calZ_v(1, W_{\xi_v}, W_{\xi_v^{-1}}^\vee, \calF_{\xi_v}({\varphi_v}))}{Z_v \left(\frac{1}{2}, \phi_{v}, F_{\xi_v}(\varphi_{v})\right)}} & \mbox{ if $v$ is finite},\\
\displaystyle{\frac{\calZ_v(1, W_{\xi_v}, \itPi_v^\vee({\rm diag}(-1,-1,1,1))W_{\xi_v^{-1}}^\vee, \calF_{\xi_v}({\varphi_v}))}{Z_v \left(\frac{1}{2}, \phi_{v}, F_{\xi_v}(\varphi_{v})\right)} }& \mbox{ if $v$ is real}.
\end{cases}
\end{align*}
When $\xi_v$ is the local component at $v$ of a non-trivial additive character of $\F \backslash \A$, it follows from Proposition \ref{P:Petersson norm} that $C_{\xi_v}(\itPi_v,\varphi_v)=C_{\xi_v}(\itPi_v)$ is the non-zero constant in (\ref{E:local constants1}) and does not depend on the choice of $\varphi_v$. 

Let $\sigma \in {\rm Aut}(\C)$. 
Let $u \in \widehat{\Z}^\times$ be the unique element such that $\sigma(\psi(x)) = \psi(ux)$ for $x \in \A_f$ and any non-trivial additive character $\psi$ of $\F\backslash \A$.
For each finite place $v$ lying over a rational prime $p$ and any non-trivial additive character $\xi_v$ of $\F_v$, define the $\sigma$-linear isomorphism
\begin{align*}
t_{\sigma,v} : \mathcal{W}(\itPi_v,\xi_{U,v}) &\longrightarrow \mathcal{W}({}^\sigma\!\itPi_v,\xi_{U,v}),\\
t_{\sigma,v} W (g) & = {}^\sigma W({\rm diag}(u_p^{-3},u_p^{-2},1,u_p^{-1})g).
\end{align*}
Note that ${}^\sigma\!\xi_v = \xi_v^{u_p}$.
Then $t_{\sigma, v}W_{\xi_v}$ is the normalized Whittaker newform of ${}^\sigma\!\itPi_v$ with respect to $\xi_{U,v}$ (cf.\,\S\,\ref {SS:2.3}
). 
We define $\sigma$-linear isomorphism $t_{\sigma,v} : \mathcal{W}(\itPi_v^\vee,\xi_{U,v}) \rightarrow \mathcal{W}({}^\sigma\!\itPi_v^\vee,\xi_{U,v})$ in the same way.
By the Chinese remainder theorem, there exists $a \in \bigcap_{p\mid {\rm N}_{\F/\Q}({\rm cond}(\itPi)\cdot {\rm cond}(\psi))}\Z_p^\times \cap \Q_{>0}^\times$ such that 
\[
au_p = s_p^2
\] for some $s_p \in \Z_p^\times$ for $p \mid {\rm N}_{\F/\Q}({\rm cond}(\itPi)\cdot {\rm cond}(\psi))$.  Let $v \mid {\rm cond}(\itPi)\cdot {\rm cond}(\psi)$ lying over a prime $p$. We have
\begin{align*}
\frac{\sigma C_{\psi_v^a}(\itPi_v)}{C_{\psi_v^{au_p}}({}^\sigma\!\itPi_v, {}^\sigma\!\varphi_v)}&=\frac{L(1,{}^\sigma\!\itPi_v,{\rm Ad})}{\sigma L(1,\itPi_v,{\rm Ad})}\cdot \frac{\sigma \calZ_v(1, W_{\psi_v^a}, W_{\psi_v^{-a}}^\vee, \calF_{\psi_v^a}({\varphi_v}))}{\calZ_v(1, t_{\sigma,v}W_{\psi_v^{au_p}}, t_{\sigma,v}W_{\psi_v^{-au_p}}^\vee, \calF_{\psi_v^{au_p}}({}^\sigma\!\varphi_v))}\\
&\times \frac{Z_v \left(\frac{1}{2}, {}^\sigma\!\phi_{v}, F_{\psi_v^{au_p}}({}^\sigma\!\varphi_{v})\right)}{\sigma Z_v \left(\frac{1}{2}, \phi_{v}, F_{\psi_v^a}(\varphi_v)\right)}.
\end{align*}
Note that ${}^\sigma W_{\psi_v^a} = t_{\sigma,v}W_{\psi_v^{au_p}}$ by definition.
It follows from Lemmas \ref{L:ab doubling}-(2), \ref{L:ab2}-(2), \ref{L:3.1}, and \ref{L:3.2}-(2) that the above ratio is equal to $1$. Similarly, we also have
\[
C_{\psi_v^{au_p}}({}^\sigma\!\itPi_v, {}^\sigma\!\varphi_v) = C_{\psi_v}({}^\sigma\!\itPi_v, \varphi_v') = C_{\psi_v}({}^\sigma\!\itPi_v),
\]
where
\[
\varphi_v' = \omega_{\psi_v}\left( \bp s_p^{-1} {\bf 1}_4 & 0 \\ 0 & s_p {\bf 1}_4\ep,1\right){}^\sigma\!\varphi_v.
\]
Indeed, since $au_p \in (\Z_p^{\times})^2$, by Lemma \ref{L:3.2}-(1) we have
\[
\frac{C_{\psi_v^{au_p}}({}^\sigma\!\itPi_v, {}^\sigma\!\varphi_v)}{C_{\psi_v}({}^\sigma\!\itPi_v, \varphi_v')} = \frac{\calZ_v(1, t_{\sigma,v}W_{\psi_v^{au_p}}, t_{\sigma,v}W_{\psi_v^{-au_p}}^\vee, \calF_{\psi_v}(\varphi_v'))}{\calZ_v(1, t_{\sigma,v} W_{\psi_v}, t_{\sigma,v} W_{\psi_v^{-1}}^\vee, \calF_{\psi_v}(\varphi_v'))}.
\]
Note that 
\[
t_{\sigma,v}W_{\psi_v^{au_p}} (g) = t_{\sigma,v} W_{\psi_v}({\rm diag}(s_p^2,1,s_p^{-4},s_p^{-2})g)
\]
and
\[
\calF_{\psi_v}(\varphi_v')(\eta({\rm diag}(s_p^2,1,s_p^{-4},s_p^{-2}),{\rm diag}(s_p^2,1,s_p^{-4},s_p^{-2}))g,s) = \calF_{\psi_v}(\varphi_v')(\eta g,s)
\]
for $g \in {\bf G}(\F_v)$. We see that the above ratio is also equal to $1$. We conclude that 
\[
\sigma \left(\prod_{v \mid {\rm cond}(\itPi)\cdot {\rm cond}(\psi)} C_{\psi_v^a}(\itPi_v)\right) = \prod_{v \mid {\rm cond}(\itPi)\cdot {\rm cond}(\psi)} C_{\psi_v}({}^\sigma\!\itPi_v).
\]
It then follows from Theorem \ref{L:3.4} that
\[
\sigma \left(\prod_{v \mid {\rm cond}(\itPi)\cdot {\rm cond}(\psi)} C_{\psi_v}(\itPi_v)\right) = \prod_{v \mid {\rm cond}(\itPi)\cdot {\rm cond}(\psi)} C_{\psi_v}({}^\sigma\!\itPi_v).
\]
Finally, note that for $v \nmid \infty\cdot{\rm cond}(\itPi)\cdot {\rm cond}(\psi)$, we have $C_{\psi_v}(\itPi_v) = C_{\psi_v}({}^\sigma\!\itPi_v)=1$. This completes the proof.
\end{proof}

By Propositions \ref{P:Petersson norm} and \ref{T:3.6}, our main result Theorem \ref{T:main} then follows from the following result.  

\begin{thm}
Let $v$ be a real place and $\psi_v$ the standard additive character of $\F_v$. We have
\[
C_{\psi_v}(\itPi_v) \in \Q^\times.
\]
\end{thm}

\begin{proof}
Applying \cite[Theorem 5.7]{Shin2012} to $\GSO_{2,2}(\A_\Q)$, there exist cohomological irreducible cuspidal automorphic representations $\tau_1$ and $\tau_2$ of $\GL_2(\A_\Q)$ with central character inverse to each other such that 
\[
\tau_{1,\infty} = D(\lambda_{1,v}-\lambda_{2,v}) \otimes |\mbox{ }|^{{\sf w}/2},\quad \tau_{2,\infty} = D(\lambda_{1,v}+\lambda_{2,v}) \otimes |\mbox{ }|^{-{\sf w}/2}
\]
and $\tau_1 \neq \tau_2^\vee$. We regard $\tau_1 \boxtimes \tau_2$ as an irreducible cuspidal automorphic representation of $\GSO_{2,2}(\A_\Q)$ and consider its global theta lift $\itPi'=\theta(\tau_1 \boxtimes \tau_2)$ to $\GSp_4(\A_\Q)$.
Then $\itPi'$ is an irreducible motivic globally generic cuspidal automorphic representation of $\GSp_4(\A_\Q)$ such that 
\[
\itPi'_\infty = \itPi_v.
\]
Let $\psi_0$ be the standard additive character of $\Q \backslash \A_\Q$. Thus $\psi_{\infty,0} = \psi_v$. 
By Propositions \ref{P:Petersson norm} and \ref{T:3.6} applied to $\psi_0$ and $\itPi'$, we have
\[
\sigma\left(\frac{L(1,\itPi',{\rm Ad})}{\zeta_\Q(2)\zeta_\Q(4)\cdot \Vert f_{\itPi'} \Vert}\cdot C_{\psi_{\infty,0}}(\itPi_\infty')\right) = \frac{L(1,{}^\sigma\!\itPi',{\rm Ad})}{\zeta_\Q(2)\zeta_\Q(4)\cdot \Vert f_{{}^\sigma\!\itPi'} \Vert}\cdot C_{\psi_{\infty,0}}(\itPi_\infty')
\]
for all $\sigma \in {\rm Aut}(\C)$.
On the other hand, we proved in Theorem \ref{T:endoscopic case} below that Theorem \ref{T:main} holds for $\itPi'$.
It follows that
\[
\sigma C_{\psi_{\infty,0}}(\itPi_\infty') = C_{\psi_{\infty,0}}(\itPi_\infty')
\]
for all $\sigma \in {\rm Aut}(\C)$.
Hence $C_{\psi_v}(\itPi_v) = C_{\psi_{\infty,0}}(\itPi_\infty') \in \Q$.
This completes the proof.
\end{proof}

\section{Petersson norms of endoscopic lifts}\label{S:endoscopic}

The purpose of this section is to prove Theorem \ref{T:main} for endoscopic lifts based on Rallis inner product formula. For simplicity, we assume $\F=\Q$ in this section.

\subsection{Cohomological cusp forms on $\GL_2$}

Let $\tau = \bigotimes_v \tau_v$ be a cohomological irreducible cuspidal automorphic representation of $\GL_2(\A)$. There exist $\kappa \in \Z_{\geq 2}$ and ${\sf w}\in\Z$ with $\kappa \equiv {\sf w}\,({\rm mod}\,2)$ such that 
\[
\tau_\infty = D(\kappa) \otimes |\mbox{ }|_\infty^{{\sf w}/2}.
\] 
Here $D(\kappa)$ denotes the discrete series representation of $\GL_2(\R)$ with minimal weight $\kappa$ and central character ${\rm sgn}^\kappa$.
Let $\tau^+$ (resp.\,$\tau^-$) be the space of holomorphic (resp.\,anti-holomorphic) cusp forms in $\tau$.
For a non-trivial additive character $\psi_v$ of $\Q_v$, let $\mathcal{W}(\tau_v,\psi_v)$ be the space of Whittaker functions of $\tau_v$ with respect to $\psi_v$.
When $v=p$ is finite, for $\sigma \in {\rm Aut}(\C)$, we define the $\sigma$-linear isomorphism 
\begin{align*}
t_{\sigma,p} : \mathcal{W}(\tau_p,\psi_p) &\longrightarrow \mathcal{W}({}^\sigma\!\tau_p,\psi_p),\\
t_{\sigma,p} W (g) & = {}^\sigma W({\rm diag}(u_p^{-1},1)g),
\end{align*}
where $u_p \in \Z_p^\times$ is the unique element such that $\sigma(\psi_p(x)) = \psi_p(u_px)$ for $x \in \Q_p$.
When $v=\infty$ and $\psi_\infty = \psi_{\infty,0}^a$, let $W_{(\pm\kappa;\,{\sf w}),\psi_\infty}\in\mathcal{W}(\tau_\infty,\psi_\infty)$ be the weight $\pm\kappa$ Whittaker function given by
\[
W_{(\pm\kappa;\,{\sf w}),\psi_\infty}(z{\bf n}(x){\bf a}(y)k_\theta) = z^{\sf w} (\pm a y)^{(\kappa+{\sf w})/2}e^{2\pi\sqrt{-1}\,a(x\pm\sqrt{-1}\,y)\pm \sqrt{-1}\,\kappa\theta}\cdot\mathbb{I}_{\R_{>0}^\times}(\pm a y)
\]
for $x \in \R$, $y,z \in \R^\times$, and $k_{\theta} = \bp \cos \theta & \sin \theta \\ -\sin \theta & \cos \theta \ep \in {\rm SO}(2)$.
Let $\psi$ be a non-trivial additive character of $\Q\backslash \A$.
For $f \in \tau$, let $W_{f,\psi}$ be the Whittaker function of $f$ with respect to $\psi$ defined by
\[
W_{f,\psi}(g) = \int_{\Q\backslash\A}f({\bf n}(x)g)\overline{\psi(x)}\,dx^{\rm Tam}.
\]
Here $dx^{\rm Tam}$ is the Tamagawa measure on $\A$.
For  $f \in \tau^\pm$, let $W_{f,\psi}^{(\infty)}\in\bigotimes_p \mathcal{W}(\tau_p,\psi_p)$ be the unique Whittaker function such that
\[
W_{f,\psi} = W_{(\pm\kappa;\,{\sf w}),\psi_\infty}\cdot W_{f,\psi}^{(\infty)}.
\]
Then we obtain the $\GL_2(\A_f)$-equivariant isomorphism 
\[
\tau^\pm \longrightarrow \bigotimes_p \mathcal{W}(\tau_p,\psi_p),\quad f\longmapsto W_{f,\psi}^{(\infty)}.
\]
For $\sigma \in {\rm Aut}(\C)$, it is well-known that the irreducible admissible representation 
\[
{}^\sigma\!\tau = \tau_\infty \otimes {}^\sigma\!\tau_f
\]
of $\GL_2(\A)$ is automorphic and cuspidal, where $\tau_f = \bigotimes_p \tau_p$.
Let 
\begin{align}\label{E:Galois conj. 1}
\tau^\pm \longrightarrow {}^\sigma\!\tau^\pm,\quad f \longmapsto {}^\sigma\!f
\end{align}
be the $\sigma$-linear $\GL_2(\A_f)$-equivariant isomorphism defined such that the diagram
\[
\begin{tikzcd}
\tau^\pm \arrow[r] \arrow[d] &{}^\sigma\!\tau^\pm\arrow[d]\\
\bigotimes_p \mathcal{W}(\tau_p,\psi_p)\arrow[r, "\bigotimes_pt_{\sigma,p}"] & \bigotimes_p \mathcal{W}({}^\sigma\!\tau_p,\psi_p)
\end{tikzcd}
\]
commutes. In other words, we have
\[
W_{{}^\sigma\!f,\psi}^{(\infty)} = \bigotimes_p t_{\sigma,p}W_{f,\psi}^{(\infty)}.
\]

\subsection{Arithmeticity of global theta lifting}\label{SS:global theta lifts}

Let $(V,Q)$ be the quadratic space over $\Q$ defined by $V={\rm M}_{2 , 2}$ and $Q[x]=\det(x)$. Let $\iota$ be the main involution on ${\rm M}_{2 , 2}$ defined by
\[\bp x_1&x_2\\x_3&x_4\ep^{\iota}=\bp x_4 & -x_2 \\ -x_3 & x_1\ep.\] The associated symmetric bilinear form is given by $(x,y) = {\rm tr}(xy^{\iota})$. 
We have an exact sequence
\begin{align}\label{E:exact seq.}
1 \longrightarrow   \GL_1 \stackrel{\Delta}{\longrightarrow} (\GL_2\times \GL_2)
\stackrel{\rho}{\longrightarrow} {\rm GSO}(V) \longrightarrow 1,
\end{align}
where $\Delta(a)=(a{\bf 1}_2,a{\bf 1}_2)$ and $\rho(h_1,h_2){x}=h_1{ x}h_2^{-1}$ for $a \in \GL_1, h_1, h_2 \in \GL_2$, and $x \in V$. For $h_1,h_2 \in \GL_2$, we write $\rho(h_1,h_2) = [h_1,h_2]$. Note that $\nu([h_1,h_2])=\det(h_1h_2^{-1})$. 
For a non-trivial additive character $\psi_v$ of $\Q_v$, we write $\omega_{\psi_v} = \omega_{\psi_v,V,2}$ for the Weil representation of $\Sp_4(\Q_v)\times{\rm O}(V)(\Q_v)$ on $\mathcal{S}(V^2(\Q_v))$ with respect to $\psi_v$.
Let $\psi$ be a non-trivial additive character of $\Q\backslash\A$.
For $\varphi \in S(V^2(\A))$, the theta function associated to $\varphi$ with respect to $\psi$ is defined by
\[\Theta_\psi(g,h;\varphi) = \sum_{x \in V^2(\Q)}\omega_\psi(g,h)\varphi(x)\]
for $(g,h) \in {\rm G}(\Sp_4 \times {\rm O}(V))(\A)$. Let $f$ be a cusp form on ${\rm GSO}(V)(\A)$ and let $\varphi \in S(V^2(\A))$. For $g \in \GSp_4(\A)$, choose $h \in {\rm GSO}(V)(\A)$ such that $\nu(h)=\nu(g)$, and put
\[
\theta_\psi(f,\varphi)(g) = \int_{{\rm SO}(V)(\Q)\backslash {\rm SO}(V)(\A)}f(h_1h)\Theta_\psi(g,h_1h;\varphi)\,dh_1^{\rm Tam}.
\]
Here $dh_1^{\rm Tam}$ is the Tamagawa measure on ${\rm SO}(V)(\A)$.
Then $\theta_\psi(f,\varphi)$ is an automorphic form on $\GSp_4(\A)$.

Let $\tau_1$ and $\tau_2$ be cohomological irreducible cuspidal automorphic representations with central character inverse to each other. There exist $\kappa_1,\kappa_2 \in \Z_{\geq 2}$ and ${\sf w}\in\Z$ with $\kappa_1 \equiv \kappa_2 \equiv {\sf w}\,({\rm mod}\,2)$ such that 
\[
\tau_{1,\infty} = D(\kappa_1) \otimes |\mbox{ }|_\infty^{{\sf w}/2},\quad \tau_{2,\infty} = D(\kappa_2) \otimes |\mbox{ }|_\infty^{-{\sf w}/2}.
\] 
We assume $\kappa_1 \geq \kappa_2$ and regard $\tau_1\boxtimes\tau_2$ as an automorphic representation of $\GSO(V)(\A)$ via the exact sequence (\ref{E:exact seq.}).
We assume further that
\[
\tau_1 \neq \tau_2^\vee.
\]
Then the global theta lift
\[
\theta(\tau_1\boxtimes\tau_2) = \left\{\theta_\psi(f,\varphi)\,\left\vert\,f\in\tau_1\boxtimes\tau_2,\,\varphi\in S(V^2(\A))\right.\right\}
\]
is an irreducible motivic globally generic cuspidal automorphic representation of $\GSp_4(\A)$ (cf.\,\cite[Theorem 12.1]{GT2011}).
As the notation suggests, the global theta lift does not depend on the choice of $\psi$ (cf.\,\cite[Proposition 1.9]{Roberts2001}).
Write 
\[
\itPi = \theta(\tau_1\boxtimes\tau_2).
\]
Note that $\omega_\itPi = \omega_{\tau_1} = \omega_{\tau_2}^{-1}$. Moreover, $\itPi_\infty$ is a generic discrete sereis representation of $\GSp_4(\R)$ with 
\[
\itPi_\infty\vert_{\Sp_4(\R)} = D_{(\lambda_1,\lambda_2)}\oplus D_{(-\lambda_2,-\lambda_1)}
\]
and
\[
(\lambda_1,\lambda_2) = (\tfrac{\kappa_1+\kappa_2}{2},\tfrac{\kappa_2-\kappa_1}{2}).
\]
Let $\itPi_{\rm mot}$ be the space of cusp forms in $\itPi$ such that the archimedean component is a lowest weight vector of the minimal ${\rm U}(2)$-type of $D_{(-\lambda_2,-\lambda_1)}$.
For $\sigma \in {\rm Aut}(\C)$ and $\psi_p$ a non-trivial additive character of $\Q_p$, we define the $\sigma$-linear isomorphism
\begin{align*}
t_{\sigma,p} : \mathcal{W}(\itPi_p,\psi_{U,p}) &\longrightarrow \mathcal{W}({}^\sigma\!\itPi_p,\psi_{U,p}),\\
t_{\sigma,p} W (g) & = {}^\sigma W({\rm diag}(u_p^{-3},u_p^{-2},1,u_p^{-1})g),
\end{align*}
where $u_p \in \Z_p^\times$ is the unique element such that $\sigma(\psi_p(x)) = \psi_p(u_px)$ for $x \in \Q_p$.
For $v=\infty$ and $\psi_\infty$ the standard additive character of $\R$, let $W_{(\lambda_1,\lambda_2;\,{\sf w}),\psi_{U,\infty}} \in \mathcal{W}(\itPi_\infty,\psi_{U,\infty})$ be the lowest weight Whittaker function of the minimal ${\rm U}(2)$-type of $D_{(-\lambda_2,-\lambda_1)}$ with respect to $\psi_{U,\infty}$ normalized as in (\ref{E:normalization}).
For $a \in \R^\times$, we then define $W_{(\lambda_1,\lambda_2;\,{\sf w}),\psi_{U,\infty}^a} \in \mathcal{W}(\itPi_\infty,\psi_{U,\infty}^a)$ by
\[
W_{(\lambda_1,\lambda_2;\,{\sf w}),\psi_{U,\infty}^a}(g) = W_{(\lambda_1,\lambda_2;\,{\sf w}),\psi_{U,\infty}}({\rm diag}(a^3,a^2,1,a)g).
\]
Let $\psi$ be a non-trivial additive character of $\Q\backslash\A$.
For $f \in \itPi$, let $W_{f,\psi_U}$ be the Whittaker function of $f$ with respect to $\psi_U$ defined by
\[
W_{f,\psi_U}(g) = \int_{U(\Q)\backslash U(\A)}f(ug)\overline{\psi_U(u)}\,du^{\rm Tam}.
\]
Here $du^{\rm Tam}$ is the Tamagawa measure on $U(\A)$.
For  $f \in \itPi_{\rm mot}$, let $W_{f,\psi_U}^{(\infty)}\in\bigotimes_p \mathcal{W}(\itPi_p,\psi_{U,p})$ be the unique Whittaker function such that
\[
W_{f,\psi_U} = W_{(\lambda_1,\lambda_2;\,{\sf w}),\psi_{U,\infty}}\cdot W_{f,\psi_U}^{(\infty)}.
\]
Then we obtain the $\GSp_4(\A_f)$-equivariant isomorphism 
\[
\itPi_{\rm mot} \longrightarrow \bigotimes_p \mathcal{W}(\itPi_p,\psi_{U,p}),\quad f\longmapsto W_{f,\psi_U}^{(\infty)}.
\]
For $\sigma \in {\rm Aut}(\C)$, the irreducible admissible representation 
\[
{}^\sigma\!\itPi = \itPi_\infty \otimes {}^\sigma\!\itPi_f
\]
of $\GSp_4(\A)$ is automorphic, globally generic, and cuspidal by Lemma \ref{L:globally generic}, where $\itPi_f = \bigotimes_p \itPi_p$.
Let 
\begin{align}\label{E:Galois conj. 2}
\itPi_{\rm mot} \longrightarrow {}^\sigma\!\itPi_{\rm mot},\quad f \longmapsto {}^\sigma\!f
\end{align}
be the $\sigma$-linear $\GSp_4(\A_f)$-equivariant isomorphism defined such that the diagram
\[
\begin{tikzcd}
\itPi_{\rm mot} \arrow[r] \arrow[d] &{}^\sigma\!\itPi_{\rm mot}\arrow[d]\\
\bigotimes_p \mathcal{W}(\itPi_p,\psi_{U,p})\arrow[r, "\bigotimes_pt_{\sigma,p}"] & \bigotimes_p \mathcal{W}({}^\sigma\!\itPi_p,\psi_{U,p})
\end{tikzcd}
\]
commutes. In other words, we have
\begin{align}\label{E:Galois conj. 3}
W_{{}^\sigma\!f,\psi_U}^{(\infty)} = \bigotimes_p t_{\sigma,p}W_{f,\psi_U}^{(\infty)}.
\end{align}
The main result of this section is Proposition \ref{P:arithmeticity}. Roughly speaking, we show that global theta lifting commutes with the Galois actions (\ref{E:Galois conj. 1}) and (\ref{E:Galois conj. 2}).

We begin with the following formula for Whittaker functions of global theta lifts.
For each place $v$ of $\Q$, let $dh_{1,v}$ be the Haar measure on $\SO(V)(\Q_v)$ defined as in \cite[\S\,7.1]{CI2019}.

\begin{lemma}\label{L:Whittaker factorization}
Let $\psi$ be a non-trivial additive character of $\Q\backslash\A$.
Let $\varphi = \bigotimes_v\varphi_v\in S(V^2(\A))$ and $f = f_1\otimes f_2 \in \tau_1\boxtimes\tau_2$ with $W_{f_1,\psi} = \prod_vW_{1,v}$ and $W_{f_2,\psi} = \prod_vW_{2,v}$. Then we have
\[
W_{\theta_\psi(f,\varphi),\psi_U} = \zeta(2)^{-2}\prod_v \frak{W}_{\psi_v}(W_{1,v}\otimes W_{2,v},\varphi_v).
\]
Here
\[
\frak{W}_{\psi_v}(W_{1,v}\otimes W_{2,v},\varphi_v)(g_v) = \int_{\Delta N(\Q_v)\backslash \SO(V)(\Q_v)}(W_{1,v}\otimes W_{2,v})(h_{1,v}h_v)\,\omega_{\psi_v}(g_v,h_{1,v}h_v)\varphi_v({\bf x}_0,{\bf y}_0)\,d\overline{h}_{1,v}
\]
for $(g_v,h_v) \in {\rm G}(\Sp_4\times {\rm O}(V))(\Q_v)$ with $h_v \in \GSO(V)(\Q_v)$,
\[
 {\bf x}_0 =
 \begin{pmatrix}
  0 & -1 \\
  0 & 0
 \end{pmatrix},\quad {\bf y}_0={\bf a}(-1),
\]
$\Delta N = \{ [{\bf n}(x),{\bf n}(-x)] \in \SO(V) \mbox{ }\vert\mbox{ }x \in \mathbb{G}_a\}$, and $d\overline{h}_{1,v}$ is the quotient measure defined by the Haar measures on $\Delta N(\Q_v)$ and $\SO(V)(\Q_v)$.
\end{lemma}

\begin{proof}
This is \cite[Lemma 7.1]{CI2019}. The factor $\zeta(2)^{-2}$ is due to the comparison between quotient of Tamagawa measures on $\Delta N(\A)$ and $\SO(V)(\A)$ with $\prod_vd\overline{h}_{1,v}$.
\end{proof}

For $v=p$, let $\varphi_p^\circ \in S(V^2(\Q_p))$ defined by
\[
\varphi_p^\circ = \mathbb{I}_{V^2(\Z_p)}.
\]
Note that for $a \in \Q_p^\times$ and $\psi_p$ the standard additive character of $\Q_p$, we have
\begin{align}\label{E:spherical condition}
\omega_{\psi_p^a}(k,[k_1,k_2])\varphi_p^\circ = \varphi_p^\circ
\end{align}
for $k \in {\rm diag}(1,1,a,a)\GSp_4(\Z_p){\rm diag}(1,1,a,a)^{-1}$ and $(k_1,k_2) \in \GL_2(\Z_p)\times\GL_2(\Z_p)$ such that $\nu(k) = \det(k_1k_2^{-1})$.
When $\tau_{1,p}$ and $\tau_{2,p}$ are both unramified, let $W_{p,a}^\circ \in \mathcal{W}(\tau_{1,p},\psi_p^a)\otimes\mathcal{W}(\tau_{2,p},\psi_p^a)$ be the $\GL_2(\Z_p)\times\GL_2(\Z_p)$-invariant Whittaker function normalized so that
\[
W_{p,a}^\circ(({\bf a}(a^{-1}),{\bf a}(a^{-1})))=1.
\]
By (\ref{E:spherical condition}), the Whittaker function $\frak{W}_{\psi_p^a}(W_{p,a}^\circ,\varphi_p^\circ) \in \mathcal{W}(\itPi_p,\psi_{U,p}^a)$ is ${\rm diag}(1,1,a,a)\GSp_4(\Z_p){\rm diag}(1,1,a,a)^{-1}$-invariant.

\begin{lemma}\label{L:Whittaker unramified}
Let $\psi_p$ be the standard additive character of $\Q_p$ and $a \in \Q_p^\times$. 
Assume both $\tau_{1,p}$ and $\tau_{2,p}$ are unramified.
We have
\[
\frak{W}_{\psi_p^a}(W_{p,a}^\circ,\varphi_p^\circ)({\rm diag}(a^{-3},a^{-2},a^{-1},a^{-2})) = \omega_{\itPi_p}(a)^{-2}|a|_p^{-1}.
\]
\end{lemma}

\begin{proof}
The computation is similar to \cite[Lemma 7.3]{CI2019} and we leave the detail to the readers.
\end{proof}

For $v=\infty$, let $\varphi_\infty^\pm \in S(V^2(\R))$ defined by
\begin{align*}
\varphi_\infty^+(x,y) &= 2^{\lambda_1+4}(-\sqrt{-1}\,x_1-x_2-x_3+\sqrt{-1}\,x_4)^{\lambda_1}(y_1+\sqrt{-1}\,y_2-\sqrt{-1}\,y_3+y_4)^{-\lambda_2}e^{-\pi\,{\rm tr}(x{}^t \! x+y{}^t \! y)},\\
\varphi_\infty^-(x,y) &= 2^{\lambda_1+4}(\sqrt{-1}\,x_1-x_2+x_3+\sqrt{-1}\,x_4)^{\lambda_1}(-y_1+\sqrt{-1}\,y_2+\sqrt{-1}\,y_3+y_4)^{-\lambda_2}e^{-\pi\,{\rm tr}(x{}^t \! x+y{}^t \! y)}.
\end{align*}
Let $\psi_\infty$ be the standard additive character of $\R$, then we have
\[
\omega_{\psi_\infty}(1,[k_{\theta_1},k_{\theta_2}])\varphi_\infty^+ = e^{-\sqrt{-1}\,(\kappa_1\theta_1+\kappa_2\theta_2)}\varphi_\infty^+,\quad \omega_{\psi_\infty^{-1}}(1,[k_{\theta_1},k_{\theta_2}])\varphi_\infty^- = e^{-\sqrt{-1}\,(\kappa_1\theta_1-\kappa_2\theta_2)}\varphi_\infty^-
\]
for $k_{\theta_1},k_{\theta_2} \in {\rm SO}(2)$, and
\begin{align}\label{E:lowest weight condition}
\omega_{\psi_\infty^\pm}(Z,1)\cdot \varphi_\infty^\pm = -\lambda_1\cdot\varphi_\infty^\pm,\quad \omega_{\psi_\infty^\pm}(Z',1)\cdot \varphi_\infty^\pm = -\lambda_2\cdot\varphi_\infty^\pm,\quad \omega_{\psi_\infty^\pm}(N_{-},1)\cdot\varphi_\infty^\pm =0. 
\end{align}
Here $Z, Z', N_-$ are elements in $\frak{sp}_4(\R)\otimes_\R\C$ defined by
\[
Z=-\sqrt{-1} \bp 0&0&1&0 \\ 0&0&0&0 \\ -1&0&0&0\\ 0&0&0&0 \ep,\quad Z' = -\sqrt{-1} \bp 0&0&0&0 \\ 0&0&0&1 \\ 0&0&0&0\\ 0&-1&0&0 \ep,\quad N_- = \bp 0&1&0&\sqrt{-1} \\ -1&0&\sqrt{-1}&0 \\ 0&-\sqrt{-1}&0&1 \\ -\sqrt{-1}& 0 & -1 & 0\ep.
\]
Define $W_{\infty}^\pm \in \mathcal{W}(\tau_{1,\infty},\psi_\infty^\pm)\otimes\mathcal{W}(\tau_{2,\infty},\psi_\infty^\pm)$ by
\[
W_\infty^+ = W_{(\kappa_1;\,{\sf w}),\psi_\infty}\otimes W_{(\kappa_2;\,-{\sf w}),\psi_\infty},\quad W_{\infty}^- = W_{(\kappa_1;\,{\sf w}),\psi_\infty^{-1}}\otimes W_{(-\kappa_2;\,-{\sf w}),\psi_\infty^{-1}}.
\]
By (\ref{E:lowest weight condition}), the Whittaker function $\frak{W}_{\psi_\infty^\pm}(W_\infty^\pm,\varphi_\infty^\pm) \in \mathcal{W}(\itPi_\infty,\psi_{U,\infty}^\pm)$ is a lowest weight Whittaker function of the minimal ${\rm U}(2)$-type of $D_{(-\lambda_2,-\lambda_1)}$ with respect to $\psi_{U,\infty}^\pm$.

\begin{lemma}\label{L:Whittaker archimedean}
Let $\psi_\infty$ be the standard additive character of $\R$.
We have
\begin{align*}
\frak{W}_{\psi_\infty^\pm}(W_\infty^\pm,\varphi_\infty^\pm) = W_{(\lambda_1,\lambda_2;\,{\sf w}),\psi_{U,\infty}^\pm}.
\end{align*}
\end{lemma}

\begin{proof}
The assertion for $\frak{W}_{\psi_\infty}(W_\infty^+,\varphi_\infty^+)$ was proved in \cite[Lemma 7.7]{CI2019}. The computation for the other case is similar and we leave the detail to the readers.
\end{proof}

In the following lemma, we establish an explicit relation between $\mathcal{W}_{\psi_v}$ and $\mathcal{W}_{\psi_v^{a^2}}$.
\begin{lemma}\label{L:change characters}
Let $\psi_v$ be a non-trivial additive character of $\Q_v$ and $a \in \Q_v^\times$. For $\varphi \in \mathcal{S}(V^2(\Q_v))$ and $W \in \mathcal{W}(\tau_{1,v},\psi_v)\otimes\mathcal{W}(\tau_{2,v},\psi_v)$, we have
\begin{align*}
& \frak{W}_{\psi_v^{a^2}}(\ell(({\bf a}(a^2),{\bf a}(a^2)))W,\varphi)\left({\rm diag}(a^{-6},a^{-4},1,a^{-2})g\right)\\
& = \omega_{\itPi_v}(a)^{-3}|a|_v^{-2}\cdot\frak{W}_{\psi_v}\left(W,\omega_{\psi_v}\left(\bp a^{-1}{\bf 1}_2 & 0 \\ 0 & a {\bf 1}_2\ep,1\right)\varphi\right)(g) 
\end{align*}
for $g \in \GSp_4(\Q_v)$.
Here $\ell(({\bf a}(a^2),{\bf a}(a^2)))W \in \mathcal{W}(\tau_{1,v},\psi_v^{a^2})\otimes \mathcal{W}(\tau_{2,v},\psi_v^{a^2})$ is defined by
\[
\ell(({\bf a}(a^2),{\bf a}(a^2)))W(h) = W(({\bf a}(a^2),{\bf a}(a^2))h).
\]
\end{lemma}

\begin{proof}
First note that
\begin{align*}
\omega_{\psi_v^{a^2}} = \omega_{\psi_v}\left(\bp a{\bf 1}_2 & 0 \\ 0 & a^{-1} {\bf 1}_2\ep,1\right)\circ\omega_{\psi_v}\circ\omega_{\psi_v}\left(\bp a^{-1}{\bf 1}_2 & 0 \\ 0 & a {\bf 1}_2\ep,1\right)
\end{align*}
and
\begin{align*}
\omega_{\psi_v}\left({\rm diag}(b^{-1},1,b,1),[{\bf a}(b^{-1}),{\bf a}(b^{-1})]\right)\varphi'({\bf x}_0,{\bf y}_0) = |b|_v^{-2}\varphi'({\bf x}_0,{\bf y}_0)
\end{align*}
for any $\varphi' \in \mathcal{S}(V^2(\Q_v))$ and $b \in \Q_v^\times$.
Let $(g,h)\in{\rm G}(\Sp_4\times{\rm O}(V))(\Q_v)$ with $h \in \GSO(V)(\Q_v)$.
Thus we have
\begin{align*}
&\frak{W}_{\psi_v^{a^2}}(\ell(({\bf a}(a^2),{\bf a}(a^2)))W,\varphi)\left({\rm diag}(a^{-6},a^{-4},1,a^{-2})g\right)\\
& = \omega_{\itPi_v}(a)^{-3}\cdot\mathcal{W}_{\psi_v^{a^2}}(\ell(({\bf a}(a^2),{\bf a}(a^2)))W,\varphi)\left({\rm diag}(a^{-3},a^{-1},a^3,a)g\right)\\
&= \omega_{\itPi_v}(a)^{-3}|a|_v^2\int_{\Delta N(\Q_v)\backslash \SO(V)(\Q_v)} W(h_1h)\,\omega_{\psi_v^{a^2}}\left({\rm diag}(a^{-3},a^{-1},a^3,a)g,[{\bf a}(a^{-2}),{\bf a}(a^{-2})]h_1h\right)\varphi({\bf x}_0,{\bf y}_0)\,dh_1\\
&= \omega_{\itPi_v}(a)^{-3}|a|_v^2\\
&\times \int_{\Delta N(\Q_v)\backslash \SO(V)(\Q_v)} W(h_1h)\,\omega_{\psi_v}\left({\rm diag}(a^{-2},1,a^2,1)g\,{\rm diag}(a^{-1},a^{-1},a,a),[{\bf a}(a^{-2}),{\bf a}(a^{-2})]h_1h\right)\varphi({\bf x}_0,{\bf y}_0)\,dh_1\\
&= \omega_{\itPi_v}(a)^{-3}|a|_v^{-2}\cdot\frak{W}_{\psi_v}\left(W,\omega_{\psi_v}\left(\bp a^{-1}{\bf 1}_2 & 0 \\ 0 & a {\bf 1}_2\ep,1\right)\varphi\right)(g).
\end{align*}
Here the factor $|a|_v^2$ in the second equality is due to the change of variables $h_1 \mapsto [{\bf a}(a^{-2}),{\bf a}(a^{-2})]h_1$.
This completes the proof.
\end{proof}

Let $\psi$ be a non-trivial additive character of $\Q\backslash\A$ and $S$ a finite set of finite places of $\Q$. We write 
\[
\tau_{i,S} = \bigotimes_{p \in S}\tau_{i,p},\quad \itPi_S = \bigotimes_{p \in S}\itPi_p,\quad \Q_S=\prod_{p \in S}\Q_p,\quad |\mbox{ }|_S = \prod_{p \in S}|\mbox{ }|_p,\quad \psi_S = \bigotimes_{p \in S}\psi_p,\quad \omega_{\psi_S} = \bigotimes_{p \in S}\omega_{\psi_p}
\]
for $i=1,2$.
Assume $\psi$ is standard.
For $W_S \in \mathcal{W}(\tau_{1,S},\psi_S^\pm)\otimes\mathcal{W}(\tau_{2,S},\psi_S^\pm)$, let 
\[
f_{W_S} \in \tau_1^+\boxtimes\tau_2^\pm
\]
defined so that
\[
W_{f_{W_S},\psi^\pm\boxtimes\psi^\pm}^{(\infty)} = W_S \cdot \prod_{p \notin S}W_{p,\pm}^\circ.
\]
For $a \in \Q^\times$, by the definition of global Whittaker function, we easily see that
\[
W_{f_{W_S},\psi^{\pm a}\boxtimes\psi^{\pm a}}^{(\infty)} = \ell(({\bf a}(a),{\bf a}(a)))W_S\cdot \prod_{p \notin S}W_{p,\pm a}^\circ.
\]
For $\sigma \in {\rm Aut}(\C)$, let $u_\sigma \in \widehat{\Z}^\times$ be the unique element such that $\sigma(\psi(x)) = \psi(u_\sigma x)$ for all $x \in \A_f$.
We have the $\sigma$-linear $\SO(V)(\A_f)$-equivariant isomorphism
\begin{align}\label{E:Galois conj. 4}
\tau_1^+\boxtimes\tau_2^\pm \longrightarrow {}^\sigma\!\tau_1^+\boxtimes{}^\sigma\!\tau_2^\pm,\quad f \longmapsto {}^\sigma\!f
\end{align}
defined by the $\sigma$-linear isomorphisms in (\ref{E:Galois conj. 1}) for $\tau_1^+$ and $\tau_2^\pm$.

\begin{prop}\label{P:arithmeticity}
Let $\psi$ be the standard additive character of $\Q\backslash\A$ and $S$ a finite set of finite places of $\Q$ containing the prime divisors of ${\rm cond}(\itPi)$. Let $\sigma \in {\rm Aut}(\C)$ and $a \in \Q^\times_{>0}$. Assume $u_{\sigma,S} = at^2$
for some $t \in \Q_S^\times$.
For $\varphi_S \in S(V^2(\Q_S))$ and $W_S \in \mathcal{W}(\tau_{1,S},\psi_S^\pm)\otimes\mathcal{W}(\tau_{2,S},\psi_S^\pm)$, we have
\begin{align}\label{E:arithmeticity proof 0}
\begin{split}
&{}^\sigma\!\itPi^{S\cup\{\infty\}}({\rm diag}(1,1,a,a)){}^\sigma\!\theta_{\psi^\pm}\left(f_{W_S},\,\varphi_\infty^\pm\otimes\varphi_S \otimes (\otimes_{p \notin S}\varphi_p^\circ)\right)\\
& = \frac{\zeta(2)^2}{\sigma(\zeta(2))^2}\cdot \omega_{\itPi_\infty}(\sqrt{a})^{-1}\cdot{}^\sigma\!\omega_{\itPi_S}(a^2t^3)^{-1}\\
&\times \theta_{\psi^{\pm a}}\left({}^\sigma\!f_{W_S},\,\omega_{\psi_\infty^\pm}\left(\bp \sqrt{a}\,{\bf 1}_2&0 \\0 & \sqrt{a}^{-1}\,{\bf 1}_2\ep,1\right)\varphi_\infty^\pm\otimes \omega_{\psi_S^\pm}\left(\bp t^{-1}{\bf 1}_2 & 0 \\ 0 & t{\bf 1}_2 \ep,1\right){}^\sigma\!\varphi_S \otimes(\otimes_{p \notin S}\varphi_p^\circ) \right).
\end{split}
\end{align}
Here ${}^\sigma\!\itPi^{S\cup\{\infty\}} = \bigotimes_{v \notin S\cup\{\infty\}}{}^\sigma\!\itPi_v$ and the Galois action on the cusp form on the left-hand side (resp.\,right-hand side) is defined in (\ref{E:Galois conj. 2}) (resp.\,(\ref{E:Galois conj. 4})).
\end{prop}

\begin{proof}
To prove the assertion, it suffices to show that the Whittaker functions of both sides with respect to $\psi_U^\pm$ are equal.
By Lemmas \ref{L:Whittaker archimedean} and \ref{L:change characters}, we have
\begin{align}\label{E:arithmeticity proof 1}
\begin{split}
&\frak{W}_{\psi_\infty^{\pm b}}\left(\ell(({\bf a}(b),{\bf a}(b)))W_\infty^\pm,\,\omega_{\psi_\infty^\pm}\left(\bp \sqrt{b}\,{\bf 1}_2&0 \\0 & \sqrt{b}^{-1}\,{\bf 1}_2\ep,1\right)\varphi_\infty^\pm\right)({\rm diag}(b^{-3},b^{-2},1,b^{-1})g)\\
&=\omega_{\itPi_\infty}(\sqrt{b})^{-3}|b|_v^{-1}\cdot W_{(\lambda_1,\lambda_2;\,{\sf w}),\psi_{U,\infty}^\pm}
\end{split}
\end{align}
for all $b>0$.
We have
\[
W_{f_{W_S},\psi^\pm\boxtimes\psi^\pm} = W_\infty^\pm\cdot W_S\cdot \prod_{p\notin S}W_{p,\pm}^\circ.
\]
By Lemma \ref{L:Whittaker factorization} for the additive character $\psi^\pm$, (\ref{E:Galois conj. 3}), and (\ref{E:arithmeticity proof 1}) with $b=1$, we see that the Whittaker function of the global theta lift on the left-hand side of (\ref{E:arithmeticity proof 0}) with respect to $\psi_U^\pm$ is equal to 
\begin{align}\label{E:arithmeticity proof 2}
\begin{split}
\sigma(\zeta(2))^{-2}\cdot W_{(\lambda_1,\lambda_2;\,{\sf w}),\psi_{U,\infty}^\pm}(g_\infty)\cdot t_{\sigma,S}\,\frak{W}_{\psi_S^\pm}\left(W_S,\,\varphi_S\right)(g_S)\cdot\prod_{p \notin S}t_{\sigma,p}\,\frak{W}_{\psi_p^\pm}(W_{p,\pm a}^\circ,\varphi_p^\circ)(g_p\cdot{\rm diag}(1,1,a,a))
\end{split}
\end{align}
for $g \in \GSp_4(\A)$.
Let $p \notin S$. By Lemma \ref{L:Whittaker unramified}, we have
\[
\frak{W}_{\psi_p^\pm}(W_{p,\pm}^\circ,\varphi_p^\circ)(g_p\cdot{\rm diag}(1,1,a,a))= \omega_{\itPi_p}(a)^2|a|_p\cdot\frak{W}_{\psi_p^{\pm a}}(W_{p,\pm a}^\circ,\varphi_p^\circ)({\rm diag}(a^{-3},a^{-2},1,a^{-1})g_p)
\]
for $g_p \in \GSp_4(\Q_p)$.
Both $t_{\sigma,p}\,\frak{W}_{\psi_p^\pm}(W_{p,\pm}^\circ,\varphi_p^\circ)$ and $\frak{W}_{\psi_p^\pm}(t_{\sigma,p}W_{p,\pm}^\circ,\varphi_p^\circ)$ are $\GSp_4(\Z_p)$-invariant Whittaker functions of ${}^\sigma\!\itPi_p$ with respect to $\psi_{U,p}^\pm$. By evaluating these Whittaker functions at $1$, we see that they must be equal. Thus we have 
\begin{align}\label{E:arithmeticity proof 3}
\begin{split}
& t_{\sigma,p}\,\frak{W}_{\psi_p^\pm}(W_{p,\pm}^\circ,\varphi_p^\circ)(g_p\cdot{\rm diag}(1,1,a,a))\\
& = {}^\sigma\!\omega_{\itPi_p}(a)^2|a|_p\cdot\frak{W}_{\psi_p^{\pm a}}(t_{\sigma,p}W_{p,\pm a}^\circ,\varphi_p^\circ)({\rm diag}(a^{-3},a^{-2},1,a^{-1})g_p)
\end{split}
\end{align}
for $g_p \in \GSp_4(\Q_p)$.
Let $p \in S$. 
By the explicit formula for the Weil representation, we have
\[
{}^\sigma\!(\omega_{\xi_S}(g,h)\varphi_S') = \omega_{{}^\sigma\!\xi_S}(g,h){}^\sigma\!\varphi_S'
\]
for any non-trivial additive character $\xi_S$ of $\Q_S$ and $\varphi_S' \in S(V^2(\Q_S))$.
It follows that
\begin{align*}
{}^\sigma\frak{W}_{\psi_S^\pm}\left(W_S,\,\varphi_S\right) &= \frak{W}_{\psi_S^{\pm at^2}}\left({}^\sigma W_S,\,{}^\sigma\!\varphi_S\right).
\end{align*}
Therefore, by Lemma \ref{L:change characters}, we have
\begin{align}\label{E:arithmeticity proof 4}
\begin{split}
&t_{\sigma,S}\,\frak{W}_{\psi_S^\pm}\left(W_S,\,\varphi_S\right)(g_S)\\
& = {}^\sigma\frak{W}_{\psi_S^\pm}\left(W_S,\,\varphi_S\right)\left({\rm diag}(u_{\sigma,S}^{-3},u_{\sigma,S}^{-2},1,u_{\sigma,S}^{-1})g_S\right)\\
& = \frak{W}_{\psi_S^{\pm at^2}}\left({}^\sigma W_S,\,{}^\sigma\!\varphi_S\right)\left({\rm diag}(u_{\sigma,S}^{-3},u_{\sigma,S}^{-2},1,u_{\sigma,S}^{-1})g_S\right)\\
& = {}^\sigma\!\omega_{\itPi_S}(t)^{-3}|t|_S^{-2}\cdot\frak{W}_{\psi_S^{\pm a}}\left(\ell(({\bf a}(t^{-2}),{\bf a}(t^{-2}))){}^\sigma W_S,\,\omega_{\psi_S^\pm}\left(\bp t^{-1}{\bf 1}_2 & 0 \\ 0 & t {\bf 1}_2\ep,1\right){}^\sigma\!\varphi_S\right)\left({\rm diag}(a^{-3},a^{-2},1,a^{-1})g_S\right)\\
& = {}^\sigma\!\omega_{\itPi_S}(t)^{-3}|t|_S^{-2}\cdot\frak{W}_{\psi_S^{\pm a}}\left(\ell(({\bf a}(a),{\bf a}(a)))t_{\sigma,S}W_S,\,\omega_{\psi_S^\pm}\left(\bp t^{-1}{\bf 1}_2 & 0 \\ 0 & t {\bf 1}_2\ep,1\right){}^\sigma\!\varphi_S\right)\left({\rm diag}(a^{-3},a^{-2},1,a^{-1})g_S\right)
\end{split}
\end{align}
for $g_S \in \GSp_4(\Q_S)$.
Note that 
\[
|t|_S^{-2}\prod_{p \notin S}{}^\sigma\!\omega_{\itPi_p}(a)^2|a|_p = \omega_{\itPi_\infty}(a)^{-2}|a|_\infty^{-1}{}^\sigma\!\omega_{\itPi_S}(a)^{-2}
\]
by the product formula and the automorphy of ${}^\sigma\!\omega_\itPi = \omega_{\itPi_\infty}\cdot{}^\sigma\!\omega_{\itPi_f}$.
By (\ref{E:arithmeticity proof 3}) and (\ref{E:arithmeticity proof 4}), we conclude that  (\ref{E:arithmeticity proof 2}) is equal to 
\begin{align}\label{E:arithmeticity proof 5}
\begin{split}
&\sigma(\zeta(2))^{-2}\cdot \omega_{\itPi_\infty}(a)^{-2}|a|_\infty^{-1}\cdot {}^\sigma\!\omega_{\itPi_S}(a^2t^3)^{-1}\cdot W_{(\lambda_1,\lambda_2;\,{\sf w}),\psi_{U,\infty}^\pm}(g_\infty)\\
&\times\frak{W}_{\psi_S^{\pm a}}\left(\ell(({\bf a}(a),{\bf a}(a)))t_{\sigma,S}W_S,\,\omega_{\psi_S^\pm}\left(\bp t^{-1}{\bf 1}_2 & 0 \\ 0 & t {\bf 1}_2\ep,1\right){}^\sigma\!\varphi_S\right)\left({\rm diag}(a^{-3},a^{-2},1,a^{-1})g_S\right)\\
&\times\prod_{p \notin S}\frak{W}_{\psi_p^{\pm a}}(t_{\sigma,p}W_{p,\pm a}^\circ,\varphi_p^\circ)\left({\rm diag}(a^{-3},a^{-2},1,a^{-1})g_p\right)
\end{split}
\end{align}
for $g \in \GSp_4(\A)$.
Note that for $f \in \itPi$, we have
\[
W_{f,\psi_U^\pm}(g) = W_{f,\psi_U^{\pm a}}({\rm diag}(a^{-3},a^{-2},1,a^{-1})g).
\]
Since
\begin{align*}
W_{{}^\sigma\!f_{W_S},\psi^{\pm a}\boxtimes\psi^{\pm a}} = \ell(({\bf a}(a),{\bf a}(a)))W_\infty^\pm\cdot\ell(({\bf a}(a),{\bf a}(a)))t_{\sigma,S}W_S\cdot\prod_{p\notin S}t_{\sigma,p}W_{p,\pm a}^\circ,
\end{align*}
by Lemma \ref{L:Whittaker factorization} for the additive character $\psi^{\pm a}$ and (\ref{E:arithmeticity proof 1}) with $b=a$, we deduce that (\ref{E:arithmeticity proof 5}) is equal to the Whittaker function of the global theta lift on the right-hand side of (\ref{E:arithmeticity proof 0}) with respect to $\psi_U^\pm$.
This completes the proof.
\end{proof}

\subsection{Petersson norms of endoscopic lifts}


Let $\<\mbox{ },\mbox{ }\>_{\SO(V)}: (\tau_1\boxtimes\tau_2) \times (\tau_1^\vee\boxtimes\tau_2^\vee) \rightarrow \C$ and $\<\mbox{ },\mbox{ }\>_{\GSp_4}: \itPi \times \itPi^\vee \rightarrow \C$ be the Petersson bilinears pairing defined by
\begin{align*}
\<f,f'\>_{\SO(V)} &= \int_{\SO(V)(\Q)\backslash\SO(V)(\A)}f(h_1)f'(h_1)\,dh_1^{\rm Tam},\\
\<f_1,f_2\>_{\GSp_4} &= \int_{\A^\times\GSp_4(\Q)\backslash\GSp_4(\A)}f_1(g)f_2(g)\,dg^{\rm Tam}.
\end{align*}
Here $dh_1^{\rm Tam}$ and $dg^{\rm Tam}$ are the Tamagawa measures on $\SO(V)(\A)$ and $\A^\times\backslash\GSp_4(\A)$, respectively.
For each place $v$ of $\Q$, we fix a non-zero $\SO(V)(\Q_v)$-equivariant bilinear pairing
\[
\<\mbox{ },\mbox{ }\>_v: (\tau_{1,v}\boxtimes\tau_{2,v}) \times (\tau_{1,v}^\vee\boxtimes\tau_{2,v}^\vee)\longrightarrow \C. 
\]
We assume the pairings are chosen so that if $f = \bigotimes_v f_v \in \tau_1\boxtimes\tau_2$ and $f' = \bigotimes_v f_v' \in \tau_1^\vee\boxtimes\tau_2^\vee$, then $\<f_v,f_v'\>_v=1$ for almost all $v$ and
\[
\<f,f'\>_{\SO(V)} = \prod_v\<f_v,f_v'\>_v.
\]
Let $\mathcal{B}_v: \mathcal{S}(V^2(\Q_v)) \times \mathcal{S}(V^2(\Q_v)) \rightarrow \C$ be the bilinear pairing defined by
\[
\mathcal{B}_v(\varphi_v,\varphi_v') = \int_{V^2(\Q_v)}\varphi_v(x_v)\varphi_v'(x_v)\,dx_v,
\]
where $dx_v$ is defined by the Haar measure on $\Q_v$ in \S\,\ref{SS:measures}.
Note that $\mathcal{B}_v$ is equivariant under the Weil representation $\omega_{\psi_v}\otimes\omega_{\psi_v^{-1}}$ for any non-trivial additive character $\psi_v$ of $\Q_v$.
For $f_v \in \tau_{1,v}\boxtimes\tau_{2,v}$, $f_v'\in \tau_{1,v}^\vee\boxtimes\tau_{2,v}^\vee$, and $\varphi_v,\varphi_v' \in \mathcal{S}(V^2(\Q_v))$, we define the local zeta integral
\begin{align*}
Z_v(f_v,f_v',\varphi_v,\varphi_v') &= \frac{\zeta_v(2)\zeta_v(4)}{L(1,\tau_{1,v} \times \tau_{2,v})}\cdot\int_{\SO(V)(\Q_v)}\mathcal{B}_v(h_{1,v}\cdot \varphi_v,\varphi_v')\<(\tau_{1,v}\boxtimes\tau_{2,v})(h_{1,v})f_v,f_v'\>_v\,dh_{1,v}.
\end{align*}
Here $L(s,\tau_{1,v}\times\tau_{2,v})$ is the Rankin--Selberg $L$-function of $\tau_{1,v}\times\tau_{2,v}$, $dh_{1,v}$ is the Haar measure on $\SO(V)(\Q_v)$ defined as in \cite[\S\,7.1]{CI2019}, and $h_{1,v}\cdot\varphi_v(x) = \varphi_v(h_{1,v}^{-1}\cdot x)$.
Note that the integral converges absolutely (cf.\,\cite[Lemma 7.7]{GI2011}).
We recall the Rallis inner product formula in the following theorem. Let
\[
L(s,\tau_1 \times \tau_2) = \prod_v L(s,\tau_{1,v} \times \tau_{2,v})
\]
be the Rankin--Selberg $L$-function of $\tau_1 \times \tau_2$.

\begin{thm}[Rallis inner product formula]\label{T:RIF}
Let $\psi$ be a non-trivial additive character of $\Q\backslash\A$.
For $f \in \tau_1\boxtimes\tau_2$, $f'\in \tau_1^\vee\boxtimes\tau_2^\vee$, and $\varphi,\varphi' \in \mathcal{S}(V^2(\A))$ with 
$
f = \bigotimes_v f_v$, $f' = \bigotimes_v f_v'$, $\varphi = \bigotimes_v \varphi_v$, $\varphi' = \bigotimes_v \varphi_v'$,
we have
\begin{align*}
\<\theta_\psi(f,\varphi),\theta_{\psi^{-1}}(f',\varphi')\>_{\GSp_4}
= 2\zeta(2)^{-2}\cdot \frac{L(1,\tau_1 \times \tau_2)}{\zeta(2)\zeta(4)} \cdot \prod_v Z_v(f_v,f_v',\varphi_v,\varphi_v').
\end{align*}
\end{thm}

\begin{proof}
This is a special case of the Rallis inner product formula proved in \cite[Theorem 11.3]{GQT2014} (see also \cite[\S 7]{GI2011}). The factor $2\zeta(2)^{-2}$ is the ratio between the Tamagawa measure and the product measure $\prod_v dh_{1,v}$ on $\SO(V)(\A)$.
\end{proof}

For each rational prime $p$ and $\sigma \in {\rm Aut}(\C)$, we fix $\sigma$-linear isomorphisms
\[
T_{\sigma,p}: \tau_{1,p}\boxtimes\tau_{2,p}\longrightarrow {}^\sigma\!\tau_{1,p}\boxtimes{}^\sigma\!\tau_{2,p},\quad T_{\sigma,p}^\vee: \tau_{1,p}^\vee\boxtimes \tau_{2,p}^\vee\longrightarrow {}^\sigma\!\tau_{1,p}^\vee\boxtimes{}^\sigma\!\tau_{2,p}^\vee.
\]
Let $\<\mbox{ },\mbox{ }\>_{\sigma,p} : ({}^\sigma\!\tau_{1,p}\boxtimes{}^\sigma\!\tau_{2,p}) \times ({}^\sigma\!\tau_{1,p}^\vee\boxtimes{}^\sigma\!\tau_{2,p}^\vee)\longrightarrow \C$
be the bilinear pairing defined by 
\[
\<T_{\sigma,p}f_p,T_{\sigma,p}^\vee f_p'\>_{\sigma,p} = \sigma\<f_p,f_p'\>_p
\]
for $f_p \in \tau_{1,p}\boxtimes\tau_{2,p}$ and $f_p' \in \tau_{1,p}^\vee\boxtimes\tau_{2,p}^\vee$. It is easy to verify that $\<\mbox{ },\mbox{ }\>_{\sigma,p}$ is $\SO(V)(\Q_p)$-equivariant (cf.\,Lemma \ref{L:sigma matrix coeff.}).

\begin{lemma}\label{L:local zeta integral}
Let $f_p \in \tau_{1,p}\boxtimes\tau_{2,p}$, $f_p'\in \tau_{1,p}^\vee\boxtimes\tau_{2,p}^\vee$, and $\varphi_p,\varphi_p' \in S(V^2(\Q_p))$.

(1) If $p \nmid {\rm cond}(\itPi)$, $f_p$ and $f_p'$ are $\SO(V)(\Z_p)$-invariant, and $\varphi_p = \varphi_p' = \varphi_p^\circ$, then
\[
Z_p(f_p,f_p',\varphi_p,\varphi_p') = \<f_p,f_p'\>_p.
\]

(2) We have
\[
\sigma Z_p(f_p,f_p',\varphi_p,\varphi_p') = Z_p(T_{\sigma,p}f_p,T_{\sigma,p}^\vee f_p',{}^\sigma\!\varphi_p,{}^\sigma\!\varphi_p')
\]
for all $\sigma \in {\rm Aut}(\C)$.
\end{lemma}

\begin{proof}
The first assertion is a special case of \cite[Proposition 6.2]{LNM1254} (see also \cite[Proposition 3]{LR2005}).
Let $\sigma \in {\rm Aut}(\C)$.
An argument similarly to the proof of Lemma \ref{L:3.1} shows that
\[
\sigma L(1,\tau_{1,p}\times\tau_{2,p}) = L(1, {}^\sigma\!\tau_{1,p}\times{}^\sigma\!\tau_{2,p}).
\]
Note that the integral defining $Z_p$ is actually a doubling local zeta integral (cf.\,\cite[\S\,7]{GI2011}). Therefore, proceeding as in the proof of Lemma \ref{P:ab doubling} below, we have
\begin{align*}
&\sigma\left( \int_{\SO(V)(\Q_p)}\mathcal{B}_p(h_{1,p}\cdot \varphi_p,\varphi_p')\<(\tau_{1,p}\boxtimes\tau_{2,p})(h_{1,p})f_p,f_p'\>_p\,dh_{1,p}\right)\\
& = \int_{\SO(V)(\Q_p)}\sigma \mathcal{B}_p(h_{1,p}\cdot \varphi_p,\varphi_p')\sigma \<(\tau_{1,p}\boxtimes\tau_{2,p})(h_{1,p})f_p,f_p'\>_p\,dh_{1,p}.
\end{align*}
Since $\varphi_p$ and $\varphi_p'$ have compact support, $\mathcal{B}_p(h_{1,p}\cdot \varphi_p,\varphi_p')$ is a finite sum involving $\varphi_p$ and $\varphi_p'$. Thus we have
\[
\sigma \mathcal{B}_p(h_{1,p}\cdot \varphi_p,\varphi_p') = \mathcal{B}_p(h_{1,p}\cdot {}^\sigma\!\varphi_p,{}^\sigma\!\varphi_p'). 
\]
Also note that
\[
\sigma \<(\tau_{1,p}\boxtimes\tau_{2,p})(h_{1,p})f_p,f_p'\>_p = \<({}^\sigma\!\tau_{1,p}\boxtimes{}^\sigma\!\tau_{2,p})(h_{1,p})T_{\sigma,p}f_p,T_{\sigma,p}^\vee f_p'\>_{\sigma,p}
\]
by definition.
This completes the proof.
\end{proof}

\begin{thm}\label{T:endoscopic case}
Theorem \ref{T:main} holds for $\itPi = \theta(\tau_1\boxtimes\tau_2)$.
\end{thm}

\begin{proof}
For $p \nmid {\rm cond}(\itPi)$, let $f_p^\circ \in \tau_{1,p}\boxtimes\tau_{2,p}$ and $(f_p^\circ)^\vee \in \tau_{1,p}^\vee\boxtimes\tau_{2,p}^\vee$ be the $\SO(V)(\Z_p)$-invariant vectors defining the restricted tensor products $\bigotimes_v \tau_{1,v}\boxtimes\tau_{2,v}$ and $\bigotimes_v \tau_{1,v}^\vee\boxtimes\tau_{2,v}^\vee$, respectively.
Let 
\[
f_{\tau_1 \boxtimes \tau_2} = \bigotimes_v f_{\tau_1 \boxtimes \tau_2,v} \in \tau_1^+\otimes\tau_2^+,\quad f_{\tau_1^\vee \boxtimes \tau_2^\vee} = \bigotimes_v f_{\tau_1^\vee \boxtimes \tau_2^\vee,v} \in (\tau_1^\vee)^+\otimes(\tau_2^\vee)^+
\]
be the normalized newforms of $\tau_1 \boxtimes \tau_2$ and $\tau_1^\vee \boxtimes \tau_2^\vee$, respectively.
We assume $f_{\tau_{1,p}\boxtimes\tau_{2,p}} = f_p^\circ$ and $f_{\tau_{1,p}^\vee\boxtimes\tau_{2,p}^\vee} = (f_p^\circ)^\vee$ for $p \nmid {\rm cond}(\itPi)$.
The Petersson norm $\Vert f_{\tau_1\boxtimes\tau_2}\Vert$ of $f_{\tau_1\boxtimes\tau_2}$ is defined by
\[
\Vert f_{\tau_1\boxtimes\tau_2}\Vert = \<f_{\tau_1\boxtimes\tau_2},(\tau_{1,\infty}^\vee\boxtimes\tau_{2,\infty}^\vee)([{\bf a}(-1),{\bf a}(-1)])f_{\tau_1^\vee\boxtimes\tau_2^\vee}\>_{\SO(V)}.
\] 
Let $S$ be the set of prime divisors of ${\rm cond}(\itPi)$. Fix 
\begin{align*}
f &= f_{\tau_1\boxtimes\tau_2,\infty} \otimes f_S \otimes (\otimes_{p \notin S}f_p^\circ),\quad 
f' = (\tau_{1,\infty}^\vee\boxtimes\tau_{2,\infty}^\vee)([1,{\bf a}(-1)])f_{\tau_1^\vee\boxtimes\tau_2^\vee,\infty} \otimes f_S' \otimes (\otimes_{p \notin S}(f_p^\circ)^\vee),\\
\varphi &= \varphi_\infty^+ \otimes \varphi_S \otimes (\otimes_{p \notin S}\varphi_p^\circ),\quad\quad\quad\varphi' = \varphi_\infty^- \otimes \varphi_S' \otimes (\otimes_{p \notin S}\varphi_p^\circ)
\end{align*}
for some $f_S \in \tau_{1,S}\boxtimes\tau_{2,S}$, $f_S' \in \tau_{1,S}^\vee\boxtimes\tau_{2,S}^\vee$, and $\varphi_S,\varphi_S' \in S(V^2(\Q_S))$ such that
\begin{align}\label{E:endoscopic proof 0}
\<\theta_\psi(f,\varphi),\itPi_\infty^\vee({\rm diag}(-1,-1,1,1))\theta_{\psi^{-1}}(f',\varphi')\>_{\GSp_4} \neq 0.
\end{align}
It is clear that $f \in \tau_1^+ \boxtimes \tau_2^+$ and $f' \in \tau_1^+ \boxtimes \tau_2^-$. 
Note that we have the factorization of $L$-functions:
\[
L(s,\itPi,{\rm Ad}) = L(s,\tau_1 \times \tau_2)\cdot L(s,\tau_1,{\rm Ad})\cdot L(s ,\tau_2,{\rm Ad}),
\]
where $L(s,\tau_i,{\rm Ad})$ is the adjoint $L$-function of $\tau_i$ for $i=1,2$.
By Theorem \ref{T:RIF} and Lemma \ref{L:local zeta integral}-(1), we have
\begin{align}\label{E:endoscopic proof 1}
\begin{split}
&\frac{\zeta(2)^4\cdot\<\theta_\psi(f,\varphi),\itPi_\infty^\vee({\rm diag}(-1,-1,1,1))\theta_{\psi^{-1}}(f',\varphi')\>_{\GSp_4}}{\Vert f_\itPi \Vert} \\
&= 2\cdot\frac{L(1,\itPi,{\rm Ad})}{\zeta(2)\zeta(4)\cdot \Vert f_\itPi \Vert}\cdot \frac{\zeta(2)^2\cdot\Vert f_{\tau_1\boxtimes\tau_2}\Vert}{L(1,\tau_1,{\rm Ad})\cdot L(1,\tau_2,{\rm Ad})}\cdot C_\infty\cdot\frac{Z_S(f_S,f_S',\varphi_S,\varphi_S')}{\<f_{\tau_1 \boxtimes \tau_2,S},f_{\tau_1^\vee \boxtimes \tau_2^\vee,S}\>_S}.
\end{split}
\end{align}
Here
\[
C_\infty = \frac{Z_\infty\left(f_{\tau_1\boxtimes\tau_2,\infty},f_{\tau_1^\vee\boxtimes\tau_2^\vee,\infty},\varphi_\infty^+,\,\omega_{\psi_\infty^{-1}}({\rm diag}(-1,-1,1,1),[1,{\bf a}(-1)])\varphi_\infty^-\right)}{\<f_{\tau_1\boxtimes\tau_2,\infty},(\tau_{1,\infty}^\vee\boxtimes\tau_{2,\infty}^\vee)([{\bf a}(-1),{\bf a}(-1)])f_{\tau_1^\vee\boxtimes\tau_2^\vee,\infty}\>_\infty}.
\]
Let $\sigma \in {\rm Aut}(\C)$. It is easy to see that
\[
\sigma \left( \frac{\<f_1,\itPi_\infty^\vee({\rm diag}(-1,-1,1,1))f_2\>_{\GSp_4}}{\Vert f_\itPi \Vert}\right) =  \frac{\<{}^\sigma\!f_1,{}^\sigma\!\itPi_\infty^\vee({\rm diag}(-1,-1,1,1)){}^\sigma\!f_2\>_{\GSp_4}}{\Vert f_{{}^\sigma\!\itPi} \Vert}
\]
for all $f_1 \in \itPi_{\rm mot}$ and $f_2 \in \itPi^\vee_{\rm mot}$.
By the result of Sturm \cite{Sturm1989}, we have
\[
\sigma\left( \frac{L(1,\tau_1,{\rm Ad})\cdot L(1,\tau_2,{\rm Ad})}{\zeta(2)^2\cdot\Vert f_{\tau_1\boxtimes\tau_2}\Vert}\right) =  \frac{L(1,{}^\sigma\!\tau_1,{\rm Ad})\cdot L(1,{}^\sigma\!\tau_2,{\rm Ad})}{\zeta(2)^2\cdot\Vert f_{{}^\sigma\!\tau_1\boxtimes{}^\sigma\!\tau_2}\Vert}.
\]
By Lemma \ref{L:local zeta integral}, we have
\[
\sigma \left(\frac{Z_S(f_S,f_S',\varphi_S,\varphi_S')}{\<f_{\tau_1 \boxtimes \tau_2,S},f_{\tau_1^\vee \boxtimes \tau_2^\vee,S}\>_S}\right) = \frac{Z_S(T_{\sigma,S}f_S,T_{\sigma,S}^\vee f_S',{}^\sigma\!\varphi_S,{}^\sigma\!\varphi_S')}{\<T_{\sigma,S}f_{\tau_1 \boxtimes \tau_2,S},T_{\sigma,S}^\vee f_{\tau_1^\vee \boxtimes \tau_2^\vee,S}\>_{\sigma,S}}.
\]
By the Chinese remainder theorem, there exists $a \in \Q_{>0}^\times$ such that $u_{\sigma,S} = at^2$ for some $t \in \Q_S^\times$. 
Now we apply $\sigma$ to both sides of (\ref{E:endoscopic proof 1}). It then follows from Proposition \ref{P:arithmeticity} that
\begin{align}\label{E:endoscopic proof 2}
\begin{split}
&\frac{\zeta(2)^4\cdot\<\theta_{\psi^a}({}^\sigma\!f,\varphi_{\sigma,a}),{}^\sigma\!\itPi_\infty^\vee({\rm diag}(-1,-1,1,1))\theta_{\psi^{-a}}({}^\sigma\!f',\varphi_{\sigma,a}')\>_{\GSp_4}}{\Vert f_{{}^\sigma\!\itPi} \Vert} \\
&= 2\cdot\sigma\left(\frac{L(1,\itPi,{\rm Ad})}{\zeta(2)\zeta(4)\cdot \Vert f_\itPi \Vert}\right)\cdot \frac{\zeta(2)^2\cdot\Vert f_{{}^\sigma\!\tau_1\boxtimes{}^\sigma\!\tau_2}\Vert}{L(1,{}^\sigma\!\tau_1,{\rm Ad})\cdot L(1,{}^\sigma\!\tau_2,{\rm Ad})}\cdot \sigma C_\infty\cdot\frac{Z_S(T_{\sigma,S}f_S,T_{\sigma,S}^\vee f_S',{}^\sigma\!\varphi_S,{}^\sigma\!\varphi_S')}{\<T_{\sigma,S}f_{\tau_1 \boxtimes \tau_2,S},T_{\sigma,S}^\vee f_{\tau_1^\vee \boxtimes \tau_2^\vee,S}\>_{\sigma,S}}.
\end{split}
\end{align}
Here
\begin{align*}
\varphi_{\sigma,a} &= \omega_{\psi_\infty}\left(\bp \sqrt{a}\,{\bf 1}_2&0 \\0 & \sqrt{a}^{-1}\,{\bf 1}_2\ep,1\right)\varphi_\infty^+\otimes \omega_{\psi_S}\left(\bp t^{-1}{\bf 1}_2 & 0 \\ 0 & t{\bf 1}_2 \ep,1\right){}^\sigma\!\varphi_S \otimes(\otimes_{p \notin S}\varphi_p^\circ),\\
\varphi_{\sigma,a}' &= \omega_{\psi_{\infty}^{-1}}\left(\bp \sqrt{a}\,{\bf 1}_2&0 \\0 & \sqrt{a}^{-1}\,{\bf 1}_2\ep,1\right)\varphi_\infty^-\otimes \omega_{\psi_S^{-1}}\left(\bp t^{-1}{\bf 1}_2 & 0 \\ 0 & t{\bf 1}_2 \ep,1\right){}^\sigma\!\varphi_S' \otimes(\otimes_{p \notin S}\varphi_p^\circ).
\end{align*}
By the equivariance under the Weil representation, we have
\begin{align}\label{E:endoscopic proof 3}
\begin{split}
&\mathcal{B}_\infty\left(h_{1,\infty}\cdot\omega_{\psi_\infty}\left(\bp \sqrt{a}\,{\bf 1}_2&0 \\0 & \sqrt{a}^{-1}\,{\bf 1}_2\ep,1\right)\varphi_\infty^+,\,\omega_{\psi_{\infty}^{-1}}\left(\bp -\sqrt{a}\,{\bf 1}_2&0 \\0 & \sqrt{a}^{-1}\,{\bf 1}_2\ep,[1,{\bf a}(-1)]\right)\varphi_\infty^-\right)\\
&=  \mathcal{B}_\infty\left(h_{1,\infty}\cdot\varphi_\infty^+,\,\omega_{\psi_{\infty}^{-1}}\left({\rm diag}(-1,-1,1,1),[1,{\bf a}(-1)]\right)\varphi_\infty^-\right),\\
&\mathcal{B}_S\left(h_{1,S}\cdot \omega_{\psi_S}\left(\bp t^{-1}{\bf 1}_2 & 0 \\ 0 & t{\bf 1}_2 \ep,1\right){}^\sigma\!\varphi_S,\,\omega_{\psi_S^{-1}}\left(\bp t^{-1}{\bf 1}_2 & 0 \\ 0 & t{\bf 1}_2 \ep,1\right){}^\sigma\!\varphi_S'\right)\\
& = \mathcal{B}_S\left(h_{1,S}\cdot {}^\sigma\!\varphi_S,\,{}^\sigma\!\varphi_S'\right)
\end{split}
\end{align}
for all $h_{1,\infty} \in \SO(V)(\R)$ and $h_{1,S} \in \SO(V)(\Q_S)$.
We may assume the isomorphisms ${}^\sigma\!\tau_1 \boxtimes {}^\sigma\!\tau_2 \simeq \bigotimes_v {}^\sigma\!\tau_{1,v} \boxtimes {}^\sigma\!\tau_{2,v}$ and ${}^\sigma\!\tau_1^\vee \boxtimes {}^\sigma\!\tau_2^\vee \simeq \bigotimes_v {}^\sigma\!\tau_{1,v}^\vee \boxtimes {}^\sigma\!\tau_{2,v}^\vee$ are normalized so that
\begin{align*}
{}^\sigma\!f &= f_{\tau_1\boxtimes\tau_2,\infty} \otimes T_{\sigma,S}f_S \otimes (\otimes_{p \notin S}T_{\sigma,p}f_p^\circ),\\
{}^\sigma\!f' &= (\tau_{1,\infty}^\vee\boxtimes\tau_{2,\infty}^\vee)([1,{\bf a}(-1)])f_{\tau_1^\vee\boxtimes\tau_2^\vee,\infty} \otimes T_{\sigma,S}^\vee f_S' \otimes (\otimes_{p \notin S}T_{\sigma,p}^\vee(f_p^\circ)^\vee).
\end{align*}
Then, by Theorem \ref{T:RIF} again together with (\ref{E:endoscopic proof 3}), we see that the left-hand side of (\ref{E:endoscopic proof 2}) is equal to
\[
2\cdot\frac{L(1,{}^\sigma\!\itPi,{\rm Ad})}{\zeta(2)\zeta(4)\cdot \Vert f_{{}^\sigma\!\itPi} \Vert}\cdot \frac{\zeta(2)^2\cdot\Vert f_{{}^\sigma\!\tau_1\boxtimes{}^\sigma\!\tau_2}\Vert}{L(1,{}^\sigma\!\tau_1,{\rm Ad})\cdot L(1,{}^\sigma\!\tau_2,{\rm Ad})}\cdot  C_\infty\cdot\frac{Z_S(T_{\sigma,S}f_S,T_{\sigma,S}^\vee f_S',{}^\sigma\!\varphi_S,{}^\sigma\!\varphi_S')}{\<T_{\sigma,S}f_{\tau_1 \boxtimes \tau_2,S},T_{\sigma,S}^\vee f_{\tau_1^\vee \boxtimes \tau_2^\vee,S}\>_{\sigma,S}}.
\]
By our assumption (\ref{E:endoscopic proof 0}), we have $Z_S(f_S, f_S',\varphi_S,\varphi_S') \neq 0$.
We thus conclude that
\[
\sigma\left(\frac{L(1,\itPi,{\rm Ad})}{\zeta(2)\zeta(4)\cdot \Vert f_\itPi \Vert}\cdot C_\infty\right) = \frac{L(1,{}^\sigma\!\itPi,{\rm Ad})}{\zeta(2)\zeta(4)\cdot \Vert f_{{}^\sigma\!\itPi} \Vert}\cdot C_\infty.
\]
Finally, it was proved in \cite[Lemma 8.10]{CI2019} that $C_\infty \in \Q^\times$.
This completes the proof.
\end{proof}

\section{Local Zeta Integrals}\label{S:local zeta integrals}

In this section, we study the convergence and the Galois equivariant properties of the doubling local zeta integrals and the Rankin-Selberg local zeta integrals defined in (\ref{E:local zeta1}) and (\ref{E:local zeta2}), respectively. The main results are Propositions \ref{P:ab doubling} and \ref{P:ab2}.

Let $\F$ be a non-archimedean local field of characteristic zero. Let $\frak{o}$, $\varpi$, and $q$ be the maximal compact subring of $\F$, a generator of the maximal ideal of $\frak{o}$, and the cardinality of the residue field $\frak{o} / \varpi\frak{o}$, respectively. Let $|\mbox{ }|$ be the absolute value on $\F$ normalized so that $|\varpi| = q^{-1}$. Fix a non-trivial additive character $\psi$ of $\F$.

\subsection{Some representation theory of $\GSp_4$}

In this section, we recall some results for $\GSp_4$ on the unitarizability criterion of generic representations and the asymptotic behavior of Whittaker functions.

Let $Q_1$ and $Q_2$ be the standard Siegel parabolic subgroup and the standard Klingen parabolic subgroup of $\GSp_4$, respectively, defined by 
\begin{align*}
Q_1 &= \left\{ \bp * & * & *&* \\ *&*&*&*\\0&0&*&*\\0&0& *&*  \ep\in \GSp_4  \right\},\\
Q_2 &= \left\{ \bp * & * & *&* \\ 0&*&*&*\\0&0&*&0\\0&*& *&*  \ep\in \GSp_4  \right\}.
\end{align*}
We denote by $N_{Q_1}$ and $N_{Q_2}$ the unipotent radical of $Q_1$ and $Q_2$, respectively. The standatd Levi components of $Q_1$ and $Q_2$ are given by
\begin{align*}
M_{Q_1} &= \left.\left\{ m_1(A,\nu) = \bp A & 0 \\ 0 & \nu {}^t\!A^{-1} \ep  \mbox{ }\right\vert\mbox{ } A \in \GL_2,\, \nu\in\GL_1\right\},\\
M_{Q_2} &= \left.\left\{  m_2(t,g)=\bp t & 0 & 0&0 \\ 0&a&0&b\\0&0&\nu t^{-1}&0\\0&c&0&d  \ep\mbox{ }\right\vert\mbox{ } \begin{array}{l}
t \in \GL_1 \\
g=\bp a & b \\ c & d \ep  \in \GL_2\\
\nu=\det(g) \in \GL_1
\end{array}
\right\}.
\end{align*}
Let $M_{Q_i}^1 = M_{Q_i} \cap \Sp_4$ for $i=1,2$.
For a character $\chi$ of $\F^\times$ and an irreducible admissible representation $\tau$ of $\GL_2(\F)$, let 
$\tau \rtimes \chi$ be the normalized induced representation acting via the right translation $\rho$ on the space consisting of smooth functions $f : \GSp_4(\F)\rightarrow \mathcal{V}_\tau$ such that
\[
f(nm_1(A,\nu)g) = \delta_{Q_1}^{1/2}(m_1(A,\nu))\chi(\nu)\tau(A)f(g)
\]
for all $n \in N_{Q_1}(\F)$, $m_1(A,\nu) \in M_{Q_1}(\F)$, and $g \in \GSp_4(\F)$. Here $\mathcal{V}_\tau$ is the representation space of $\tau$ and $\delta_{Q_1}$ is the modulus character of $Q_1(\F)$ given by $\delta_{Q_1}(m_1(A,\nu)) = |\det(A)|^3|\nu|^{-3}$. Similarly, let 
$\chi \rtimes \tau$ be the normalized induced representation acting via the right translation $\rho$ on the space consisting of smooth functions $f : \GSp_4(\F)\rightarrow \mathcal{V}_\tau$ such that
\[
f(nm_2(t,g')g) = \delta_{Q_2}^{1/2}(m_2(t,g'))\chi(t)\tau(g')f(g)
\]
for all $n \in N_{Q_2}(\F)$, $m_2(t,g') \in M_{Q_2}(\F)$, and $g \in \GSp_4(\F)$. Here $\delta_{Q_2}$ is the modulus character of $Q_2(\F)$ given by $\delta_{Q_2}(m_2(t,g)) = |t|^4|\det(g)|^{-2}$. Note that the central characters of $\tau \rtimes \chi$ and $\chi \rtimes \tau$ are equal to $\omega_\tau \chi^2$ and $\omega_\tau\chi$, respectively. 
By the results of Sally and Tadi\'c \cite{ST1993}, any non-supercuspidal irreducible admissible generic  representation of $\GSp_4(\F)$ is the generic subrepresentation of an induced representation in one of the following types:
\begin{itemize}
\item[(I)] ${\rm Ind}_{B(\F)}^{\GL_2(\F)}(\chi_1 \boxtimes \chi_2) \rtimes \chi$ for some characters $\chi_1,\chi_2,\chi$ of $\F^\times$ such that $\chi_1 \neq |\mbox{ }|^{\pm 1}$, $\chi_2 \neq |\mbox{ }|^{\pm 1}$, and $\chi_1 \neq |\mbox{ }|^{\pm 1}\chi_2^{\pm 1}$;
\item[(IIa)] $({\rm St}\otimes \mu) \rtimes \chi$ for some characters $\mu,\chi$ of $\F^\times$ such that $\mu \neq |\mbox{ }|^{\pm 3/2}$ and $\mu^2 \neq |\mbox{ }|^{\pm 1}$;
\item[(IIIa)] $\chi \rtimes ({\rm St}\otimes \mu)$ for some characters $\mu,\chi$ of $\F^\times$ such that $\chi \neq {\bf 1}$ and $\chi \neq |\mbox{ }|^{\pm 2}$;
\item[(IVa)] $|\mbox{ }|^2 \rtimes ({\rm St}\otimes \mu)$ for some character $\mu$ of $\F^\times$;
\item[(Va)] $({\rm St} \otimes \mu) \rtimes \chi$ for some characters $\mu,\chi$ of $\F^\times$ such that $\mu^2 = |\mbox{ }|$ and $\mu \neq |\mbox{ }|^{1/2}$;
\item[(VIa)] ${\bf 1} \rtimes ({\rm St}\otimes \mu)$ for some character $\mu$ of $\F^\times$;
\item[(VII)] $\chi \rtimes \tau$ for some character $\chi$ of $\F^\times$ and irreducible supercuspidal representation $\tau$ of $\GL_2(\F)$ such that $\chi \neq {\bf 1}$ and either $\chi^2 \neq |\mbox{ }|^{\pm 2}$ or $\chi = |\mbox{ }|^{\pm 1}$ or $\tau \otimes \chi|\mbox{ }|^{\pm 1} \neq \tau$;
\item[(VIIIa)] ${\bf 1} \rtimes \tau$ for some irreducible supercuspidal representation $\tau$ of $\GL_2(\F)$;
\item[(IXa)] $\chi \rtimes \tau$ for some character $\chi$ of $\F^\times$ and irreducible supercuspidal representation $\tau$ of $\GL_2(\F)$ such that $\chi^2 = |\mbox{ }|^{2},$ $\chi \neq |\mbox{ }|$, and $\tau \otimes \chi|\mbox{ }|^{-1} = \tau$;
\item[(X)] $\tau \rtimes \chi$ for some character $\chi$ of $\F^\times$ and irreducible supercuspidal representation $\tau$ of $\GL_2(\F)$ such that $\omega_\tau \neq |\mbox{ }|^{\pm 1}$;
\item[(XIa)] $\tau \rtimes \chi$ for some character $\chi$ of $\F^\times$ and irreducible supercuspidal representation $\tau$ of $\GL_2(\F)$ such that $\omega_\tau = |\mbox{ }|$.
\end{itemize}
Here ${\rm St}$ denotes the Steinberg representation of $\GL_2(\F)$ and we follow \cite{RS2007} for the labelling of types.
\begin{lemma}\label{L:generic unitary}
Let $\itPi$ be a non-supercuspidal irreducible admissible generic  representation of $\GSp_4(\F)$. Assume $\itPi$ is unitary. If $\itPi$ is of one of the types (IIIa), (VIa), and (VIIIa), then the inducing data are unitary. In the remaining cases, the following conditions are satisfied:
\begin{itemize}
\item[(I)] $|e(\chi_1)|+|e(\chi_2)|<1$ and $\chi_1\chi_2\chi^2$ is unitary;
\item[(IIa)] $|e(\mu)|<\tfrac{1}{2}$ and $\mu\chi$ is unitary;
\item[(IVa)] $|\mbox{ }|\mu$ is unitary;
\item[(Va)] $|\mbox{ }|^{1/2}\chi$ is unitary;
\item[(VII)] $|e(\chi)| < 1$ and $\omega_\tau\chi$ is unitary;
\item[(IXa)] $|\mbox{ }|\omega_\tau$ is unitary;
\item[(X)] $|e(\omega_\tau)|<1$ and $\omega_\tau\chi^2$ is unitary;
\item[(XIa)] $|\mbox{ }|^{1/2}\chi$ is unitary.
\end{itemize}
Moreover, if $\itPi$ is of one of the types (IVa), (Va), (IXa), and (XIa), then $\itPi$ is a discrete series representation.
\end{lemma}

\begin{proof}
The conditions for unitarizability were proved in \cite[Theorem 4.4, Propositions 4.7 and 4.9]{ST1993}. The assertion for discrete series representations was proved in \cite[Theorem 4.1, Propositions 4.6 and 4.8]{ST1993}.
\end{proof}

\begin{lemma}\label{L:Whittaker asymptotic}
Let $\itPi$ be an irreducible admissible generic  representation of $\GSp_4(\F)$. There exist a finite set $\frak{X}_\itPi$ of characters of ${\bf T}(\F)$ and a positive integer $N_\itPi$ such that for any Whittaker function $W$ of $\itPi$, we have
\[
W(tk) = \delta_{{\bf B}}(t)^{1/2}\sum_{0\leq n_1 \leq N_\itPi}\sum_{0\leq n_2 \leq N_\itPi}\sum_{\eta \in \frak{X}_\itPi}\eta(t)(\log_q|a|)^{n_1}(\log_q|b|)^{n_2}\varphi_{n_1,n_2,\eta}(a,b,k)
\]
for some locally constant function $\varphi_{n_1,n_2,\eta}$ on $\F \times \F \times \GSp_4(\o)$ with compact support for $0 \leq n_1,n_2 \leq N_\itPi$ and $\eta \in \frak{X}_\itPi$. 
Here $t = {\rm diag}(ab,a,b^{-1},1) \in {\bf T}(\F)$.
Moreover, the set $\frak{X}_\itPi$ is given as follows:
for $\eta \in \frak{X}_\itPi$, 
\[
\eta({\rm diag}(ab,a,b^{-1},1)) = \eta_1(a)\eta_2(b)
\]
with $\eta_1 = |\mbox{ }|^{e(\omega_\itPi)/2}$, $\eta_2={\bf 1}$ if $\itPi$ is supercuspidal and
\begin{itemize}
\item[(I)] $(\eta_1,\eta_2) \in \{(\chi,\chi_1^{-1}),(\chi,\chi_2^{-1}),(\chi\chi_1,\chi_1),(\chi\chi_1,\chi_2^{-1}),(\chi\chi_2,\chi_1^{-1}),(\chi\chi_2,\chi_2),(\chi\chi_1\chi_2,\chi_1),(\chi\chi_1\chi_2,\chi_2)\}$;
\item[(IIa)] $(\eta_1,\eta_2) \in \{(\chi,|\mbox{ }|^{1/2}\mu^{-1}),(\mu^2\chi,|\mbox{ }|^{1/2}\mu),(|\mbox{ }|^{1/2}\mu\chi,|\mbox{ }|^{1/2}\mu^{\pm1})\}$;
\item[(IIIa)] $(\eta_1,\eta_2) \in \{(|\mbox{ }|^{1/2}\mu,\chi^{-1}),(|\mbox{ }|^{1/2}\mu,|\mbox{ }|),(|\mbox{ }|^{1/2}\mu\chi,\chi),(|\mbox{ }|^{1/2}\mu\chi,|\mbox{ }|)\}$;
\item[(IVa)] $\eta_1 = |\mbox{ }|^{5/2}\mu$, $\eta_2 = |\mbox{ }|^2$;
\item[(Va)] $\eta_1 \in \{|\mbox{ }|\chi,|\mbox{ }|^{1/2}\mu\chi\}$, $\eta_2 = |\mbox{ }|^{1/2}\mu$;
\item[(VIa)] $\eta_1 = |\mbox{ }|^{1/2}\mu$, $\eta_2 \in \{{\bf 1},|\mbox{ }|\}$;
\item[(VII)] $\eta_1 = |\mbox{ }|^{e(\omega_\itPi)/2}$, $\eta_2 = \chi^{\pm1}$;
\item[(VIIIa)] $\eta_1 = |\mbox{ }|^{e(\omega_\itPi)/2}$, $\eta_2={\bf 1}$;
\item[(IXa)] $\eta_1 = |\mbox{ }|^{e(\omega_\itPi)/2}$, $\eta_2 = \chi$;
\item[(X)] $\eta_1 \in \{\chi,\omega_\tau\chi\}$, $\eta_2 = {\bf 1}$;
\item[(XIa)] $\eta_1 = \omega_\tau\chi$, $\eta_2 = {\bf 1}$.
\end{itemize}
\end{lemma}

\begin{proof}
The assertion follows from the result of Lapid and Mao \cite[Theorem 3.1]{LM2009} for the special case $G=\GSp_4$ (see also \cite[p.\,155, Proposition 1.1.1]{Jiang1996}). The formula for $\frak{X}_\itPi$ is then a consequence of the explicit formula for the semisimplification of the normalized Jacquet module of $\itPi$ with respect to the parabolic subgroup of $\GSp_4$ determined by the cuspidal support of $\itPi$ (cf.\,\cite[Tables A.3 and A.4]{RS2007}). For types (I)-(VIa), we consider the Borel subgroup ${\bf B}$. For types (VII), (VIIIa), and (IXa) (resp.\,(X) and (XIa)), we consider the Klingen parabolic subgroup $Q_2$ (resp.\,Siegel parabolic subgroup $Q_1$).
\end{proof}

\begin{rmk}
If $\itPi$ is either supercuspidal or of one of the types (VII), (VIIIa), and (IXa), then $\varphi_{n_1,n_2,\eta}=0$ when $|a|$ is sufficiently small. Thus in this case, $\eta_1$ can be any character. We take $\eta_1 = |\mbox{ }|^{e(\omega_\itPi)/2}$ so that the estimation in the proof of Proposition \ref{P:ab2} is more uniform.
\end{rmk}

The following lemma is on the asymptotic behavior of matrix coefficients for $\GL_2$ and will be used in the proof of Proposition \ref{P:ab doubling}.

\begin{lemma}\label{L:matrix coeff. GL_2}
Let $\tau$ be an irreducible admissible generic  representation of $\GL_2(\F)$ and $\phi$ a matrix coefficient of $\tau$. There exist characters $\eta_1,\eta_2$ of $\F^\times$ and locally constant functions $\mathit{\Phi}_1,\mathit{\Phi}_2$ on $\F \times \GL_2(\o) \times \GL_2(\o)$ such that
\[
\phi(k_1{\bf a}(t)k_2) = \eta_1(t)|t|^{1/2}\mathit{\Phi}_1(t,k_1,k_2) + \eta_2(t)|t|^{1/2}\mathit{\Phi}_2(t,k_1,k_2)
\]
for $(t,k_1,k_2) \in (\o \smallsetminus \{0\}) \times \GL_2(\o) \times \GL_2(\o)$.
Moreover, if $\tau = {\rm Ind}_{B(\F)}^{\GL_2(\F)}(\chi_1\boxtimes \chi_2)$ for some characters $\chi_1,\chi_2$, then $\eta_1,\eta_2 \in \{\chi_1,\chi_2\}$.
If $\tau = {\rm St} \otimes \mu$ for some character $\mu$, then $\eta_1=\eta_2 = \mu|\mbox{ }|^{1/2}$.
If $\tau$ is supercuspidal, then $\eta_1=\eta_2 = |\mbox{ }|^{e(\omega_\tau)/2}$.
\end{lemma}

\begin{proof}
This is well-known. Indeed, the assertion follows from the explicit computation of the normalized Jacquet module of $\tau$ with respect to $B$ (cf.\,\cite[\S\,8.12]{GH2011}). 
\end{proof}

\begin{rmk}
If $\tau$ is supercuspidal, then $\mathit{\Phi}_1=\mathit{\Phi}_2=0$ when $|t|$ sufficiently small. Thus in this case, $\eta_1$ and $\eta_2$ can be any characters. 
We take $\eta_1=\eta_2 = |\mbox{ }|^{e(\omega_\tau)/2}$ so that the estimation in the proof of Proposition \ref{P:ab doubling} is more uniform.
\end{rmk}

\subsection{Intertwining operators}

Let $B_8$ be the standard Borel subgroup of $\Sp_8$ defined by
\[
B_8 = \left\{\begin{pmatrix}
   * & * & *  & *  & * & * & * & * \\
   0 & * & *  & *  & * & * & * & * \\
   0 & 0 & *  & *  & * & * & * & * \\
   0 & 0 & 0  & *  & * & * & * & * \\
   0 & 0 & 0  & 0  & * & 0 & 0 & 0 \\
   0 & 0 & 0  & 0  & * & * & 0 & 0 \\
   0 & 0 & 0  & 0  & * & * & * & 0 \\
   0 & 0 & 0  & 0  & * & * & * & *
 \end{pmatrix} \in \Sp_8\right\}.
\]
Denote by $N_8$ the unipotent of $B_8$ and $T_8 \subset B_8$ the maximal torus of $\Sp_8$ consisting of diagonal matrices. Let $W_8=N_{\Sp_8}(T_8)/T_8$ be the Weyl group of $T_8$ in $\Sp_8$. For $i=1,2,3,4$, let $\epsilon_i : T_8 \rightarrow \GL_1$ be the algebraic character defined by
\[
\epsilon_i({\rm diag}(t_1,t_2,t_3,t_4, t_1^{-1}, t_2^{-1}, t_3^{-1}, t_4^{-1})) = t_i.
\]
The set of positive roots for $(\Sp_8,T_8)$ is given by
\[
\{\epsilon_i \pm \epsilon_j, 2\epsilon_k \, \vert \, 1 \leq i<j \leq 4, \,k=1,2,3,4\}.
\]
For each positive root $\epsilon$, we normalize the associated embedding $\iota_\epsilon : \SL_2 \rightarrow \Sp_8$ as follows:
\begin{align}\label{E:unipotent element}
\begin{split}
\iota_{\epsilon_i + \epsilon_j} ({\bf n}(x)) &= \bp {\bf 1}_4 & x(E_{i,j}+E_{j,i}) \\ 0 & {\bf 1}_4\ep,\quad \iota_{\epsilon_i - \epsilon_j} ({\bf n}(x)) = \bp {\bf 1}_4+xE_{i,j} & 0 \\ 0 & {\bf 1}_4-xE_{j,i}\ep,\\
 \iota_{2\epsilon_i} ({\bf n}(x))  &= \bp {\bf 1}_4 & xE_{i,i} \\ 0 & {\bf 1}_4\ep.
\end{split}
\end{align}
Here $E_{i,j} \in {\rm M_{4,4}}$ is the matrix with $(i,j)$-entry equal to $1$ and zero otherwise.
Let $N_\epsilon \subset \Sp_8$ be the image of the unipotent radial of $B$ under $\iota_\epsilon$ and identify $N_\epsilon$ with $\mathbb{G}_a$ via $\iota_\epsilon$.
For a character $\chi$ of $T_8(\F)$, let $I(\chi)$ be the normalized induced representation acting via the right translation $\rho$ on the space consisting of smooth functions $f : \Sp_8(\F)\rightarrow \C$ such that
\[
f(ntg)=\delta_{B_8}^{1/2}(t)\chi(t)f(g)
\]
for $n \in N_8(\F)$, $t \in T_8(\F)$, and $g \in \Sp_8(\F)$. The modulus character $\delta_{B_8}$ of $B_8(\F)$ is given by 
\[
\delta_{B_8}({\rm diag}(t_1,\cdots,t_4,t_1^{-1},\cdots,t_4^{-1}))=|t_1|^8|t_2|^6|t_3|^4|t_4|^2.
\] 
For $w \in W_8$, we define the intertwining operator
\begin{align}\label{E:intertwining integral}
\begin{split}
&M_{w} : I(\chi)\longrightarrow I(\chi^{w}),\\
&M_{w}f(g) = \int_{N_w(\F)}f(wng)\,dn.
\end{split}
\end{align}
Here $\chi^w(t) = \chi(wtw^{-1})$ and
\[
N_w = \prod_{\scriptstyle{\epsilon>0} \atop \scriptstyle{w\epsilon<0}}N_\epsilon.
\]
The Haar measure $dn$ is normalized so that ${\rm vol}(N_w(\o),dn)=1$. The integral is absolutely convergent if $\chi$ belongs to some open subset and can be meromorphically continued to all $\chi$. If we write $w = w_1\cdots w_\ell$ into a reduced decomposition, then we have
\begin{align}\label{E:intertwining decomp.}
M_w = M_{w_\ell}\circ \cdots \circ M_{w_1}.
\end{align}

The following lemma is on the analytic and Galois equivariant properties of the intertwining integrals in the simplest case. For characters $\chi_1, \chi_2$ of $\F^\times$, let ${\rm ind}_{B(\F)}^{\GL_2(\F)}(\chi_1\boxtimes\chi_2)$ be the (non-normalized) induced representation on the space consisting of smooth functions $f : \GL_2(\F)\rightarrow \C$ such that 
\[
f({\bf n}(x){\bf a}(t_1){\bf d}(t_2)g) = \chi_1(t_1)\chi_2(t_2)f(g)
\]
for $x \in \F$, $t_1,t_2 \in \F^\times$, and $g \in \GL_2(\F)$.

\begin{lemma}\label{L:GL_2}
Let $\chi_1,\chi_2$ be characters of $\F^\times$ and $f \in {\rm ind}_{B(\F)}^{\GL_2(\F)}(\chi_1\boxtimes\chi_2)$.
The intertwining integral 
\[
\int_\F f({\bf w}{\bf n}(x)) \,dx 
\]
is absolutely convergent if $e(\chi_1\chi_2^{-1})>1$. 
We have
\[
\sigma\left( \int_\F f({\bf w}{\bf n}(x)) \,dx \right) = \int_\F \sigma(f({\bf w}{\bf n}(x))) \,dx 
\]
for all $\sigma \in {\rm Aut}(\C)$ when both sides are absolutely convergent.
\end{lemma}

\begin{proof}
Indeed, we have
\[
\int_{\F}f({\bf w}{\bf n}(x))\,dx = \int_{|x| \leq q^N}f({\bf w}{\bf n}(x))\,dx + \chi_1(-1)\chi_2(-1)f(1)\int_{|x| > q^{N}}\chi_1^{-1}\chi_2(x)\,dx.
\]
Here $N$ is sufficiently large so that 
\[
f\left( \bp 1&0\\x&1 \ep \right) = f(1)
\]
for all $|x|< q^{-N}$. The first integral is a finite sum and the second integral converges for $e(\chi_1\chi_2^{-1})>1$. For $\sigma \in {\rm Aut}(\C)$, we have
\[
\sigma\left(\int_{|x| \leq q^N}f({\bf w}{\bf n}(x))\,dx\right) = \int_{|x| \leq q^N}\sigma(f({\bf w}{\bf n}(x)))\,dx
\]
and 
\[
\sigma\left(\int_{|x| > q^{N}}\chi_1^{-1}\chi_2(x)\,dx\right) = \int_{|x| > q^{N}}{}^\sigma\!\chi_1^{-1}{}^\sigma\!\chi_2(x)\,dx
\]
for $e(\chi_1\chi_2^{-1})>1$ and $e({}^\sigma\!\chi_1{}^\sigma\!\chi_2^{-1})>1$.
The assertions then follow at once.
\end{proof}

Let $\chi$ be a character of $T_8(\F)$. For $\sigma \in {\rm Aut}(\C)$, we have the $\sigma$-linear isomorphism 
\[
I(\chi) \longrightarrow I({}^\sigma\!\chi),\quad f \longmapsto {}^\sigma\!f
\]
with ${}^\sigma\!f(g) = \sigma(f(g))$ for $g \in \Sp_8(\F)$.

\begin{corollary}\label{C:Galois intertwining}
Let $\chi$ be a character of $T_8(\F)$ and $w \in W_8$. 
Then we have
\[
\sigma(M_wf(g)) = M_w{}^\sigma\!f(g)
\]
for all $\sigma \in {\rm Aut}(\C)$, $f \in I(\chi)$, and $g \in \Sp_8(\F)$ when both sides are absolutely convergent.
\end{corollary}

\begin{proof}
The assertion is an immediate consequence of (\ref{E:intertwining decomp.}) and Lemma \ref{L:GL_2}.
\end{proof}

\subsection{Doubling local zeta integrals}

Let $H_1 = \GL_4$ and $H_2 = \GL_2\times \Sp_4$. We define embeddings
\begin{align*}
\iota_1:H_1 &\longrightarrow \Sp_8,\quad g \longmapsto \bp g & 0 \\ 0 & {}^t\!g^{-1}\ep,\\
\iota_2:H_2 &\longrightarrow \Sp_8,\quad \left(\bp a&b \\c & d\ep,\bp A & B \\ C & D \ep \right) \longmapsto  
\bp
a&0&b&0&0&0&0&0\\
0&a_{11}&0&a_{12}&0&b_{11}&0&b_{12}\\
c&0&d&0&0&0&0&0\\
0&a_{21}&0&a_{22}&0&b_{21}&0&b_{22}\\
0&0&0&0&a'&0&b'&0\\
0&c_{11}&0&c_{12}&0&d_{11}&0&d_{12}\\
0&0&0&0&c'&0&d'&0\\
0&c_{21}&0&c_{22}&0&d_{21}&0&d_{22}\\
\ep.
\end{align*}
Here $\bp a' & b' \\ c' & d'\ep = \bp a & c \\ b & d\ep^{-1}$, $A=(a_{ij})$, $B=(b_{ij})$, $C=(c_{ij})$, and $D=(d_{ij})$. Recall we have identify $\Sp_4 \times \Sp_4$ as a subgroup of $\Sp_8$ vie the embedding (\ref{E:doubling embedd.}). The image of $M_{Q_i}^1 \times M_{Q_i}^1 \subset \Sp_4 \times \Sp_4$ in $\Sp_8$ factors through the embedding $
\iota_i$, thus induces an embedding from $M_{Q_i}^1 \times M_{Q_i}^1$ into $H_i$. We identify $M_{Q_i}^1 \times M_{Q_i}^1$ as a subgroup of $H_i$ in this way. Let $R_i$ be the maximal parabolic subgroup of $H_i$ defined by
\begin{align*}
R_1 & =  \left\{ \bp * & * & *&* \\ *&*&*&*\\0&0&*&*\\0&0& *&*  \ep\in H_1  \right\},\\
R_2 & = \left\{ \left( \bp * & * \\ 0 & *\ep ,\bp * & * & *&* \\ *&*&*&*\\0&0&*&*\\0&0& *&*  \ep\right) \in H_2\right\}.
\end{align*}
Put
\begin{align*}
\xi_i = \begin{cases}
\bp 0 & {\bf 1}_2 \\ {\bf 1}_2 & -{\bf 1}_2\ep \in H_1(\F) & \mbox{ if $i=1$},\\ 
\left(\bp 0 & 1 \\ 1 & -1\ep ,\begin{pmatrix}
  0 & 0 & -\frac{1}{2} & \frac{1}{2}  \\
  \frac{1}{2}  & \frac{1}{2} & 0 & 0 \\
  1 & - 1 & 0 & 0 \\
  0 & 0 & 1 & 1
 \end{pmatrix}\right)\in H_2(\F) & \mbox{ if $i=2$}.
\end{cases}
\end{align*}
Note that we have $\xi_i (m,m) \xi_i^{-1}\in R_i(\F)$
for all $m \in M_{Q_i}(\F)$.

Let $F \in I(s)$ be a holomorphic section, where $I(s)$ is the degenerate principal series representation defined in \S\,\ref{SS:local zeta 1}. For $i=1,2$, we define the intertwining integral 
\[
{\Psi}_i(g,s;F) = \int_{N_{Q_i}(\F)}F(\delta(n,1)\iota_i(\xi_i^{-1}g),s)\,dn,
\]
where $g \in H_i(\F)$.
By \cite[Proposition 1]{LR2005}, the integral $\Psi_i(g,s;F)$ converges absolutely for ${\rm Re}(s)$ sufficiently large and admits meromorphic continuation to $s \in \C$. Moreover, $\Psi_i(\mbox{ };F) \in {\rm Ind}_{R_i(\F)}^{H_i(\F)}(\mu_{i,s})$, where $\mu_{i,s}$ is the character of the standard Levi component of $R_i(\F)$ given by
\begin{align*}
\mu_{1,s}\left(\bp A & 0 \\ 0 & D \ep\right) &= |\det(A)|^{s+3/2}|\det(D)|^{-s+3/2},\\
\mu_{2,s}\left(\bp a & 0 \\ 0 & d \ep, \bp A & 0 \\ 0 & {}^t\!A^{-1} \ep\right) &= |a|^{s+2}|d|^{-s+2}|\det(A)|^{s}.
\end{align*}
In the following lemma, we explicitly determine the region of convergence and prove the Galois equivariant property of the intertwining integrals.

\begin{lemma}\label{L:4.3 2}
Let ${F} \in I(s)$ be a holomorphic section.
For $i=1,2$, the integral $\Psi_i(g,s;F)$ converges absolutely for ${\rm Re}(s)>-\tfrac{1}{2}$ and satisfies the Galois equivariant property
\[
\sigma\Psi_i(g,\tfrac{n}{2};F) = \Psi_i(g,\tfrac{n}{2};{}^\sigma\!F)
\]
for all $\sigma \in {\rm Aut}(\C)$, $g \in H_i(\F)$, and odd positive integers $n$.
\end{lemma}

\begin{proof}

Let $w_1,w_2 \in W_8$ be Weyl elements defined by 
\[
w_1 = \iota_{2\epsilon_4}({\bf w})\iota_{\epsilon_3-\epsilon_4}({\bf w})\iota_{2\epsilon_4}({\bf w}),\quad w_2 = \iota_{2\epsilon_4}({\bf w})\iota_{\epsilon_3-\epsilon_4}({\bf w})\iota_{\epsilon_2-\epsilon_3}({\bf w}).
\]
Here $\iota_\epsilon$ is the embedding defined in (\ref{E:unipotent element}) for each positive root $\epsilon$.
A direct calculation shows that
\begin{align*}
\delta (N_{Q_1}(\F),1) \delta^{-1} &\in P(\F)\cdot  \iota_1 \bp 0&0&0&1\\ 0&0&1&0 \\ 0&1&0&0 \\ 1&0&0&0\ep^{-1}w_1 N_{w_1}(\F) w_1^{-1}\iota_1 \bp 0&0&0&1\\ 0&0&1&0 \\ 0&1&0&0 \\ 1&0&0&0\ep,\\
\delta (N_{Q_2}(\F),1) \delta^{-1} &\in P(\F)\cdot \iota_1 \bp 0&0&1&0\\ 0&0&0&1 \\ 0&1&0&0 \\ 1&0&0&0\ep^{-1}w_2 N_{w_2}(\F) w_2^{-1}\iota_1 \bp 0&0&1&0\\ 0&0&0&1 \\ 0&1&0&0 \\ 1&0&0&0\ep.
\end{align*}
Note that $F \in I(s) \subset I(\chi_s)$, where $\chi_s$ is the character of $T_8(\F)$ defined by
\[
\chi_s({\rm diag}(t_1,\cdots,t_4,t_1^{-1},\cdots,t_4^{-1})) = |t_1|^{s-3/2}|t_2|^{s-1/2}|t_3|^{s+1/2}|t_4|^{s+3/2}.
\]
Therefore, we have
\begin{align*}
\Psi_1(g,s;F) &= M_{w_1}F \left( w_1^{-1}\iota_1 \bp 0&0&0&1\\ 0&0&1&0 \\ 0&1&0&0 \\ 1&0&0&0\ep \delta \iota_1 (\xi_1^{-1}g)\right),\\
\Psi_2(g,s;F) &= M_{w_2}F \left(w_2^{-1}\iota_1\bp 0&0&1&0\\ 0&0&0&1 \\ 0&1&0&0 \\ 1&0&0&0\ep \delta \iota_2 (\xi_2^{-1}g)\right)
\end{align*}
for $g \in H_i(\F)$.
Here $M_{w_i} : I(\chi_s) \rightarrow I(\chi_s^{w_i})$ is the intertwining operator defined in (\ref{E:intertwining integral}) for $i=1,2$. 
The absolute convergence for ${\rm Re}(s) > -\tfrac{1}{2}$ then follows immediately from Lemma \ref{L:GL_2}. Indeed, we have
\begin{align*}
\chi_s^{\iota_{2\epsilon_4}({\bf w})} ({\rm diag}(t_1,\cdots,t_4,t_1^{-1},\cdots,t_4^{-1})) &= |t_1|^{s-3/2}|t_2|^{s-1/2}|t_3|^{s+1/2}|t_4|^{-s-3/2},\\
\chi_s^{\iota_{2\epsilon_4}({\bf w})\iota_{\epsilon_3-\epsilon_4}({\bf w})} ({\rm diag}(t_1,\cdots,t_4,t_1^{-1},\cdots,t_4^{-1})) &= |t_1|^{s-3/2}|t_2|^{s-1/2}|t_3|^{-s-3/2}|t_4|^{s+1/2}.
\end{align*}
Hence $M_{w_1}F = M_{\iota_{2\epsilon_4}({\bf w})}\circ M_{\iota_{\epsilon_3-\epsilon_4}({\bf w})}\circ M_{\iota_{2\epsilon_4}({\bf w})}F$ and $M_{w_2}F = M_{\iota_{\epsilon_2-\epsilon_3}({\bf w})}\circ M_{\iota_{\epsilon_3-\epsilon_4}({\bf w})}\circ M_{\iota_{2\epsilon_4}({\bf w})}F$ are absolutely convergent for 
\[
{\rm Re}(s) > \max\{-\tfrac{3}{2},-1,-\tfrac{1}{2}\}=-\tfrac{1}{2}.
\]
For $\sigma \in {\rm Aut}(\C)$, we have ${}^\sigma\!\chi_{n/2} = \chi_{n/2}$ and ${}^\sigma\!F \vert_{s=n/2} \in I(\tfrac{n}{2}) \subset I(\chi_{n/2})$ for all odd integers $n$. The Galois equivariant property for $\Psi_i$ then follows from 
Corollary \ref{C:Galois intertwining}.
This completes the proof.

\end{proof}

Let $\itPi$ be an irreducible admissible representation of $\GSp_4(\F)$.

\begin{lemma}\label{L:sigma matrix coeff.}
Let $\phi$ be a matrix coefficient of $\itPi$. Then ${}^\sigma\!\phi$ is a matrix coefficient of ${}^\sigma\!\itPi$ for all $\sigma \in {\rm Aut}(\C)$.
\end{lemma}

\begin{proof}
Let $\sigma \in {\rm Aut}(\C)$. Let $\mathcal{V}_\itPi$ and $\mathcal{V}_{\itPi^\vee}$ be the representation spaces of $\itPi$ and $\itPi^\vee$, respectively, and fix $\sigma$-linear isomorphisms $t: \mathcal{V}_\itPi \rightarrow \mathcal{V}_\itPi$ and $t^\vee: \mathcal{V}_{\itPi^\vee} \rightarrow \mathcal{V}_{\itPi^\vee}$. The representations ${}^\sigma\!\itPi$ and ${}^\sigma\!\itPi^\vee$ are realized on $\mathcal{V}_\itPi$ and $\mathcal{V}_{\itPi^\vee}$, respectively, with actions defined in (\ref{E:sigma action}).
Let $\<\mbox{ },\mbox{ }\>$ be a non-zero equivariant bilinear pairing for $\itPi \times \itPi^\vee$ realized on $\mathcal{V}_\itPi \times \mathcal{V}_{\itPi^\vee}$.
Define a bilinear pairing $\<\mbox{ },\mbox{ }\>'$ on $\mathcal{V}_\itPi \times \mathcal{V}_{\itPi^\vee}$ by
\[
\<v,v^\vee\>' = \sigma\<t^{-1}v,(t^\vee)^{-1}v^\vee\>
\]
for $v \in \mathcal{V}_\itPi$ and $v^\vee \in \mathcal{V}_{\itPi^\vee}$. It is easy to verify that $\<\mbox{ },\mbox{ }\>'$ defines an equivariant pairing for  ${}^\sigma\!\itPi \times {}^\sigma\!\itPi^\vee$. Assume $\phi(g) = \<\itPi(g)v,v^\vee\>$ for some $v,v^\vee$. Then ${}^\sigma\!\phi(g) = \<tv,tv^\vee\>'$ is a matrix coefficient of ${}^\sigma\!\itPi$. This completes the proof.
\end{proof}

\begin{prop}\label{P:ab doubling}
Let $\phi$ be a matrix coefficient of $\itPi$ and $F \in I(s)$ be a holomorphic section. The local zeta integral $Z(s,\phi,F)$ converges absolutely for ${\rm Re}(s)$ sufficiently large, and satisfies the Galois equivariant property
\[
\sigma Z(\tfrac{n}{2},\phi,F) = Z(\tfrac{n}{2},{}^\sigma\!\phi,{}^\sigma\!F)
\]
for all $\sigma \in {\rm Aut}(\C)$ and sufficiently large odd integers $n$.
Assume that $\itPi$ is essentially unitary and generic, then $Z(s,\phi,F)$ converges absolutely for ${\rm Re}(s)\geq\tfrac{1}{2}$.
\end{prop}

\begin{proof}

If $\itPi$ is supercuspidal, then 
\[
Z(s,\phi,F) = \sum_{i=1}^n F(\delta(g_i,1),s)\phi(g_i)
\] 
for some $g_1,\cdots,g_n$ depending only on the support of $\phi$ and on a sufficiently small open compact subgroup of $\GSp_4(\F)$ which stabilizes $F$. The assertions then follow at once.

Suppose that $\itPi$ is a subrepresentation of an induced representation of the form $\tau \rtimes \chi$ or $\chi \rtimes \tau$ for some irreducible admissible generic  representation $\tau$ of $\GL_2(\F)$ and some character $\chi$ of $\F^\times$. 
Let $\eta_1$ and $\eta_2$ be the characters of $\F^\times$ depending on $\tau$ described in Lemma \ref{L:matrix coeff. GL_2}.
Let $Q=Q_1$ (resp,~$Q=Q_2$) in if $\itPi\subset \tau \rtimes \chi$ (resp.\,$\itPi \subset\chi \rtimes \tau$). 
We write
\[
\Psi =  \Psi_i,\quad \iota = \iota_i,\quad \xi=\xi_i, \quad \mu_{s}=\mu_{i,s}, \quad H=H_i,\quad R=R_i
\]
if $Q=Q_i$.
Fix a non-zero equivariant bilinear pairing $\<\mbox{ },\mbox{ }\>$ on $\tau \times \tau^\vee$.
We may assume that
\[
\phi(g) = \int_{\Sp_4(\o)}\<f(kg),f^\vee(k)\>\,dk
\]
for some $f$ and $f^\vee$.
For ${\rm Re}(s) > -\tfrac{1}{2}$, we have
\begin{align*}
Z(s,\phi,F) &= \int_{\Sp_4(\F)}F(\delta(g,1),s)\int_{\Sp_4(\o)}\<f(kg),f^\vee(k)\>\,dkdg\\
& = \int_{\Sp_4(\o)}\int_{\Sp_4(\F)} F(\delta(g,k),s)\<f(g),f^\vee(k)\>\,dgdk\\
& = \int_{\Sp_4(\o)^2}dk\int_{M_Q^1(\F)}dm\int_{N_Q(\F)}dn\, \delta_{Q}(m)^{-1}F(\delta(nmk_1,k_2),s)\<f(mk_1),f^\vee(k_2)\>\\
& = \int_{\Sp_4(\o)^2}dk\int_{M_Q^1(\F)}dm\, \delta_{Q}(m)^{-1}\Psi(\xi(m,1),s;\rho((k_1,k_2))F)\<f(mk_1),f^\vee(k_2)\>.
\end{align*}
Here we use (\ref{E:diagonal}) in the second line and Lemma \ref{L:4.3 2} in the fourth line.
Also note that we have identify $M_Q^1 \times M_Q^1$ as a subgroup of $H$ via $\iota$.
By the Cartan decomposition, for all $\varphi \in L^1(M_Q^1(\F))$, we have
\begin{align}\label{E:Q_1 measure}
\int_{M_Q^1(\F)}\varphi(m)\,dm = \int_{\F^\times}d^\times t_1 \int_{|t_2| \leq 1}d^\times t_2 \int_{\GL_2(\o)^2}dk\,{}^\sharp\!\left( \GL_2(\o){\bf a}(t_2)\GL_2(\o) / \GL_2(\o)\right) \varphi(k_1t_1{\bf a}(t_2)k_2) 
\end{align}
if $Q=Q_1$; and
\begin{align}\label{E:Q_2 measure}
\int_{M_Q^1(\F)}\varphi(m)\,dm = \int_{\F^\times}d^\times t_1\int_{|t_2| \leq 1}d^\times t_2 \int_{\SL_2(\o)^2}dk\, {}^\sharp\!\left( \SL_2(\o){\bf m}(t_2)\SL_2(\o) / \SL_2(\o)\right)\varphi(t_1,k_1{\bf m}(t_2)k_2)
\end{align}
if $Q=Q_2$. Here we identify $\GL_2$ and $\GL_1 \times \SL_2$ with $M_{Q_1}^1$ and $M_{Q_2}^1$ via the embeddings 
\[
A \longmapsto \bp A & 0 \\ 0 & {}^t\!A^{-1}\ep
,\quad
\left(t,\bp a & b \\ c & d\ep\right) \longmapsto \bp t & 0 & 0&0 \\ 0&a&0&b\\0&0& t^{-1}&0\\0&c&0&d  \ep.
\]
Note that for any $\Psi \in {\rm Ind}_{R(\F)}^{H(\F)}(\mu_s)$, we have
\begin{align}\label{E:Q_1 Iwasawa}
\begin{split}
\Psi(\xi(t_1{\bf a}(t_2),1)g,s) 
 &= |t_1^2t_2|^{s+5/2}\Psi\left(\bp 1&0&0&0 \\ 0&1&0&0 \\ t_1t_2&0&-1&0 \\ 0&t_1&0&-1 \ep g,s\right)\\
 &= |t_1^2t_2|^{-s+1/2} \Psi\left(\bp 0&0&1&0 \\ 0&0&0&1 \\ 1&0&-t_1^{-1}t_2^{-1}&0 \\ 0&1&0&-t_1^{-1} \ep g,s \right)\\
 &= |t_1|^3|t_2|^{s+5/2}\Psi\left(\bp 1&0&0&0 \\  0&0&0&1 \\ t_1t_2&0&-1&0 \\ 0&1&0&-t_1^{-1} \ep g,s \right)
\end{split}
\end{align}
for $t_1,t_2 \in \F^\times, g \in H(\F)$ if $Q=Q_1$; and 
\begin{align}\label{E:Q_2 Iwasawa}
\begin{split}
\Psi(\xi ( (t_1,{\bf m}(t_2))      ,1)g,s) &= 
|t_1|^{s+5/2}|t_2|^{s+3/2}\Psi\left(\left( \bp1 & 0 \\ t_1 & -1 \ep,  \bp 0&0&0&1 \\ 1&0&0&0 \\ t_2&-1&0&0 \\ 0&0&1&t_2\ep \right)g,s \right)\\
& = |t_1|^{-s+3/2}|t_2|^{s+3/2}\Psi\left(\left( \bp 0 & 1 \\ 1 & -t_1^{-1}\ep,  \bp 0&0&0&1 \\ 1&0&0&0 \\ t_2&-1&0&0 \\ 0&0&1&t_2\ep \right)g,s\right)
\end{split}
\end{align}
for $t_1,t_2 \in \F^\times, g \in H(\F)$ if $Q=Q_2$.
Assume $Q=Q_1$. By Lemma \ref{L:matrix coeff. GL_2}, there exist locally constant functions $\mathit{\Phi}_1,\mathit{\Phi}_2$ on $\F \times \GL_2(\o)\times \GL_2(\o) \times \Sp_4(\o) \times \Sp_4(\o)$ such that 
\[
\<f((k_1{\bf a}(t)k_2,1)k_1'),f^\vee(k_2')\> = \eta_1(t) |t|^2\mathit{\Phi}_1(t,k_1,k_2,k_1',k_2') + \eta_2(t) |t|^2\mathit{\Phi}_2(t,k_1,k_2,k_1',k_2')
\]
for $(t,k_1,k_2,k_1',k_2') \in (\o\smallsetminus \{0\}) \times \GL_2(\o) \times \GL_2(\o) \times \Sp_4(\o) \times \Sp_4(\o)$.
Combining with (\ref{E:Q_1 measure}) and (\ref{E:Q_1 Iwasawa}), for ${\rm Re}(s)>-\tfrac{1}{2}$, we see that the integral $Z(s,\phi,F)$ is a finite sum of integrals of the following forms: 
\begin{align*}
I_1(\tau,\eta,\varphi_1,\varphi_2,\varphi_3,s)&=\int_{(\F^\times)^2}|t_1|^{2s+2}\omega_\tau(t_1)|t_2|^{s+1/2}\eta(t_2)\varphi_1(t_1)\varphi_2(t_2)\varphi_3(t_1t_2)\,d(t_1,t_2),\\
I_2(\tau,\eta,\varphi_1,\varphi_2,\varphi_3,s)&=\int_{(\F^\times)^2}|t_1|^{2s+2}\omega_\tau(t_1)^{-1}|t_2|^{-s-3/2}\eta(t_2)\varphi_1(t_1)\varphi_2(t_2)\varphi_3(t_1t_2^{-1})\,d(t_1,t_2),\\
I_3(\tau,\eta,\varphi_1,\varphi_2,\varphi_3,s)&=\int_{(\F^\times)^2}\omega_\tau(t_1)^{-1}|t_2|^{s+1/2}\eta(t_2)\varphi_1(t_1)\varphi_2(t_2)\varphi_3(t_1^{-1}t_2)\,d(t_1,t_2),
\end{align*}
where $\eta\in \{\eta_1,\eta_2\}$ and $\varphi_i$ is a locally constant function on $\F$ with compact support for $i=1,2,3$. 
The integrals $I_1, I_2, I_3$ converge absolutely for
\begin{align}\label{E:Q_1 region 1}
{\rm Re}(s) > \max\{-1+\tfrac{1}{2}|e(\omega_\tau)|,-\tfrac{1}{2}-e(\eta),-\tfrac{1}{2} + e(\omega_\tau)-e(\eta)\}.
\end{align}
For an odd positive integer $n$ such that $s=\tfrac{n}{2}$ belongs to the above region of convergence, it is easy to verify that
\begin{align}\label{E:Q_1 Galois}
\begin{split}
\sigma I_1(\tau,\eta,\varphi_1,\varphi_2,\tfrac{n}{2}) & = I_1({}^\sigma\!\tau,{}^\sigma\!\eta,{}^\sigma\!\varphi_1,{}^\sigma\!\varphi_2,\tfrac{n}{2}),\\
\sigma I_2(\tau,\eta,\varphi_1,\varphi_2,\varphi_3,\tfrac{n}{2}) & = I_2({}^\sigma\!\tau,{}^\sigma\!\eta,{}^\sigma\!\varphi_1,{}^\sigma\!\varphi_2,{}^\sigma\!\varphi_3,\tfrac{n}{2}),  \\
\sigma I_3(\tau,\eta,\varphi_1,\varphi_2,\varphi_3,\tfrac{n}{2}) & =I_3({}^\sigma\!\tau,{}^\sigma\!\eta,{}^\sigma\!\varphi_1,{}^\sigma\!\varphi_2,{}^\sigma\!\varphi_3,\tfrac{n}{2})
\end{split}
\end{align}
for all $\sigma \in {\rm Aut}(\C)$.
Assume $Q=Q_2$. By Lemma \ref{L:matrix coeff. GL_2}, there exist locally constant functions $\mathit{\Phi}_1',\mathit{\Phi}_2'$ on $\F \times \SL_2(\o)\times \SL_2(\o) \times \Sp_4(\o) \times \Sp_4(\o)$ such that 
\[
\<f(((1,k_1{\bf m}(t)k_2),1)k_1'),f^\vee(k_2')\> = \omega_\tau^{-1}\eta_1^2(t) |t|\mathit{\Phi}_1'(t,k_1,k_2,k_1',k_2') + \omega_\tau^{-1}\eta_2^2(t) |t|\mathit{\Phi}_2'(t,k_1,k_2,k_1',k_2')
\]
for $(t,k_1,k_2,k_1',k_2') \in (\o\smallsetminus \{0\}) \times \SL_2(\o) \times \SL_2(\o) \times \Sp_4(\o) \times \Sp_4(\o)$.
Combining with (\ref{E:Q_2 measure}) and (\ref{E:Q_2 Iwasawa}),  for ${\rm Re}(s)>-\tfrac{1}{2}$, we see that the integral $Z(s,\phi,F)$ is a finite sum of integrals of the following forms: 
\begin{align*}
I_4(\chi,\eta,\varphi_1,\varphi_2,s) & = \int_{(\F^\times)^2}|t_1|^{s+1/2}\chi(t_1)|t_2|^{s+1/2}\eta(t_2)\varphi_1(t_1)\varphi_2(t_2)\,d(t_1,t_2),\\
I_5(\chi,\eta,\varphi_1,\varphi_2,s) & = \int_{(\F^\times)^2}|t_1|^{s+1/2}\chi(t_1)^{-1}|t_2|^{s+1/2}\eta(t_2)\varphi_1(t_1)\varphi_2(t_2)\,d(t_1,t_2),
\end{align*}
where $\eta\in \{\omega_\tau^{-1}\eta_1^2,\omega_\tau^{-1}\eta_2^2\}$ and $\varphi_i$ is a locally constant function on $\F$ with compact support for $i=1,2$.
The integrals $I_4, I_5$ converge absolutely for
\begin{align}\label{E:Q_2 region 1}
{\rm Re}(s) > \max\{-\tfrac{1}{2}+|e(\chi)|,-\tfrac{1}{2}+e(\omega_\tau)-2e(\eta_1),-\tfrac{1}{2}+e(\omega_\tau)-2e(\eta_2)\}.
\end{align}
For an odd positive integer $n$ such that $s=\tfrac{n}{2}$ belongs to the above region of convergence, it is easy to verify that
\begin{align}\label{E:Q_2 Galois}
\begin{split}
\sigma I_4(\chi,\eta,\varphi_1,\varphi_2,\tfrac{n}{2})& = I_4({}^\sigma\!\chi,{}^\sigma\!\eta,{}^\sigma\!\varphi_1,{}^\sigma\!\varphi_2,\tfrac{n}{2}),\\
\sigma I_5(\chi,\eta,\varphi_1,\varphi_2,\tfrac{n}{2}) & = I_5({}^\sigma\!\chi,{}^\sigma\!\eta,{}^\sigma\!\varphi_1,{}^\sigma\!\varphi_2,\tfrac{n}{2}) 
\end{split}
\end{align}
for all $\sigma \in {\rm Aut}(\C)$.
Let $\sigma \in {\rm Aut}(\C)$ and $n$ an odd positive integer such that $s=\tfrac{n}{2}$ belongs to the region of convergence.
We conclude from (\ref{E:Q_1 Galois}) and (\ref{E:Q_2 Galois}) that
\[
\sigma Z(\tfrac{n}{2},\phi,F) = \int_{\Sp_4(\o)^2}dk\int_{M_Q^1(\F)}dm\, \delta_{Q}(m)^{-1}\sigma \Psi(\xi(m,1),\tfrac{n}{2};\rho((k_1,k_2))F) \sigma \<f(mk_1),f^\vee(k_2)\>.
\]
Note that
\[
\sigma \Psi(\xi(m,1),\tfrac{n}{2};\rho((k_1,k_2))F) = \Psi(\xi(m,1),\tfrac{n}{2};\rho((k_1,k_2)){}^\sigma\!F)
\]
for $m \in M_Q^1, (k_1,k_2) \in \Sp_4(\o)^2$ by Lemma \ref{L:4.3 2} and
\[
{}^\sigma\!\phi(g) = \int_{\Sp_4(\o)}\sigma\<f(kg),f^\vee(k)\>.
\]
Therefore, we have
\[
\sigma Z(\tfrac{n}{2},\phi,F) = Z(\tfrac{n}{2},{}^\sigma\!\phi,{}^\sigma\!F).
\]

Assume $\itPi$ is non-supercuspidal, essentially unitary, and generic. Note that 
\[
(\tau\rtimes \chi)\otimes |\mbox{ }|^t = \tau\rtimes \chi|\mbox{ }|^t,\quad (\chi \rtimes \tau)\otimes|\mbox{ }|^t = \chi\rtimes (\tau\otimes|\mbox{ }|^t)
\]
for $t \in \C$.
By Lemma \ref{L:generic unitary}, the inequalities (\ref{E:Q_1 region 1}) and (\ref{E:Q_2 region 1}) are satisfied when ${\rm Re}(s)\geq \tfrac{1}{2}$ except when $\itPi$ is either of type (IVa) or (IXa) or (XIa). Indeed, if $\itPi$ is of one of the types (I), (II), (Va), and (X), then $Q=Q_1$ and (\ref{E:Q_1 region 1}) is satisfied for ${\rm Re}(s)\geq\tfrac{1}{2}$. 
If $\itPi$ is of one of the types (IIIa), (VIa), (VII), and (VIIIa), then $Q=Q_2$ and (\ref{E:Q_2 region 1}) is satisfied for ${\rm Re}(s)\geq \tfrac{1}{2}$. 
For the remaining types (IVa), (IXa), and (XIa), $\itPi$ is a discrete series representation. Hence the integral $Z(s,\phi,F)$ converges absolutely for ${\rm Re}(s) \geq -\tfrac{1}{2}$ by \cite[Lemma 9.5]{GI2014}.
This completes the proof.



\end{proof}

\subsection{Local zeta integrals for $\GSp_4 \times \GSp_4$}

We write
\begin{align*}
n^-(w,y,x,u,v) &= \iota_{\epsilon_3+\epsilon_4}\bp 1 & 0 \\ w & 1\ep \iota_{2\epsilon_3}\bp 1 & 0 \\ y & 1\ep \iota_{\epsilon_2+\epsilon_4}\bp 1 & 0 \\ x & 1\ep \iota_{\epsilon_1+\epsilon_4}\bp 1 & 0 \\ u & 1\ep \iota_{\epsilon_2+\epsilon_3}\bp 1 & 0 \\ v & 1\ep 
\end{align*}
and $n^-(w,y,x) = n^-(w,y,x,0,0)$. Here we follow (\ref{E:unipotent element}) for the notation.

Let $\mathcal{I}(s)$ be the degenerate principal series representation defined in \S\,\ref{SS:local zeta 2}.
For $\alpha, \beta \in \F$ and $\mathcal{F} \in \mathcal{I}(s)$, define the twisted intertwining integral
\[
I(\alpha,\beta;\mathcal{F},\psi) = \int_{\F^3}\mathcal{F}(n^-(w,y,x,0,0))\psi(-\alpha x + \beta y)\,dw\,dy\,dx.
\]

\begin{lemma}\label{L:4.1}
Let $u,v \in \F$.

(1) Suppose $u \in \F^\times$. We have
\[
I(\alpha,\beta;\rho(n^-(0,0,0,u,0))\mathcal{F},\psi) = |u|^{-s-1} I\left( u\alpha,\beta; \rho\left( \iota_{\epsilon_1+\epsilon_4}  \bp 0 & -1 \\ 1 & u^{-1} \ep \right)\mathcal{F} ,\psi \right).
\]

(2) Suppose $v \in \F^\times$. We have
\[
I(\alpha,\beta;\rho(n^-(0,0,0,0,v))\mathcal{F},\psi) = |v| ^{-2s-2}I\left(v\alpha,v^2\beta; \rho\left(\iota_{\epsilon_2+\epsilon_3}  \bp 0 & -1 \\ 1 & v^{-1}\ep \right)\mathcal{F},\psi \right).
\]

(3) Suppose $u,v \in \F^\times$. We have
\[
I(\alpha,\beta;\rho(n^-(0,0,0,u,v))\mathcal{F},\psi) = |u|^{-s-1}|v| ^{-2s-2} I\left( uv\alpha,v^2\beta; \rho\left( \iota_{\epsilon_1+\epsilon_4} \bp 0 & -1 \\ 1 & u^{-1} \ep \iota_{\epsilon_2+\epsilon_3}\bp 0 & -1 \\ 1 & v^{-1} \ep\right)\mathcal{F} ,\psi \right).
\]
\end{lemma}

\begin{proof}

Note that
\begin{align}
 \bp 1 & 0 \\ x & 1\ep =  {\bf m}(x^{-1}){\bf n}(x) \bp 0 & -1 \\ 1 & x^{-1}\ep, \label{E:4.1}\\
n^-(w,y,x) \iota_{\epsilon_1+\epsilon_4}({\bf n}(u)) n^-(w,y,x)^{-1} &\in \mathcal{P}(\F) \cap {\rm ker}(\delta_{\mathcal{P}}), \label{E:4.2}\\
n^-(w,y,x) \iota_{\epsilon_2+\epsilon_3}({\bf n}(v)) n^-(w-vxy,y,x)^{-1} &\in \mathcal{P}(\F) \cap {\rm ker}(\delta_{\mathcal{P}}), \label{E:4.3}.
\end{align}
One can easily verify that the first (resp.\,second) assertion follows from (\ref{E:4.1}) and (\ref{E:4.2}) (resp.\,(\ref{E:4.3})). The third assertion is a direct consequence of (1) and (2). We leave the detail to the readers.
\end{proof}

\begin{lemma}\label{L:4.2}
Let $\mathcal{F} \in \mathcal{I}(s)$ be a holomorphic section.
The integral $I(\alpha,\beta;\mathcal{F},\psi)$ converges absolutely for ${\rm Re}(s)>-1$ and satisfies the Galois equivariant property
\[
\sigma\left(I(\alpha,\beta;\mathcal{F},\psi)\vert_{s=n}\right) = I(\alpha,\beta;{}^\sigma\!\mathcal{F},{}^\sigma\!\psi)\vert_{s=n}
\]
for all $\sigma \in {\rm Aut}(\C)$ and integers $n \geq 0$. Moreover, for ${\rm Re}(s) > -1$, the integral $I(\alpha,\beta;\rho(n^-(0,0,0,u,v))\mathcal{F},\psi)$ as a function in $(\alpha,\beta,u,v) \in (\F^\times)^4$ is a finite sum of functions of the form
\begin{align*}
&\varphi_{1,1}(\alpha)\varphi_{1,2}(\beta)\varphi_{1,3}(u)\varphi_{1,4}(v) + \varphi_{2,1}(u\alpha)\varphi_{2,2}(\beta)\varphi_{2,3}(u^{-1})\varphi_{2,4}(v)\\
&+\varphi_{3,1}(v\alpha)\varphi_{3,2}(v^2\beta)\varphi_{3,3}(u)\varphi_{3,4}(v^{-1}) + \varphi_{4,1}(uv\alpha)\varphi_{4,2}(v^2\beta)\varphi_{4,3}(u^{-1})\varphi_{4,4}(v^{-1}),
\end{align*}
where $\varphi_{i,j}$ is a locally constant function on $\F^\times$ so that $\varphi_{i,j}(x)=0$ for $|x|$ sufficiently large and there there exist $c_{i,j} \in \C$ and character $\chi_{i,j}$ of $\F^\times$ such that
$\varphi_{i,j}(x) = c_{i,j}\cdot\chi_{i,j}(x)$
for $|x|$ sufficiently small.
\end{lemma}

\begin{proof}
We rewrite the integral $I(\alpha,\beta;\mathcal{F},\psi)$ into $8$ terms according to whether $|x|$, $|y|$, and $|w|$ are sufficiently large or not.
Note that
\begin{align*}
n^-(w,y,0)\iota_{\epsilon_2+\epsilon_4}({\bf n}(x))n^-(-w,-y,0) &\in \mathcal{P}(\F)\cap {\rm ker}(\delta_{\mathcal P}),\\
n^-(w,0,0)\iota_{2\epsilon_3}({\bf n}(y))n^-(-w,0,0) &\in \mathcal{P}(\F)\cap {\rm ker}(\delta_{\mathcal P}).
\end{align*}
It follows that for all sufficiently larger integers $N_1$, $N_2$, and $N_3$ depending only on $\mathcal{F}$, we have \begin{align*}
&I(\alpha,\beta;\mathcal{F},\psi)\\
&=\int_{|x|  \leq  q^{N_1}}dx\int_{|y|  \leq  q^{N_2}}dy\int_{|w|  \leq  q^{N_3}}dw \,\mathcal{F}(n^-(w,y,x))\psi(-\alpha x + \beta y)\\
&+\int_{|x|  >  q^{N_1}}|x| ^{-s-2}\psi(-\alpha x) \,dx\int_{|y|  \leq  q^{N_2}}dy\int_{|w|  \leq  q^{N_3}}dw \,\mathcal{F}\left(n^-(w,y,0) \iota_{\epsilon_2+\epsilon_4}({\bf w})\right)\psi(\beta y)\\
&+\int_{|y|  >  q^{N_2}}|y| ^{-s-2}\psi(\beta y)\,dy\int_{|x|  \leq  q^{N_1}}dx\int_{|w|  \leq  q^{N_3}}dw\, \mathcal{F}\left(n^-(w,0,x)\iota_{2\epsilon_3}({\bf w})\right)\psi(-\alpha x)\\
&+\int_{|x|  \leq  q^{N_1}}dx\int_{|y|  \leq  q^{N_2}}dy\int_{|w|  >  q^{N_3}}dw\,|w| ^{-s-3} \mathcal{F}\left(\iota_{\epsilon_3+\epsilon_4}\bp 0 & -1 \\ 1 & w^{-1} \ep n^-(0,y,x)\right)\psi(-\alpha x + \beta y)\\
&+\int_{|x|  >  q^{N_1}}|x| ^{-s-2}\psi(-\alpha x)\,dx\\
&\quad\times\int_{|y|  \leq  q^{N_2}}dy\int_{|w|  >  q^{N_3}}dw\,|w| ^{-s-3}\mathcal{F}\left(\iota_{\epsilon_3+\epsilon_4}\bp 0 & -1 \\ 1 & w^{-1} \ep  \iota_{2\epsilon_3}\bp  1&  \\ y & 1\ep\iota_{\epsilon_2+\epsilon_4}({\bf w}) \right)\psi(\beta y)\\
&+\int_{|y|  >  q^{N_2}}|y| ^{-s-2}\psi(\beta y)\,dy\\
&\quad \times \int_{|x|  \leq  q^{N_1}}dx\int_{|w|  >  q^{N_3}}dw\,|w| ^{-s-3}\mathcal{F}\left( \iota_{\epsilon_3+\epsilon_4}\bp 0 & -1 \\ 1 & w^{-1} \ep\iota_{\epsilon_2+\epsilon_4}\bp 1 & \\ x & 1\ep \iota_{2\epsilon_3}({\bf w})\right)\psi(-\alpha x)\\
&+\int_{|x|  >  q^{N_1}}|x| ^{-s-2}\psi(-\alpha x)dx\int_{|y|  >  q^{N_2}}|y| ^{-s-2}\psi(\beta y)\,dy\\
&\quad \times \int_{|w|  \leq  q^{N_3}}dw\,\mathcal{F}\left( \iota_{\epsilon_3+\epsilon_4}\bp  1&  \\ w & 1\ep\iota_{2\epsilon_3}({\bf w})\iota_{\epsilon_2+\epsilon_4}({\bf w}) \right)\\
&+\int_{|x|  >  q^{N_1}}|x| ^{-s-2}\psi(-\alpha x)\,dx\int_{|y|  >  q^{N_2}}|y| ^{-s-2}\psi(\beta y)\,dy\int_{|w|  >  q^{N_3}}|w| ^{-s-3}\,dw\\
&\quad \times \mathcal{F}\left( \iota_{\epsilon_3+\epsilon_4}({\bf w})\iota_{2\epsilon_3}({\bf w})\iota_{\epsilon_2+\epsilon_4}({\bf w})\right).
\end{align*}
We see that $I(\alpha,\beta;\mathcal{F},\psi)$ is absolutely convergent for ${\rm Re}(s)>-1$. Note that
\begin{align*}
\int_{|x| \leq q^{N}}\psi(ax)\,dx &= q^N\cdot\mathbb{I}_{\varpi^{N+d}\o}(a),\\
\int_{|x|\geq q^N}|x|^{-s-2}\psi(ax)\,dx &= \left[ \tfrac{1-q^{-1}}{1-q^{-s-1}}\cdot(q^{-N(s+1)}-q^{(d-1)(s+1)}|a|^{s+1})- q^{(d-1)(s+1)-1}|a|^{s+1}\right]\cdot\mathbb{I}_{\varpi^{N+d-1}\o}(a),
\end{align*}
where $\varpi^d\o$ is the largest fractional ideal of $\F$ on which $\psi$ is trivial.
Let ${\rm Re}(s)>-1$.
Combining with Lemma \ref{L:4.1}, we deduce that each term of $I(\alpha,\beta;\rho(n^-(0,0,0,u,v))\mathcal{F},\psi)$ as a function in $(\alpha,\beta,u,v) \in (\F^\times)^4$ is equal to a finite sum of functions satisfying the conditions in the lemma. The Galois equivariant property then also follows at once. Indeed, for an integer $n$, we have
$\sigma(|\mbox{ }|^n) = |\mbox{ }|^n$ and ${}^\sigma\!\mathcal{F}(g,n) = \sigma(\mathcal{F}(g,n))$
for all $g \in \Sp_4(\F)$ by definition. 
This completes the proof.
\end{proof}

\begin{lemma}\label{L:4.3}

Let $\varphi_1,\cdots,\varphi_7$ be locally constant functions on $\F^\times$ so that $\varphi_{i}(x)=0$ for $|x|$ sufficiently large and there there exist $c_{i} \in \C$, character $\chi_{i}$ of $\F^\times$, and integer $m_i \geq 0$ such that
$\varphi_{i}(x) = c_{i}\cdot\chi_{i}(x)(\log_q|x|)^{m_i}$
for $|x|$ sufficiently small.
Let $I_1(\varphi_1,\cdots,\varphi_7)$ be the integral defined by
\[
I_1(\varphi_1,\cdots,\varphi_7) = \int_{(\F^\times)^4} \varphi_1(a) \varphi_2(b) \varphi_3(a v^{-2}) \varphi_4(b u^{-1} v) \varphi_5(bu^{-1})\varphi_6(u)\varphi_7(v)\, d(a,b,u,v).
\]
Then we have
\[
\sigma I_1(\varphi_1,\cdots,\varphi_7) = I_1({}^\sigma\!\varphi_1,\cdots,{}^\sigma\!\varphi_7)
\]
for all $\sigma \in {\rm Aut}(\C)$ when both sides are absolutely convergent.
Similar assertion holds for the following integrlas:
\begin{align*}
I_2(\varphi_1,\cdots,\varphi_6) &= \int_{(\F^\times)^4} \varphi_1(a) \varphi_2(b) \varphi_3(a v^{-2}) \varphi_4(buv) \varphi_5(u)\varphi_6(v)\, d(a,b,u,v),\\
I_3(\varphi_1,\cdots,\varphi_6) &= \int_{(\F^\times)^4} \varphi_1(a) \varphi_2(b) \varphi_3(a v^{2}) \varphi_4(bu^{-1}v^{-1}) \varphi_5(u)\varphi_6(v)\, d(a,b,u,v),\\
I_4(\varphi_1,\cdots,\varphi_7) &= \int_{(\F^\times)^4} \varphi_1(a) \varphi_2(b) \varphi_3(a v^{2}) \varphi_4(buv^{-1}) \varphi_5(bv^{-1})\varphi_6(u)\varphi_7(v)\, d(a,b,u,v).
\end{align*}
\end{lemma}

\begin{proof}

We recall a type of local integral of the form
\begin{align}\label{E:integral 2}
\int_{\F^\times}
\chi(x) (\log_q|x|)^m\cdot \varphi(x)\,d^\times x,
\end{align}
where $\varphi$ is a locally constant function on $\F$ with compact support, $\chi$ is a character of $\F^\times$, and $m \geq 0$ is an integer. The integral converges absolutely when $e(\chi)>0$. In this case, it is easy to verify that the integral satisfies the Galois equivariant property
\[
\sigma\left(\int_{\F^\times}\chi(x) (\log_q|x|)^m\cdot\varphi(x)\,d^\times x\right) = \int_{\F^\times}{}^\sigma\!\chi(x) (\log_q|x|)^m\cdot {}^\sigma\!\varphi(x)\,d^\times x
\]
for all $\sigma \in {\rm Aut}(\C)$ (cf.\,\cite[Proposition A]{Grobner2018}) when both sides are absolutely convergent.

We only consider the integral $I_1(\varphi_1,\cdots,\varphi_7)$. The assertion for other three integrals $I_2,I_3,I_4$ can be proved in a similar way and we omit it.
First we consider the case when $\varphi_6$ and $\varphi_7$ vanish unless $|u|$ and $|v|$ are sufficiently small.
Then the integral $I_1(\varphi_1,\cdots,\varphi_7)$ is a finite sum of integrals of the form
\[
\int_{(\F^\times)^2} \varphi_1'(a) \varphi_3'(a v^{-2})\varphi_6'(v)\, d(a,v)\cdot  \int_{(\F^\times)^2}\varphi_2'(b) \varphi_4'(b u^{-1})  \varphi_5'(u)\,d(b,u),
\]
where $\varphi_i'$ satisfies the same conditions of the lemma and $\varphi_5',\varphi_6'$ vanish unless $|u|$ and $|v|$ are sufficiently small.
Firstly we consider the integral
\begin{align}\label{E:integral}
\int_{(\F^\times)^2} \varphi_1'(a) \varphi_3'(a v^{-2})\varphi_6'(v)\, d(a,v).
\end{align}
When the support of $\varphi_1'$ is contained in a bounded set away from zero, the above integral is a finite sum of integrals of the form (\ref{E:integral 2}). 
Similar assertion holds when the support of $\varphi_6'$ is contained in a bounded set away from zero.
Therefore, we may assume 
\[
\varphi_1'(a) = \chi(a)(\log_q|a|)^{n}\mathbb{I}_{\varpi^{N_1}\o}(a),\quad \varphi_6'(v) = \mu(v)(\log_q|v|)^{m}\mathbb{I}_{\varpi^{N_2}\o}(v)
\]
for some characters $\chi,\mu$ of $\F^\times$ and integers $n,m,N_1,N_2 \geq 0$. Furthermore, by the condition on $\varphi_3'$, we may assume that either 
\begin{align}\label{E:Schwartz 1}
\varphi_3' = \mathbb{I}_{c+\varpi^{N_3}\o}
\end{align}
for some $c \in \F$ and integer $N_3$ with $c \notin \varpi^{N_3}\o$ or 
\begin{align}\label{E:Schwartz 2}
\varphi_3'(x) = \nu(x)(\log_q|x|)^{r} \mathbb{I}_{\varpi^{N_3}\o}(x)
\end{align}
for some character $\nu$ of $\F^\times$ and integers $r,N_3 \geq 0$.  Write 
\[
\mathbb{I}_{\varpi^{N_2}\o} = \mathbb{I}_{\varpi^{N_4}\o} + (\mathbb{I}_{\varpi^{N_2}\o} - \mathbb{I}_{\varpi^{N_4}\o})
\]
for some sufficiently large $N_4$ such that $c \in \varpi^{-2N_4+N_1}\o$ (resp.~$\varpi^{N_3} \in \varpi^{-2N_4+N_1}\o$) if $\varphi_3'$ is a function of the form (\ref{E:Schwartz 1}) (resp.~(\ref{E:Schwartz 2})). 
Then the integral (\ref{E:integral}) over $v \in \varpi^{N_2}\o \smallsetminus \varpi^{N_4}\o$ is a finite sum of integrals of the form (\ref{E:integral 2}).
Suppose that $\varphi_3'$ is a function of the form (\ref{E:Schwartz 1}). Then the integral (\ref{E:integral}) over $v \in \varpi^{N_4}\o$ is equal to 
\begin{align*}
&\int_{|a| \leq q^{-N_1}}\int_{|v| \leq q^{-N_4}} \chi(a)(\log_q|a|)^{n}\mu(v)(\log_q|v|)^m \mathbb{I}_{c+\varpi^{N_3}\o}(av^{-2})\,d(a,v)\\
&=\int_{|a-c| \leq q^{-N_3}} \int_{\varpi^{N_4}\o}\chi(av^2)\mu(v)(\log_q|v|)^m(\log_q|cv^2|)^n\,d(a, v). 
\end{align*}
Now we assume that $\varphi_3'$ is a function of the form (\ref{E:Schwartz 2}). Then the integral (\ref{E:integral}) over $v \in \varpi^{N_4}\o$ is equal to 
\begin{align*}
&\int_{|a| \leq q^{-N_1}}\int_{|v| \leq q^{-N_4}} \chi(a)(\log_q|a|)^{n}\mu(v)(\log_q|v|)^m \mathbb{I}_{\varpi^{N_3}\o}(av^{-2})\,d(a,v)\\
&=\int_{|a| \leq q^{-N_3}|v|^2}\int_{|v| \leq q^{-N_4}} \chi(a)(\log_q|a|)^{n}\mu(v)(\log_q|v|)^m \,d(a,v)\\
&=\int_{|a| \leq q^{-N_3}}\int_{|v| \leq q^{-N_4}} \chi(av^2)(\log_q|av^2|)^{n}\mu(v)(\log_q|v|)^m \,d(a,v).
\end{align*}
In any case, we conclude that the integral (\ref{E:integral}) satisfies the Galois equivariant property, since it is a finite sum of products of integrals of the form (\ref{E:integral 2}). By a similar argument, the integral 
\[
\int_{(\F^\times)^2}\varphi_2'(b) \varphi_4'(b u^{-1})  \varphi_5'(u)\,d(b,u)
\]
also satisfies the Galois equivariant property.

Now we consider the remaining cases. If $\varphi_6$ vanishes when $|u|$ is sufficiently small, then the integral $I_1(\varphi_1,\cdots,\varphi_7)$ is a finite sum of integrals of the form
\begin{align}\label{E:type 1}
\int_{(\F^\times)^3} \varphi_1'(a) \varphi_2'(b) \varphi_3'(a v^{-2}) \varphi_4'(bv)\varphi_5'(v)\, d(a,b,v),
\end{align}
where $\varphi_i'$ satisfies the same conditions of the lemma. 
If $\varphi_7$ vanishes when $|v|$ is sufficiently small, then the integral $I_1(\varphi_1,\cdots,\varphi_7)$ is a finite sum of integrals of the form
\begin{align}\label{E:type 2}
\int_{(\F^\times)^3} \varphi_1'(a) \varphi_2'(b) \varphi_3'(bu^{-1}) \varphi_4'(u)\, d(a,b,u),
\end{align}
where $\varphi_i'$ satisfies the same conditions of the lemma.  Proceeding similarly as in the previous paragraph, one can prove that the integrals of types (\ref{E:type 1}) and (\ref{E:type 2}) also satisfy the Galois equivariant property. We leave the detail to the readers. This completes the proof.
\end{proof}

Let $\itPi$ be an irreducible admissible generic representation of $\GSp_4(\F)$.

\begin{prop}\label{P:ab2}
Let $W_1 \in \mathcal{W}(\itPi ,\psi_U)$, $W_2 \in \mathcal{W}(\itPi^\vee ,\psi_U^{-1})$, and $\mathcal{F} \in \mathcal{I}(s)$ be a holomorphic section. The local zeta integral
$\mathcal{Z}(s,W_1,W_2,\mathcal{F})$
converges absolutely for ${\rm Re}(s)$ sufficiently large, and satisfies the Galois equivariant property
\[
\sigma \mathcal{Z}(n,W_1,W_2,\mathcal{F}) = \mathcal{Z}(n,{}^\sigma W_1, {}^\sigma W_2,{}^\sigma\!\mathcal{F})
\]
for all $\sigma \in {\rm Aut}(\C)$ and sufficiently large odd integers $n$.
Assume $\itPi$ is essentially unitary, then $\mathcal{Z}(s,W_1,W_2,\mathcal{F})$ converges absolutely for ${\rm Re}(s) \geq 1$.
\end{prop}

\begin{proof}

By Lemma \ref{L:Whittaker asymptotic}, we may assume 
\begin{align*}
W_1(utk) &= \psi_U(u)\delta_{{\bf B}}(t)^{1/2}\eta(t)(\log_q|a|)^{n_1}(\log_q|b|)^{n_2}\varphi_1(a,b,k),\\
W_2(utk) &= {\psi}_U^{-1}(u)\delta_{{\bf B}}(t)^{1/2}\eta'(t)(\log_q|a|)^{n_1'}(\log_q|b|)^{n_2'}\varphi_2(a,b,k)
\end{align*}
for some $\eta \in \frak{X}_\itPi$ and $\eta' \in \frak{X}_{\itPi^\vee}$, some integers $0 \leq n_1,n_2,n_1',n_2' \leq N_\itPi$, and locally constant functions $\varphi_1,\varphi_2$ on $\F \times \F \times \GSp_4(\o)$ with compact support. Here $u \in U(\F)$, $t={\rm diag}(ab,a,b^{-1},1) \in {\bf T}(\F)$, and $k \in \GSp_4(\o)$. Write 
\[
\eta({\rm diag}(ab,a,b^{-1},1)) = \eta_1(a)\eta_2(b), \quad \eta'({\rm diag}(ab,a,b^{-1},1)) = \eta_1'(a)\eta_2'(b).
\]
In the notation as in the proof of \cite[Lemma 9.5]{CI2019}, we have
\[
k_a = e(\eta_1)-e(\eta_1')-e(\omega_\itPi),\quad k_b = e(\eta_2),\quad k_c = 2e(\eta_1')+e(\omega_\itPi),\quad k_d = e(\eta_2').
\]
Following the same argument as in the proof of that lemma, we see that the integral $\mathcal{Z}(s,W_1,W_2,\mathcal{F})$ converges absolutely for 
\begin{align*}
&{\rm Re}(s)>\\
&  \max\{ -e(\eta_1)+\tfrac{1}{2}e(\omega_\itPi)-1, -e(\eta_2)-1, -2e(\eta_1')-e(\omega_\itPi)-1, -e(\eta_2')-1, -e(\eta_1)-\tfrac{1}{2}e(\eta_2')+\tfrac{1}{2}e(\omega_\itPi)-1, \\
&\quad\quad\quad\quad\quad\quad\quad\quad\quad\quad\quad\quad\quad\quad  -2e(\eta_1)-2e(\eta_1')-1, -e(\eta_2)-e(\eta_2')-1,-2e(\eta_1')+e(\eta_2')-e(\omega_\itPi)-1\}.
\end{align*}
When $\itPi$ is essentially unitary, by Lemmas \ref{L:generic unitary} and \ref{L:Whittaker asymptotic}, one can verify case by case that the above inequality holds for ${\rm Re}(s)\geq 1$.

Now we show that the integral $\calZ(s, W_1,W_2, \calF)$ satisfies the Galois equivariant property.
Write 
\[
W_i(tk) = \delta_{{\bf B}}(t)^{1/2}\mathit{\Phi}_i(a,b,k)
\]
for $i=1,2$, $t={\rm diag}(ab,a,b^{-1},1) \in {\bf T}(\F)$, and $k \in \GSp_4(\o)$.
Let 
\[(\GSp_4(\o)\times \GSp_4(\o))^\circ = \left\{ (k_1,k_2) \in \GSp_4(\o)\times \GSp_4(\o) \mbox{ }\vert\mbox{ } \nu(k_1)=\nu(k_2)\right\}.\]
We define $({\bf B}\times {\bf B})^\circ$ and $({\bf T} \times {\bf T})^\circ$ in a similar way.
We have
\begin{align*}
\calZ(s, W_1,W_2, \calF) &= 
 \int_{Z_{H}(\F) \tilde{U}(\F) \backslash {\bf G}(\F)}
 \calF(\eta g, s) (W_1\otimes W_2)(g) \, dg\\
 &=\int_{(\GSp_4(\o)\times \GSp_4(\o))^\circ}\int_{Z_{H}(\F)\backslash({\bf T} \times {\bf T})^\circ(\F)} \delta_{({\bf B}\times {\bf B})^\circ}(t)^{-1} (W_1\otimes W_2)(tk)\\
 &\quad\times \int_{U'(\F)\backslash U(\F)} \calF(\eta(u,1)tk)\psi_U(u) \,du\,dt\,dk.
\end{align*}
Note that
\[
\left\{ \left({\rm diag}(ab,a,b^{-1},1), {\rm diag}(cd,c,c^{-1}d^{-1}a,c^{-1}a)\right) \mbox{ }\vert\mbox{ }a,b,c,d \in \F^\times\right\}
\]
is a set of representatives for $Z_{H}(\F)\backslash({\bf T} \times {\bf T})^\circ(\F)$. 
Let
\[
t_1={\rm diag}(ab,a,b^{-1},1),\quad t_2={\rm diag}(cd,c,c^{-1}d^{-1}a,c^{-1}a).
\]
For 
\[
u = \bp 1 & x & *&w \\0&1&*&y\\0&0&1&0\\0&0&-x&1 \ep \in U(\F),
\]
a direct calculation gives that
\begin{align*}
\eta (u,1)(t_1,t_2) \eta^{-1} = p(u,t_1,t_2)n^-(-ac^{-2}d^{-1}w,-ac^{-2}y,ac^{-1}d^{-1}x,abc^{-1}d^{-1}-1,ac^{-1}-1)
\end{align*}
for some $p(u,t_1,t_2) \in \mathcal{P}(\F)$ with $\delta_{\mathcal{P}}(p(u,t_1,t_2)) = |a^{3/2}bc^{-1}| ^6$.
Therefore,
\begin{align*}
&\int_{U'(\F)\backslash U(\F)} \calF(\eta(u,1)tk)\psi_U(u) \,du \\
&= \delta_{({\bf B}\times {\bf B})^\circ}((t_1,t_2))^{1/2}|a^{3/2}bc^{-1}| ^{s+1}\\
&\times \int_{\F^3} \mathcal{F}(n^-(w,y,x,abc^{-1}d^{-1}-1,ac^{-1}-1)\eta k ) \psi(-a^{-1}cdx + a^{-1}c^2y) \,dw\,dy\,dx\\
& = \delta_{({\bf B}\times {\bf B})^\circ}((t_1,t_2))^{1/2}|a^{3/2}bc^{-1}| ^{s+1} \cdot I(a^{-1}cd,a^{-1}c^2 ; \rho(n^-(0,0,0, abc^{-1}d^{-1}-1,ac^{-1}-1)\eta k )\mathcal{F},\psi).
\end{align*}
We conclude that 
\begin{align*}
\calZ(s, W_1,W_2, \calF) &= \int_{(\GSp_4(\o)\times \GSp_4(\o))^\circ} dk\int_{(\F^\times)^4} d(a,b,c,d)\,|a^{3/2}bc^{-1}| ^{s+1} \omega_\itPi(a^{-1}c)\mathit{\Phi}_1(a,b,k_1)\mathit{\Phi}_2(a^{-1}c^2,d,k_2) \\
&\quad \times I(a^{-1}cd,a^{-1}c^2 ; \rho(n^-(0,0,0,abc^{-1}d^{-1}-1,ac^{-1}-1)\eta k )\mathcal{F},\psi)\\
& = \int_{(\GSp_4(\o)\times \GSp_4(\o))^\circ} dk\int_{(\F^\times)^4} d(a,b,u,v)\,|a^{1/2}bv| ^{s+1} \mathit{\Phi}_1(a,b,k_1)\omega_\itPi(v)^{-1}\mathit{\Phi}_2(av^{-2},bu^{-1}v,k_2) \\
&\quad \times I(bu^{-1},av^{-2} ; \rho(n^-(0,0,0,u-1,v-1)\eta k )\mathcal{F},\psi).
\end{align*}
Here 
$d(a,b,u,v)$ is the Haar measure on $(\F^\times)^4$ with ${\rm vol}((\o^\times)^4,d(a,b,u,v))=1$.
By Lemma \ref{L:4.2}, for $k=(k_1,k_2) \in (\GSp_4(\o)\times \GSp_4(\o))^\circ $ and ${\rm Re}(s)$ sufficiently large, the integral
\begin{align*}
&\int_{(\F^\times)^4}d(a,b,u,v)\,|a^{1/2}bv| ^{s+1} \omega_\itPi(v)^{-1}\mathit{\Phi}_1(a,b,k_1)\mathit{\Phi}_2(av^{-2},bu^{-1}v,k_2)\\
&\times I(bu^{-1},av^{-2} ; \rho(n^-(0,0,0,u-1,v-1)\eta k )\mathcal{F},\psi)
\end{align*}
is a finite sum of integrals of the forms $I_1, I_2,I_3,I_4$ in Lemma \ref{L:4.3}. Let $\sigma \in {\rm Aut}(\C)$ and $n$ an odd integer so that the above integrals are all absolutely converge, we have
\[
\sigma(|a^{1/2}bv|^{n+1}) = |a^{1/2}bv|^{n+1}
\]
and 
\[
\sigma \left(I(\alpha,\beta;\rho(n^-(0,0,0,u-1,v-1)\eta k )\mathcal{F},\psi)\vert_{s=n} \right) = I(\alpha,\beta;\rho(n^-(0,0,0,u-1,v-1)\eta k ){}^\sigma\!\mathcal{F},{}^\sigma\!\psi)\vert_{s=n}
\]
by Lemma \ref{L:4.2}. It then follows from the Galois equivariant property proved in Lemma \ref{L:4.3} that 
\begin{align*}
&\sigma \calZ(n, W_1,W_2, \calF) \\
& = \int_{(\GSp_4(\o)\times \GSp_4(\o))^\circ} dk\int_{(\F^\times)^4} d(a,b,u,v)\,|a^{1/2}bv| ^{n+1} \omega_{{}^\sigma\!\itPi}(v)^{-1}{}^\sigma\!\mathit{\Phi}_1(a,b,k_1){}^\sigma\!\mathit{\Phi}_2(av^{-2},bu^{-1}v,k_2) \\
&\quad \times I(bu^{-1},av^{-2} ; \rho(n^-(0,0,0,u-1,v-1)\eta k ){}^\sigma\!\mathcal{F},{}^\sigma\!\psi).
\end{align*}
Finally, note that
\[
{}^\sigma W_1(t_1k_1){}^\sigma W_2(t_2k_2) = \delta_{\mathbf{B}}(t_1)^{1/2}\delta_{\mathbf{B}}(t_2)^{1/2}\omega_{{}^\sigma\!\itPi}(a^{-1}c){}^\sigma\!\mathit{\Phi}_1(a,b,k_1){}^\sigma\!\mathit{\Phi}_2(a^{-1}c^2,d,k_2)
\]
for $t_1 = {\rm diag}(ab,a,b^{-1},1), t_2={\rm diag}(cd,c,c^{-1}d^{-1}a,c^{-1}a) \in {\bf T}(\F)$ and $k_1,k_2 \in \GSp_4(\o)$. Therefore, we have
\[
\sigma \calZ(n, W_1,W_2, \calF) = \calZ(n, {}^\sigma W_1,{}^\sigma W_2, {}^\sigma\!\calF).
\]
This completes the proof.
\end{proof}

\end{document}